\documentclass[11pt]{article}
\usepackage{amssymb}

\usepackage{amsmath}
\usepackage{theorem}
\usepackage{here}
\usepackage[dvipdfmx]{color}

\newtheorem{thm}{Theorem}[section]
\newtheorem{lem}{Lemma}[section]

\newtheorem{prp}{Proposition}[section]
\theorembodyfont{\rmfamily}
\newtheorem{algo}{Algorithm}[section]

\newtheorem{remark}{Remark}[section]
\makeatletter

\@addtoreset{equation}{section}
\makeatother

\usepackage{comment} %
\usepackage{bm}
\usepackage[pdftex]{graphicx}

\def\ta{{\tau}}
\def\xt{{\tilde x}}

\def\pd{\partial}

\def\Var{{\rm Var}}

\def\al{{\alpha}}
\def\be{{\beta}}
\def\ga{{\gamma}}
\def\de{{\delta}}
\def\ep{{\varepsilon}}
\def\la{{\lambda}}
\def\si{{\sigma}}

\def\th{{\theta}}

\def\ka{{\kappa}}

\def\bbe{{\text{\boldmath $\beta$}}}

\def\bla{{\text{\boldmath $\lambda$}}}

\def\bth{{\text{\boldmath $\theta$}}}

\def\bpsi{{\text{\boldmath $\psi$}}}

\def\bnu{{\text{\boldmath $\nu$}}}

\def\tah{{\hat \tau}}

\def\lah{{\hat \la}}

\def\alt{{\tilde \al}}
\def\bet{{\tilde \be}}

\def\lat{{\tilde \la}}

\def\bbeh{{\widehat \bbe}}

\def\blah{{\widehat \bla}}

\def\bbet{{\widetilde \bbe}}

\def\blat{{\widetilde \bla}}

\def\sis{{\small \bsi}}

\def\De{{\Delta}}

\def\Ga{{\Gamma}}

\def\La{{\Lambda}}

\def\bLa{{\text{\boldmath $\La$}}}

\def\bPsi{{\text{\boldmath $\Psi$}}}

\def\b{{\text{\boldmath $b$}}}

\def\e{{\text{\boldmath $e$}}}

\def\v{{\text{\boldmath $v$}}}

\def\x{{\text{\boldmath $x$}}}
\def\y{{\text{\boldmath $y$}}}

\def\A{{\text{\boldmath $A$}}}

\def\I{{\text{\boldmath $I$}}}

\def\P{{\text{\boldmath $P$}}}
\def\Q{{\text{\boldmath $Q$}}}

\def\X{{\text{\boldmath $X$}}}

\def\yh{{\hat y}}

\def\yt{{\tilde y}}

\def\st{{\tilde s}}

\def\tr{{\rm tr\,}}

\def\[{{\text{\boldmath $[$}}}
\def\]{{\text{\boldmath $]$}}}

\def\/{{\Bigr/\!\!}}

\def\1r{{\rm (1)}}
\def\2r{{\rm (2)}}
\def\3r{{\rm (3)}}
\def\4r{{\rm (4)}}
\def\5r{{\rm (5)}}

\def\non{{\nonumber}}
\DeclareMathOperator*{\bdiag}{{\bf Diag\,}}
\def\sis{{{\si ^2}}}
\def\tas{{{\ta ^2}}}

\def\latsu{\underline{\lat ^2}}
\def\latso{\overline{\lat ^2}}
\def\latsm{\lat _{*}^{2}}
\def\tats{\tilde{\ta } ^2}
\def\tahs{\tah ^2}
\def\yd{y^{\dag }}
\def\ytd{\yt ^{\dag }}

\begin{document}
\title{Geometric Ergodicity of Gibbs Algorithms for a Normal Model With a Global-Local Shrinkage Prior}
\author{
Yasuyuki Hamura\footnote{Graduate School of Economics, Kyoto University, 
Yoshida-Honmachi, Sakyo-ku, Kyoto, 606-8501, JAPAN. 
\newline{
E-Mail: yasu.stat@gmail.com}} \
}
\maketitle
\begin{abstract}
We consider Gibbs samplers for a normal linear regression model with a global-local shrinkage prior and show that they produce geometrically ergodic Markov chains. 
First, under the horseshoe local prior and a three-parameter beta global prior under some assumptions, we prove geometric ergodicity for a Gibbs algorithm in which it is relatively easy to update the global shrinkage parameter. 
Second, we consider a more general class of global-local shrinkage priors. 
Under milder conditions, geometric ergodicity is proved for two- and three-stage Gibbs samplers based on rejection sampling. 
We also construct a practical rejection sampling method in the horseshoe case. 
Finally, a simulation study is performed to compare proposed and existing methods.

\par\vspace{4mm}
\end{abstract}

\section{Introduction}
\label{sec:introduction}
Since the seminal works of %
\cite{cps2009, cps2010}, global-local shrinkage priors and, in particular, horseshoe priors have been used by many scholars for Bayesian modeling and inference. 
The theoretical properties of these models related to efficient and robust estimation, and variable selection have also been extensively studied; cf. %
\cite{ps2011} and %
\cite{bdpw2019}. 

Because of the nonconjugacy of horseshoe models, Markov chain Monte Carlo (MCMC) algorithms are actually used while applying them. 
However, the efficiency or geometric ergodicity of such algorithms has not been fully investigated %
(see \cite{jh2001} for the notion and implications of geometric ergodicity). 
For example, %
\cite{job2020} and %
\cite{bbj2021} considered MCMC algorighms related to horseshoe priors but assumed that the prior distribution of some shrinkage parameter is truncated, as pointed out by %
\cite{bkp2022}. 
\cite{bkp2022} comprehensively examined Gibbs samplers for horseshoe and related models, and succeeded in removing the assumption of a truncated support. 
However, the authors roughly assumed that a global shrinkage parameter has a finite $(p / 2)$-th negative moment, where $p$ is the number of regression coefficients. 
Also, Metropolis steps were performed in practice in updating the global shrinkage prameter. 

In the first half of this study, we consider the same horseshoe model as in Section 2.1 of %
\cite{bkp2022} and a different MCMC algorighm from them, in which it is relatively easy to update the global shrinkage parameter. 
We show that the Markov chain based on the algorithm is geometrically ergodic under some assumptions. 
In particular, we prove geometric ergodicity for some untruncated three-parameter beta prior for the global shrinkage parameter which does not have a finite $(p / 5)$-th negative moment. 

In the second half of this study, we consider different algorithms to remove the negative moment condition on the global shrinkage parameter. 
They are based on a reparameterized model where all parameters appear in the likelihood factor. 
Consequently, we obtain a geometric ergodicity result which is also applicable when we use different priors from the horse prior. 
Since the algorithms use rejection sampling, we discuss approaches to construnct efficient Accept-Reject algorithms. 

The remainder of the article is organized as follows. 
In Section \ref{subsec:model}, our horseshoe model and algorithm are presented. 
A geometric ergodicity result for the algorithm is given in Section \ref{subsec:ergodicity}. 
In Section \ref{sec:general}, we consider a more general class of global-local shrinkage priors and, in the horseshoe case, construct a practical sampling method. 
In Section \ref{sec:sim}, we perform a simulation study to compare proposed and existing methods.

\section{The Horseshoe Case}
\label{sec:horseshoe} 
\subsection{The model and an algorithm}
\label{subsec:model} 
Consider the horseshoe model of Section 2.1 of %
\cite{bkp2022} and suppose that we observe the following: 
\begin{align}
&y_i \sim {\rm{N}} ( y_i | {\x _i}^{\top } \bbe , \sis ) \non 
\end{align}
for $i = 1, \dots , n$ and use the horseshoe prior: 
\begin{align}
&\be _k \sim {\rm{N}} ( \be _k | 0, \sis \tas {\la _k}^{2} ) \text{,} \non \\
&\la _k \sim {2 \over \pi } {1 \over 1 + {\la _k}^2} \text{,} \non 
\end{align}
for each $k = 1, \dots , p$, where $\y = ( y_i )_{i = 1}^{n}$, $\bbe = ( \be _k )_{k = 1}^{p}$, $\X = ( \x _1 , \dots , \x _n )^{\top }$ and where: 
\begin{align}
&\tas \sim \pi _{\ta } ( \tas ) \text{,} \non \\
&\sis \sim {\rm{IG}} ( \sis | a' , b' ) \text{,} \non 
\end{align}
for $a' , b' > 0$. 
By making the change of variables $\la _{k}^{2} = {\la _k}^2$ and $\lat _{k}^{2} =\ta ^2 \la _{k}^{2}$,  we obtain: 
\begin{align}
&y_i \sim {\rm{N}} ( y_i | {\x _i}^{\top } \bbe , \sis ) \text{,} \quad i = 1, \dots , n \text{,} \non \\
&\be _k \sim {\rm{N}} ( \be _k | 0, \sis \lat _{k}^{2} ) \text{,} \quad k = 1, \dots , p \text{,} \non \\
&\lat _{k}^{2} \sim {1 \over \pi } {( \tas )^{1 / 2} \over ( \lat _{k}^{2} )^{1 / 2}} {1 \over \tas + \lat _{k}^{2}} \text{,} \quad k = 1, \dots , p \text{,} \non \\
&\tas \sim \pi _{\ta } ( \tas ) \text{,} \non \\
&\sis \sim {\rm{IG}} ( \sis | a' , b' ) \text{.} \non 
\end{align}
Note that as in %
\cite{bkp2022}, 
\begin{align}
&{( \tas )^{1 / 2} \over ( \lat _{k}^{2} )^{1 / 2}} {1 \over \tas + \lat _{k}^{2}} = {( \tas )^{1 / 2} \over ( \lat _{k}^{2} )^{3 / 2}} {1 \over 1 + \tas / \lat _{k}^{2}} = {( \tas )^{1 / 2} \over ( \lat _{k}^{2} )^{3 / 2}} \int_{0}^{\infty } {1 \over {\nu _k}^2} \exp \Big\{ - {1 \over \nu _k} \Big( 1 + {\tas \over \lat _{k}^{2}} \Big) \Big\} d{\nu _k} \text{.} \non 
\end{align}
Then, the posterior distribution is equivalent to the density of $\bbe $, $\sis $, ${\blat {}^2} = ( \lat _{k}^{2} )_{k = 1}^{p}$, $\tas $, and the additional latent variables $\bnu = ( \nu _k )_{k = 1}^{p}$ proportional to: 
\begin{align}
&\Big\{ \prod_{i = 1}^{n} {1 \over ( \sis )^{1 / 2}} \exp \Big( - {( y_i - {\x _i}^{\top } \bbe )^2 \over 2 \sis } \Big) \Big\} \Big\{ \prod_{k = 1}^{p} {1 \over ( \sis )^{1 / 2} ( \lat _{k}^{2} )^{1 / 2}} \exp \Big( - {{\be _k}^2 \over 2 \sis \lat _{k}^{2}} \Big) \Big\} \non \\
&\times \Big[ \prod_{k = 1}^{p} {( \tas )^{1 / 2} \over ( \lat _{k}^{2} )^{3 / 2}} {1 \over {\nu _k}^2} \exp \Big\{ - {1 \over \nu _k} \Big( 1 + {\tas \over \lat _{k}^{2}} \Big) \Big\} \Big] \pi _{\ta } ( \tas ) {1 \over ( \sis )^{1 + a'}} \exp \Big( - {b' \over \sis } \Big) \text{.} \label{eq:joint_posterior} 
\end{align}

We now construct a two-stage Gibbs sampler for the parameters distributed as equation (\ref{eq:joint_posterior}). 
First, the full conditional distribution of ${\blat {}^2}$ is proportional to: 
\begin{align}
&\prod_{k = 1}^{p} \Big[ {1 \over ( \lat _{k}^{2} )^2} \exp \Big\{ - {1 \over \lat _{k}^{2}} \Big( {{\be _k}^2 \over 2 \sis } + {\tas \over \nu _k} \Big) \Big\} \Big] \text{.} \non 
\end{align}
Next, $( \bbe , \sis )$ and $( \bnu , \tas )$ are conditionally independent given ${\blat {}^2}$. 
The full conditional of $( \bnu , \tas )$ is proportional to: 
\begin{align}
&\pi _{\ta } ( \tas ) \prod_{k = 1}^{p} \Big[ {( \tas )^{1 / 2} \over {\nu _k}^2} \exp \Big\{ - {1 \over \nu _k} \Big( 1 + {\tas \over \lat _{k}^{2}} \Big) \Big\} \Big] \non \\
&= \pi _{\ta } ( \tas ) \Big\{ \prod_{k = 1}^{p} {( \tas )^{1 / 2} \over 1 + \tas / \lat _{k}^{2}} \Big\} \prod_{k = 1}^{p} \Big[ \Big( 1 + {\tas \over \lat _{k}^{2}} \Big) {1 \over {\nu _k}^2} \exp \Big\{ - {1 \over \nu _k} \Big( 1 + {\tas \over \lat _{k}^{2}} \Big) \Big\} \Big] \text{.} \non 
\end{align}
Meanwhile, the full conditional of $( \bbe , \sis )$ is proportional to: 
\begin{align}
&{1 \over ( \sis )^{n / 2 + p / 2 + 1 + a'}} \exp \Big[ - {1 \over 2 \sis } \Big\{ \sum_{i = 1}^{n} ( y_i - {\x _i}^{\top } \bbe )^2 + \sum_{k = 1}^{p} {{\be _k}^2 \over \lat _{k}^{2}} + 2 b' \Big\} \Big] \non 
\end{align}
and is given by: 
\begin{align}
&{\rm{IG}} ( \sis | n / 2 + a' , \{ \| \y \| ^2 - \y ^{\top } \X ( \X ^{\top } \X + \bPsi )^{- 1} \X ^{\top } \y \} / 2 + b' ) \non \\
&\times {\rm{N}}_p ( \bbe | ( \X ^{\top } \X + \bPsi )^{- 1} \X ^{\top } \y , \sis ( \X ^{\top } \X + \bPsi )^{- 1} ) \text{,} \non 
\end{align}
where 
\begin{align}
\bPsi &= \bdiag_{k = 1}^{p} {1 \over \lat _{k}^{2}} \text{,} \non 
\end{align}
since 
\begin{align}
&\sum_{i = 1}^{n} ( y_i - {\x _i}^{\top } \bbe )^2 + \sum_{k = 1}^{p} {{\be _k}^2 \over \lat _{k}^{2}} \non \\
&= \{ \bbe - ( \X ^{\top } \X + \bPsi )^{- 1} \X ^{\top } \y \} ^{\top } ( \X ^{\top } \X + \bPsi ) \{ \bbe - ( \X ^{\top } \X + \bPsi )^{- 1} \X ^{\top } \y \} \non \\
&\quad + \| \y \| ^2 - \y ^{\top } \X ( \X ^{\top } \X + \bPsi )^{- 1} \X ^{\top } \y \text{.} \non 
\end{align}

\begin{algo}
\label{algo:lat} 
The parameters $\bbe $, $\sis $, ${\blat {}^2}$, $\tas $, and $\bnu $ are updated in the following way. 
\begin{itemize}
\item
Sample 
\begin{align}
\sis &\sim {\rm{IG}} ( \sis | n / 2 + a' , \{ \| \y \| ^2 - \y ^{\top } \X ( \X ^{\top } \X + \bPsi )^{- 1} \X ^{\top } \y \} / 2 + b' ) \non 
\end{align}
and sample 
\begin{align}
\bbe &\sim {\rm{N}}_p ( \bbe | ( \X ^{\top } \X + \bPsi )^{- 1} \X ^{\top } \y , \sis ( \X ^{\top } \X + \bPsi )^{- 1} ) \text{.} \non 
\end{align}
Sample 
\begin{align}
\tas &\sim p( \tas | {\blat {}^2} ) = f( \tas ; {\blat {}^2} ) / \int_{0}^{\infty } f(t; {\blat {}^2} ) dt \text{,} \non 
\end{align}
where 
\begin{align}
f(t; \v ) &= \pi _{\ta } (t) \prod_{k = 1}^{p} {t^{1 / 2} \over v_k + t} \non 
\end{align}
for $t \in (0, \infty )$ for $\v \in (0, \infty )^p$, and sample 
\begin{align}
\nu _k &\sim %
{\rm{IG}} ( \nu _k | 1, 1 + \tas / \lat _{k}^{2} ) \non 
\end{align}
for each $k = 1, \dots , p$. 
\item
Sample 
\begin{align}
&\lat _{k}^{2} \sim %
{\rm{IG}} ( \lat _{k}^{2} | 1, {\be _k}^2 / (2 \sis ) + \tas / \nu _k ) \non 
\end{align}
for each $k = 1, \dots , p$. 
\end{itemize}
\end{algo}

In many cases, the conditional distribution of $\tas $ can be sampled from using rejection sampling. 
For example, suppose that the following holds: 
\begin{align}
\pi _{\ta } ( \tas ) %
\propto {( \tas )^{a - 1} \over (c + \tas )^{a + b}} \non 
\end{align}
for some $a, b, c > 0$. 
Let $\th = \log \tas \in \mathbb{R}$. 
Then, we have: 
\begin{align}
p( \th | {\blat {}^2} ) &\propto {( e^{\th } )^a \over (c + e^{\th } )^{a + b}} \prod_{k = 1}^{p} {( e^{\th } )^{1 / 2} \over \lat _{k}^{2} + e^{\th }} \non 
\end{align}
and 
\begin{align}
{\pd \over \pd \th } \log p( \th | {\blat {}^2} ) &= a - (a + b) {e^{\th } \over c + e^{\th }} + {p \over 2} - \sum_{k = 1}^{p} {e^{\th } \over \lat _{k}^{2} + e^{\th }} \text{,} \non 
\end{align}
which is decreasing in $\th $. 
Therefore, %
\cite{gw1992} method is applicable.

\subsection{Geometric ergodicity}
\label{subsec:ergodicity} 
Here, we show that Algorithm \ref{algo:lat} is efficient under some assumptions. 
For simplicity, we only consider the case where we use a (scaled) three-parameter beta %
\cite[see][]{adc2011} prior for $\tas $. 
Specifically, we let $a, b, c > 0$ and suppose that: 
\begin{align}
&\pi _{\ta } ( \tas ) = {1 \over C(a, b, c)} {( \tas )^{a - 1} \over (c + \tas )^{a + b}} \non 
\end{align}
for all $\tas \in (0, \infty )$, where $C(a, b, c) = \int_{0}^{\infty } \{ t^{a - 1} / (c + t)^{a + b} \} dt$. 
The class of the three-parameter beta distributions includes the invariant and inverse gamma distributions as limiting cases. 
Additionally, the half-Cauchy distribution is obtained by setting $a = b = 1 / 2$ and $c = 1$. 
\cite{ps2011} argued that $\pi _{\ta }$ ``should have substantial mass near zero'', which corresponds to smaller values of $a$ and $c$. 
For the theory in this section, the specific form of $\pi _{\ta }$ is not as important as the negative moment condition on $\tas $ or how small $a$ can be. 

\begin{thm}
\label{thm:lat} 
\hfill
\begin{itemize}
\item[{\rm{(i)}}]
The Markov chain based on Algorithm \ref{algo:lat} is geometrically ergodic if $a > p / 2$. 
\item[{\rm{(ii)}}]
Assume that $a \le p / 2$; $a \notin \mathbb{N}$ if $p \in 2 \mathbb{N}$; and $a \notin \mathbb{N} - 1 / 2$ if $p \in 2 \mathbb{N} - 1$. 
Suppose that 
\begin{align}
a / p > \log 2 - 1 / 2 \text{.} \label{eq:condition_ii} 
\end{align}
Then, the Markov chain based on Algorithm \ref{algo:lat} is geometrically ergodic. 
\end{itemize}
\end{thm}

Since $\exp (7 / 10) > 1 + 7 / 10 + (7 / 10)^2 / 2 + (7 / 10)^3 / 6 %
> 2$, we have $\log 2 - 1 / 2 < 1 / 5$. 
Therefore, (\ref{eq:condition_ii}) holds if $a / p \ge 1 / 5$. 
In contrast, the assumption of Remark 2.1 of %
\cite{bkp2022}, who do not restrict attention to the three-parameter beta priors, is roughly that $\tas $ has a finite $(p / 2)$-th negative moment. 
Other authors assume that the prior distributions of $\tas $ and/or ${\blat {}^2}$ are truncated %
(see Section 1 of \cite{bkp2022}). 
The technical requirement that $a \notin \mathbb{N}$ or $a \notin \mathbb{N} - 1 / 2$ is very weak in practice. 
For example, we can set $a = 0.51$ instead of $a = 0.5$. 
Finally, as discussed in Remark \ref{rem:small_a} of the Supplementary Material, it may be difficult to significantly relax the restrictive assumption (\ref{eq:condition_ii}).

\section{The General Case}
\label{sec:general} 
\subsection{The model and algorithms} 
Here, we consider using reparameterization and rejection sampling to sample from the posterior distribution when the local prior is not necessarily horseshoe. 
In the horseshoe model of Section 2.1 of %
\cite{bkp2022}, 
\begin{align}
&y_i \sim {\rm{N}} ( y_i | {\x _i}^{\top } \bbe , \sis ) \text{,} \quad i = 1, \dots , n \text{,} \non \\
&\be _k \sim {\rm{N}} ( \be _k | 0, \sis \tas {\la _k}^{2} ) \text{,} \quad k = 1, \dots , p \text{,} \non \\
&\la _k \sim {2 \over \pi } {1 \over 1 + {\la _k}^2} \text{,} \quad k = 1, \dots , p \text{,} \non \\
&\tas \sim \pi _{\ta } ( \tas ) \text{,} \non \\
&\sis \sim {\rm{IG}} ( \sis | a' , b' ) \text{.} \non 
\end{align}
We now suppose that instead of being positive, $\bla = ( \la _k )_{k = 1}^{p}$ is distributed on $\mathbb{R} ^p$. 
Additionally, we replace the (half-)Cauchy prior with a general local density $\pi _{\la } ( \la )$, $\la \in \mathbb{R}$. 
Then, the posterior density is: 
\begin{align}
p( \bbe , \sis , \tas , \bla | \y ) &\propto \Big\{ \prod_{i = 1}^{n} {\rm{N}} ( y_i | {\x _i}^{\top } \bbe , \sis ) \Big\} \Big[ \prod_{k = 1}^{p} \{ {\rm{N}} ( \be _k | 0, \sis \tas {\la _k}^{2} ) \pi _{\la } ( \la _k ) \} \Big] {\rm{IG}} ( \sis | a' , b' ) \pi _{\ta } ( \tas ) \text{.} \non 
\end{align}

In the remainder of this section, we suppress the dependence on $\y $. 
By making the change of variables $\bbet = ( \bet _k )_{k = 1}^{p} = ( \be _k / \{ ( \tas )^{1 / 2} \la _k \} )_{k = 1}^{p}$, we have: 
\begin{align}
p( \bbet , \sis , \tas , \bla ) &\propto \Big\{ \prod_{i = 1}^{n} {\rm{N}} ( y_i | ( \tas )^{1 / 2} {\x _i}^{\top } ( \bla \circ \bbet ), \sis ) \Big\} \Big[ \prod_{k = 1}^{p} \{ {\rm{N}} ( \bet _k | 0, \sis ) \pi _{\la } ( \la _k ) \} \Big] {\rm{IG}} ( \sis | a' , b' ) \pi _{\ta } ( \tas ) \text{.} \label{eq:general_joint} 
\end{align}
The shrinkage parameters $\tas $ and $\bla $ only appear in the exponents of the likelihood and in the prior factors. 
This feature is not seen in the original parameterization and simplifies the analysis by making it possible to dispense with energy functions corresponding to negative moments. 
Such a transformation is discussed, for example, at the end of Section 15.5 of %
\cite{gcsdvr2013}. 

Because directly sampling from the multivariate distribution $p( \bla | \bbet , \sis , \tas )$ is prohibitively difficult when $p$ is larger (unless $\X ^{\top } \X $ is diagonal), we introduce latent variables such that $\la _1 , \dots , \la _p$ are conditionally independent in the augmented model. 
Note that 
\begin{align}
&\prod_{i = 1}^{n} {\rm{N}} ( y_i | {\x _i}^{\top } \bbe , \sis ) \non \\
&= {1 / | \X ^{\top } \X |^{1 / 2} \over (2 \pi )^{(n - p) / 2} } {\exp \{ - \y ^{\top } \Q _{\X } \y / (2 \sis ) \} \over ( \sis )^{(n - p) / 2}} {\rm{N}}_p ( \bbe | ( \X ^{\top } \X )^{- 1} \X ^{\top } \y , \sis ( \X ^{\top } \X )^{- 1} ) \non \\
&= {1 / | \X ^{\top } \X |^{1 / 2} \over (2 \pi )^{(n - p) / 2} } {\exp \{ - \y ^{\top } \Q _{\X } \y / (2 \sis ) \} \over ( \sis )^{(n - p) / 2}} \non \\
&\quad \times \int_{\mathbb{R} ^p} {\rm{N}}_p ( \bbe | \bth , \sis ( d \I ^{(p)} )) {\rm{N}}_p ( \bth | ( \X ^{\top } \X )^{- 1} \X ^{\top } \y , \sis \{ ( \X ^{\top } \X )^{- 1} - d \I ^{(p)} \} ) d\bth \non 
\end{align}
for all $\bbe \in \mathbb{R} ^p$, where $\P _{\X } = \X ( \X ^{\top } \X )^{- 1} \X ^{\top }$ and $\Q _{\X } = \I ^{(n)} - \P _{\X }$, and where $d$ is any positive number making $d \I ^{(p)}$ smaller than $( \X ^{\top } \X )^{- 1}$. 
Then, by introducing latent variable $\bth = ( \th _k )_{k = 1}^{p} \in \mathbb{R} ^p$, we have: 
\begin{align}
p( \bla | \bth , \bbet , \sis , \tas ) %
&\propto \prod_{k = 1}^{p} \Big( \exp \Big[ - {\{ ( \tas )^{1 / 2} \la _k \bet _k - \th _k \} ^2 \over 2 d \sis } \Big] \pi _{\la } ( \la _k ) \Big) \label{eq:conditional_bla} 
\end{align}
and 
\begin{align}
&p( \bth | \bbet , \sis , \tas , \bla ) \non \\
&\propto \exp \Big( - {1 \over 2 \sis } \Big[ {\| \bbeh - \bth \| ^2 \over d} + \{ \bth - ( \X ^{\top } \X )^{- 1} \X ^{\top } \y \} ^{\top } \{ ( \X ^{\top } \X )^{- 1} - d \I ^{(p)} \} ^{- 1} \{ \bth - ( \X ^{\top } \X )^{- 1} \X ^{\top } \y \} \Big] \Big) \non \\
&\propto {\rm{N}}_p ( \bth | [ \I ^{(p)} / d + \{ ( \X ^{\top } \X )^{- 1} - d \I ^{(p)} \} ^{- 1} ]^{- 1} [ \bbeh / d + \{ ( \X ^{\top } \X )^{- 1} - d \I ^{(p)} \} ^{- 1} ( \X ^{\top } \X )^{- 1} \X ^{\top } \y ], \non \\
&\quad \sis [ \I ^{(p)} / d + \{ ( \X ^{\top } \X )^{- 1} - d \I ^{(p)} \} ^{- 1} ]^{- 1} ) \text{,} \non 
\end{align}
where $\bbeh = ( \tas )^{1 / 2} ( \bla \circ \bbet )$. 
Meanwhile, by marginalizing $\bth $, we obtain: 
\begin{align}
&p( \bbet , \sis , \tas | \bla ) \non \\
&\propto {\rm{N}}_p ( \bbet | \{ \I ^{(p)} + \tas \bLa ( \X ^{\top } \X ) \bLa \} ^{- 1} \bLa \X ^{\top } \y ( \tas )^{1 / 2} , \sis \{ \I ^{(p)} + \tas \bLa ( \X ^{\top } \X ) \bLa \} ^{- 1} ) \non \\
&\quad \times {\rm{IG}} ( \sis | n / 2 + a' , \y ^{\top } [ \I ^{(n)} - \X \bLa \{ \I ^{(p)} / \tas + \bLa ( \X ^{\top } \X ) \bLa \} ^{- 1} \bLa \X ^{\top } ] \y / 2 + b' ) \non \\
&\quad \times {1 \over ( \y ^{\top } [ \I ^{(n)} - \X \bLa \{ \I ^{(p)} / \tas + \bLa ( \X ^{\top } \X ) \bLa \} ^{- 1} \bLa \X ^{\top } ] \y / 2 + b' )^{n / 2 + a'}} {\pi _{\ta } ( \tas ) \over | \I ^{(p)} + \tas \bLa ( \X ^{\top } \X ) \bLa |^{1 / 2}} \text{,} \non 
\end{align}
where $\bLa = \bdiag_{k = 1}^{p} \la _k$.

\begin{algo}
\label{algo:pc-two} 
The parameters $\bth $, $\bbet $, $\sis $, $\tas $, and $\bla $ are updated in the following way. 
\begin{itemize}
\item
Sample 
\begin{align}
\bla &\sim p( \bla | \bth , \bbet , \sis , \tas ) \non \\
&\propto \prod_{k = 1}^{p} \Big( \exp \Big[ - {\{ ( \tas )^{1 / 2} \la _k \bet _k - \th _k \} ^2 \over 2 d \sis } \Big] \pi _{\la } ( \la _k ) \Big) \text{.} \non 
\end{align}
\item
Sample 
\begin{align}
\tas &\sim p( \tas | \bla ) \non \\
&\propto {\pi _{\ta } ( \tas ) / | \I ^{(p)} + \tas \bLa ( \X ^{\top } \X ) \bLa |^{1 / 2} \over ( \y ^{\top } [ \I ^{(n)} - \X \bLa \{ \I ^{(p)} / \tas + \bLa ( \X ^{\top } \X ) \bLa \} ^{- 1} \bLa \X ^{\top } ] \y / 2 + b' )^{n / 2 + a'}} \text{,} \non 
\end{align}
sample 
\begin{align}
\sis &\sim p( \sis | \tas , \bla ) \non \\
&= {\rm{IG}} ( \sis | n / 2 + a' , \y ^{\top } [ \I ^{(n)} - \X \bLa \{ \I ^{(p)} / \tas + \bLa ( \X ^{\top } \X ) \bLa \} ^{- 1} \bLa \X ^{\top } ] \y / 2 + b' ) \text{,} \non 
\end{align}
sample 
\begin{align}
\bbet &\sim p( \bbet | \sis , \tas , \bla ) \non \\
&= {\rm{N}}_p ( \bbet | \{ \I ^{(p)} + \tas \bLa ( \X ^{\top } \X ) \bLa \} ^{- 1} \bLa \X ^{\top } \y ( \tas )^{1 / 2} , \sis \{ \I ^{(p)} + \tas \bLa ( \X ^{\top } \X ) \bLa \} ^{- 1} ) \text{,} \non 
\end{align}
and sample 
\begin{align}
\bth &\sim p( \bth | \bbet , \sis , \tas , \bla ) \non \\
&= {\rm{N}}_p ( \bth | \A ^{- 1} \{ ( \tas )^{1 / 2} ( \bla \circ \bbet ) / d + \b \} , \sis \A ^{- 1} ) \text{,} \non 
\end{align}
where 
\begin{align}
\A &= \I ^{(p)} / d + \{ ( \X ^{\top } \X )^{- 1} - d \I ^{(p)} \} ^{- 1} \quad \text{and} \quad \b = \{ ( \X ^{\top } \X )^{- 1} - d \I ^{(p)} \} ^{- 1} ( \X ^{\top } \X )^{- 1} \X ^{\top } \y \text{.} \non 
\end{align}
\end{itemize}
\end{algo}

As discussed on page 14 of %
\cite{bkp2022}, sampling $\tas \sim p( \tas | \bla )$ can be difficult in general. 
To deal with this difficulty, we consider using a partially collapsed %
(see \cite{vp2008}) three-stage Gibbs sampler. 
Note that: 
\begin{align}
p( \tas | \bla , \bbet , \sis ) %
&\propto {\rm{N}}_p (( \tas )^{1 / 2} ( \bla \circ \bbet ) | ( \X ^{\top } \X )^{- 1} \X ^{\top } \y , \sis ( \X ^{\top } \X )^{- 1} ) \pi _{\ta } ( \tas ) \non \\
&\propto %
\exp \{ - \tas \bbet ^{\top } \bLa ( \X ^{\top } \X ) \bLa \bbet / (2 \sis ) + ( \tas )^{1 / 2} \bbet ^{\top } \bLa \X ^{\top } \y / \sis \} \pi _{\ta } ( \tas ) \text{.} \non 
\end{align}

\begin{algo}
\label{algo:pc-three} 
The parameters $\bth $, $\bbet $, $\sis $, $\tas $, and $\bla $ are updated in the following way. 
\begin{itemize}
\item
Sample 
\begin{align}
\bla &\sim p( \bla | \bth , \bbet , \sis , \tas ) \non \\
&\propto \prod_{k = 1}^{p} \Big( \exp \Big[ - {\{ ( \tas )^{1 / 2} \la _k \bet _k - \th _k \} ^2 \over 2 d \sis } \Big] \pi _{\la } ( \la _k ) \Big) \text{.} \non 
\end{align}
\item
Sample 
\begin{align}
\tas &\sim p( \tas | \bla , \bbet , \sis ) \non \\
&\propto \exp \{ - \tas \bbet ^{\top } \bLa ( \X ^{\top } \X ) \bLa \bbet / (2 \sis ) + ( \tas )^{1 / 2} \bbet ^{\top } \bLa \X ^{\top } \y / \sis \} \pi _{\ta } ( \tas ) \text{.} \non 
\end{align}
\item
Sample 
\begin{align}
\sis &\sim p( \sis | \tas , \bla ) \non \\
&= {\rm{IG}} ( \sis | n / 2 + a' , \y ^{\top } [ \I ^{(n)} - \X \bLa \{ \I ^{(p)} / \tas + \bLa ( \X ^{\top } \X ) \bLa \} ^{- 1} \bLa \X ^{\top } ] \y / 2 + b' ) \text{,} \non 
\end{align}
sample 
\begin{align}
\bbet &\sim p( \bbet | \sis , \tas , \bla ) \non \\
&= {\rm{N}}_p ( \bbet | \{ \I ^{(p)} + \tas \bLa ( \X ^{\top } \X ) \bLa \} ^{- 1} \bLa \X ^{\top } \y ( \tas )^{1 / 2} , \sis \{ \I ^{(p)} + \tas \bLa ( \X ^{\top } \X ) \bLa \} ^{- 1} ) \text{,} \non 
\end{align}
and sample 
\begin{align}
\bth &\sim p( \bth | \bbet , \sis , \tas , \bla ) \non \\
&= {\rm{N}}_p ( \bth | \A ^{- 1} \{ ( \tas )^{1 / 2} ( \bla \circ \bbet ) / d + \b \} , \sis \A ^{- 1} ) \text{,} \non 
\end{align}
where 
\begin{align}
\A &= \I ^{(p)} / d + \{ ( \X ^{\top } \X )^{- 1} - d \I ^{(p)} \} ^{- 1} \quad \text{and} \quad \b = \{ ( \X ^{\top } \X )^{- 1} - d \I ^{(p)} \} ^{- 1} ( \X ^{\top } \X )^{- 1} \X ^{\top } \y \text{.} \non 
\end{align}
\end{itemize}
\end{algo}

When $\pi _{\ta }$ is gamma, $\sqrt{\tas } | ( \bla , \bbet , \sis )$ is distributed as a power truncated normal (PTN) %
(\cite{hpx2022}) 
random variable (see Appendix C of He et al. (2022) for an efficient PTN generator). 
If $\pi _{\ta }$ is a truncated gamma, we would use a modified version of their PTN sampler in which the proposal distribution is tuned using a Newton-Raphson algorithm. 

The conditional densities of $\la _1 , \dots , \la _p$ given by (\ref{eq:conditional_bla}) are also not necessarily standard distributions. 
However, we can sample from them using rejection sampling since they are simple one-dimensional distributions. 
As an example, we consider the horseshoe case here, wherein: 
\begin{align}
&p( \la _k | \bth , \bbet , \sis , \tas ) \propto \exp \Big[ - {\{ ( \tas )^{1 / 2} \la _k \bet _k - \th _k \} ^2 \over 2 d \sis } \Big] {1 \over 1 + {\la _k}^2} \non 
\end{align}
for all $k = 1, \dots , p$. 
Therefore, constructing an Accept-Reject algorithm is sufficient to sample from the class of densities given by 
\begin{align}
f(x | r, s) &= g(x; r, s) / \int_{- \infty }^{\infty } g( \xt ; r, s) d\xt \text{,} \quad x \in \mathbb{R} \text{,} \quad (r, s) \in (0, \infty ) \times \mathbb{R} \text{,} \non 
\end{align}
where 
\begin{align}
g(x; r, s) &= {e^{- x^2 / 2} \over r^2 + (x - s)^2} \text{,} \quad x \in \mathbb{R} \text{,} \non 
\end{align}
for $(r, s) \in (0, \infty ) \times \mathbb{R}$. 
Because such an algorithm will be used in each scan of the Gibbs sampler, the rejection constant should be guaranteed to be bounded as a function of the parameters $(r, s) \in (0, \infty ) \times \mathbb{R}$. 
We achieve this by partitioning the support of the target density into subsets. 
We do not examine the optimization of the rejection constant in this study.

\begin{prp}
\label{prp:normal-cauchy} 
There exists a positive constant $C > 0$ not depending on anything which we can actually calculate such that for all $(r, s) \in (0, \infty ) \times \mathbb{R}$, there exists a normalized density $h_{r, s} \colon \mathbb{R} \to (0, \infty )$ which we can easily evaluate and sample from such that for all $x \in \mathbb{R}$, we have 
\begin{align}
f(x | r, s) &\le C h_{r, s} (x) \text{.} \non 
\end{align}
\end{prp}

\subsection{Geometric ergodicity}
We can prove geometric ergodicity without assuming a negative moment condition on $\tas $. 
This is in contrast to the case of Section \ref{sec:horseshoe}. 

\begin{thm}
\label{thm:general_ergodicity} 
\hfill
\begin{itemize}
\item[{\rm{(i)}}]
The Markov chain based on Algorithm \ref{algo:pc-two} is geometrically ergodic if the following holds: 
\begin{align}
&\int_{0}^{\infty } t^{\ep } \pi _{\ta } (t) dt < \infty \non 
\end{align}
for some $\ep > 0$. 
\item[{\rm{(ii)}}]
The Markov chain based on Algorithm \ref{algo:pc-three} is geometrically ergodic if $\tas \le \overline{\tas }$ for some $0 < \overline{\tas } < \infty $. 
\end{itemize}
\end{thm}

The assumption of part (i) is satisfied by virtually all priors for $\tas $. 
It is sufficient to check the condition for an arbitrarily small $\ep > 0$. 
The assumption of part (ii) is that the support of $\tas $ is bounded above. 
Compatible with guideline (ii) of %
\cite{ps2011}, this assumption is better than the negative moment conditions of Theorem \ref{thm:lat}.

\subsection{An improved sampler}
\label{subsec:improved} 
Algorithm \ref{algo:pc-three} %
produces $\bla \sim p( \bla | \bth , \bbet , \sis , \tas )$, $\tas \sim p( \tas | \bla , \bbet , \sis )$, and then $( \bth , \bbet , \sis ) \sim p( \bth , \bbet , \sis | \tas , \bla )$. 
In order to construct a practically more efficient algorithm, we first consider sampling $( \bla , \bbet ) \sim p( \bla , \bbet | \bth , \sis , \tas )$, %
$( \tas , \sis ) \sim p( \tas , \sis | \bla , \bbet )$, and $( \bth , \bbet , \sis ) \sim p( \bth , \bbet , \sis | \tas , \bla )$. 
Removing the redundant, intermediate draw of $\sis $, we summarize this algorithm as follows. 

\begin{algo}
\label{algo:pc-new} 
The parameters $\bth $, $\bbet $, $\sis $, $\tas $, and $\bla $ are updated in the following way. 
\begin{itemize}
\item
Sample 
\begin{align}
\bla &\sim p( \bla | \bth , \sis , \tas ) \non 
\end{align}
and sample 
\begin{align}
\bbet &\sim p( \bbet | \bth , \sis , \tas , \bla ) \text{.} \non 
\end{align}
\item
Sample 
\begin{align}
\tas &\sim p( \tas | \bla , \bbet ) \text{.} \non 
\end{align}
\item
Sample 
\begin{align}
\sis &\sim p( \sis | \tas , \bla ) \text{,} \non 
\end{align}
sample 
\begin{align}
\bbet &\sim p( \bbet | \sis , \tas , \bla ) \text{,} \non 
\end{align}
and sample 
\begin{align}
\bth &\sim p( \bth | \bbet , \sis , \tas , \bla ) \text{,} \non 
\end{align}
\end{itemize}
\end{algo}

Introducing the latent variable $\bth $ makes the autocorrelation %
higher. 
For instance, when $0 < d \ll 1$ with $\tas $ fixed, we have $\bth \approx ( \tas )^{1 / 2} ( \bla \circ \bbet )$, which means that knowledge of any two of the variables $\bth $, $\bbet $, and $\bla $ almost uniquely determines the remaining one. 
However, we could reduce the autocorrelation by using the above algorithm instead of Algorithm \ref{algo:pc-three}. 
This is because only one of the three variables is conditioned on in the first and third steps of Algorithm \ref{algo:pc-new} so that the joint distribution of the current and new sets of the variables is less connected. 
Algorithm \ref{algo:pc-new} can also be seen as a partially collapsed version of the simpler two-stage Gibbs sampler in which we alternately sample $( \bbet , \bla )$ and $( \bth , \sis , \tas )$. 

The third step of Algorithm \ref{algo:pc-new} can be performed as described in Section \ref{sec:general}. 
For the second step, 
\begin{align}
&\int_{0}^{\infty } p( \tas , \sis | \bla , \bbet ) d\sis \non \\
&\propto \int_{0}^{\infty } {\exp \{ - \y ^{\top } \Q _{\X } \y / (2 \sis ) \} \over ( \sis )^{(n - p) / 2}} {\rm{N}}_p ( \bbeh | ( \X ^{\top } \X )^{- 1} \X ^{\top } \y , \sis ( \X ^{\top } \X )^{- 1} ) {\rm{N}}_p ( \bbet | \bm{0} ^{(p)} , \sis \I ^{(p)} ) p( \tas , \sis ) d\sis \non \\
&= \pi _{\ta } ( \tas ) / \Big( b' + {\y ^{\top } \Q _{\X } \y + \{ \bbeh - ( \X ^{\top } \X )^{- 1} \X ^{\top } \y \} ^{\top } ( \X ^{\top } \X ) \{ \bbeh - ( \X ^{\top } \X )^{- 1} \X ^{\top } \y \} + \| \bbet \| ^2 \over 2} \Big) ^{(n - p) / 2 + a'} \text{,} \non 
\end{align}
where $\bbeh = ( \tas )^{1 / 2} ( \bla \circ \bbet )$. 
Therefore, if $\pi _{\ta } ( \tas ) \propto 1(0 < \tas \le \overline{\tas } ) ( \tas )^{1 / 2 - 1} / (1 + \tas )$ for some $0 < \overline{\tas } < \infty $, we can sample $( \tas )^{1 / 2}$ using rejection sampling. 
The rejection constant is $1 + \overline{\tas }$ and the dominating density is a truncated $t$ density. 

For the first step, 
\begin{align}
p( \bbet | \bth , \sis , \tas , \bla ) %
&\propto \prod_{k = 1}^{p} {\rm{N}} \Big( \bet _k \Big| {\th _k ( \tas )^{1 / 2} \la _k / d \over \tas {\la _k}^2 / d + 1}, {\sis \over \tas {\la _k}^2 / d + 1} \Big) \text{.} \non 
\end{align}
Meanwhile, 
\begin{align}
&\prod_{k = 1}^{p} \int_{- \infty }^{\infty } p( \bet _k , \la _k | \bth , \sis , \tas ) d{\bet _k} \propto \prod_{k = 1}^{p} \Big( {\pi _{\la } ( \la _k ) \over (1 + \tas {\la _k}^2 / d)^{1 / 2}} \exp \Big[ {1 \over 2 \sis } {\{ \th _k ( \tas )^{1 / 2} \la _k / d \} ^2 \over 1 + \tas {\la _k}^2 / d} \Big] \Big) \text{.} \non 
\end{align}
Therefore, in the horseshoe case, 
the conditional densities of $\la _1 , \dots , \la _p$ are within the class of densities given by 
\begin{align}
&f(x | r, s) = {1 \over 1 + x^2} {1 \over (r + x^2 )^{1 / 2}} \exp \Big( s {x^2 \over r + x^2} \Big) / \int_{- \infty }^{\infty } {1 \over 1 + \xt ^2} {1 \over (r + \xt ^2 )^{1 / 2}} \exp \Big( s {\xt ^2 \over r + \xt ^2} \Big) d\xt \text{,} \non 
\end{align}
$x \in \mathbb{R}$, $(r, s) \in (0, \infty )^2$. 

As in Section \ref{sec:general}, we require an accept-reject algorithm and the rejection constant should be bounded. 

\begin{prp}
\label{prp:new} 
There exists a positive constant $C > 0$ not depending on anything which we can actually calculate such that for all $(r, s) \in (0, \infty )^2$, there exists a normalized density $h_{r, s} \colon \mathbb{R} \to (0, \infty )$ which we can easily evaluate and sample from such that for all $x \in \mathbb{R}$, we have 
\begin{align}
f(x | r, s) &\le C h_{r, s} (x) \text{.} \non 
\end{align}
\end{prp}

In practice, in order to reduce the correlation between $\tas $ and $\sum_{k = 1}^{p} \la _k  \bet _k$, we can introduce a working parameter $\al \in (0, \infty )$ by making the change of variables $( \tas , \bla ) = ( \tahs / \al ^2 , \al \blah )$. 
This increases the number of the stages of our algorighm (unless $\al $ and some other parameters can be sampled jointly). 
Meanwhile, reparameterization may or may not improve convergence in general.

More specifically, we can generate 
\begin{align}
\bla &\sim p( \bla | \al , \bth , \sis , \tahs ) \propto \prod_{k = 1}^{p} \Big( {\pi _{\la } ( \la _k ) \over \{ 1 + ( \tahs / \al ^2 ) {\la _k}^2 / d \} ^{1 / 2}} \exp \Big[ {1 \over 2 \sis } {{\th _k}^2 ( \tahs / \al ^2 ) {\la _k}^2 / d^2 \over 1 + ( \tahs / \al ^2 ) {\la _k}^2 / d} \Big] \Big) \text{.} \non 
\end{align}
to sample $\blah = \bla / \al \sim p( \blah | \al , \bth , \sis , \tahs )$. 
Also, we can easily derive a conditional distribution of $\tahs $ from that of $\tas $ given above. 
Finally, under the prior $\al \sim \pi _{\al } ( \al ) = 2 \{ \al ^{2 \ep _{\al } - 1} / (1 + \al ^2 )^{\ep _{\al } + \ep _{\al }} \} / B( \ep _{\al } , \ep _{\al } )$ with $\ep _{\al } = {10}^{- 4}$ in the horseshoe case, we have 
\begin{align}
p( \al | \bth , \bbet , \sis , \tas , \bla ) &\propto %
{\al ^{p + 2 \ep _{\al }} \over (1 + \al ^2 )^{\ep _{\al } + \ep _{\al }}} {1 \over \al ^2 + \tahs } \prod_{k = 1}^{p} {1 \over 1 + {\lah _k}^2 \al ^2} \text{,} \non 
\end{align}
which implies that the full conditional of $\xi = \log ( \al ^2 )$ is log-concave: 
\begin{align}
p( \xi | \bth , \bbet , \sis , \tas , \bla ) &\propto {e^{\{ (p + 1) / 2 + \ep _{\al } \} \xi } \over (1 + e^{\xi } )^{\ep _{\al } + \ep _{\al }}} {1 \over e^{\xi } + \tahs } \prod_{k = 1}^{p} {1 \over 1 + {\lah _k}^2 e^{\xi }} \text{;} \non 
\end{align}
therefore, we can use the method of \cite{gw1992}.

\section{A Simulation Study}
\label{sec:sim} 
Here, we compare the first method of Section \ref{sec:horseshoe} (new1), the last method of Section \ref{subsec:improved} (new2), the modified JOB method of \cite{bkp2022} (mjob), and the original, unmodified JOB method of \cite{job2020} (ujob). 
For new2, we ignore the upper bound $\tas \le \overline{\tas }$ by letting $\overline{\tas } \to \infty $, so that the posterior distiruibution is exactly the same as that of the other methods. 
For new2, mjob, and ujob, we use rejection sampling based on Lemma \ref{lem:naive} of Section \ref{sec:naive_rejection} of the Supplementary Material, instead of following \cite{bkp2022} to perform Metropolis steps, in updating $\tas $, which makes it easier to see how different parameterizations compare in efficiency. 
(We note that our approaches based on data augmentation with $\bth $ are more suited for rejection sampling since we do not have to compute determinants.) 
For ujob, we use the method of \cite{d2021} instead of that of Appendix S1 of \cite{job2020} in updating $( \la _{k}^{2} )_{k = 1}^{p}$; see Section \ref{sec:rejecton_expint} of the Supplementary Material of this paper. 
We focus on effective sample sizes in this paper; computation times will be highly dendent upon the details of the algorithms and how codes are written, especially when rejection sampling is used.

We basically use the setting of Section 2.3 of \cite{bkp2022} but we set $p \le n$ because our methods have been developed for this case. 
Also, we set the hyperparameters $( a' , b' )$ equal to $(1 / 10, 1 / 10)$. 
We assume the half Cauchy prior distribution for $\sqrt{\tas }$ and $\sqrt{\la _{1}^{2}} , \dots , \sqrt{\la _{p}^{2}}$. 
Five (instead of ten) datasets are generated as follows: $n = 100$; $p = 25, 75$; $\bbe ^0 = ((- 1 + [ \{ 3 - (- 1) \} / 11] k)_{k = 1, \dots , 10} , ( \bm{0} ^{(p - 10)} )^{\top } )^{\top }$; $\X \sim \{ {\rm{N}} (0, 1) \} ^{n \times p}$; $\y \sim {\rm{N}}_n ( \X \bbe ^0 , \I ^{(n)} / 100)$. 
We generate $2,000$ samples after discarding $100$ samples as burn-in. 

\begin{figure}[H]
\centering
\includegraphics[width = 0.65\linewidth]{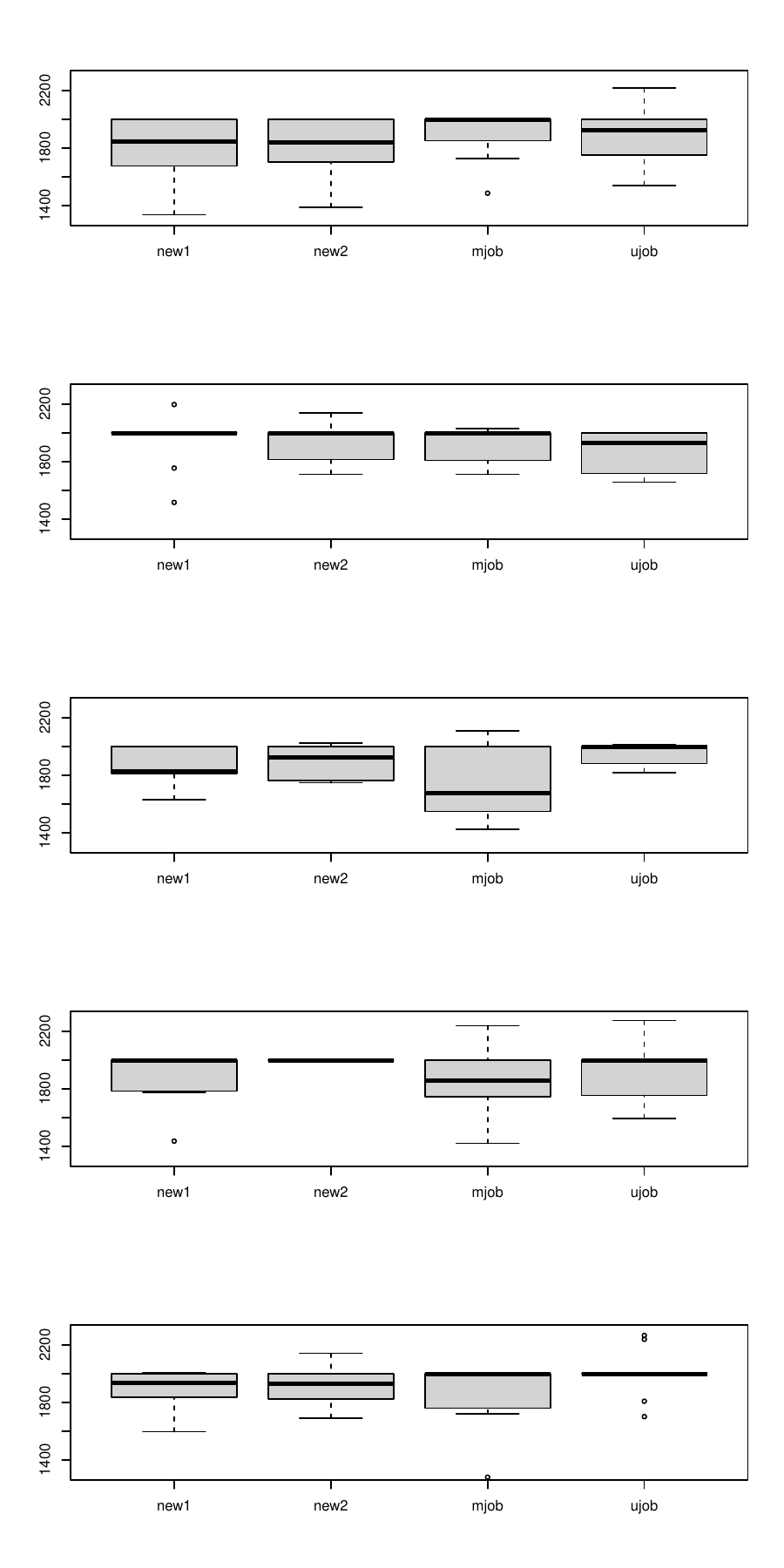}
\caption{Boxplots of effective sample sizes of $\be _1 , \dots , \be _{10}$ for the methods of Section \ref{sec:horseshoe} (new1), Section \ref{subsec:improved} (new2), \cite{bkp2022} (mjob), and \cite{job2020} (ujob) based on five datasets with $(n, p) = (100, 25)$. }
\label{fig:be_small} 
\end{figure}

\begin{figure}[H]
\centering
\includegraphics[width = 0.65\linewidth]{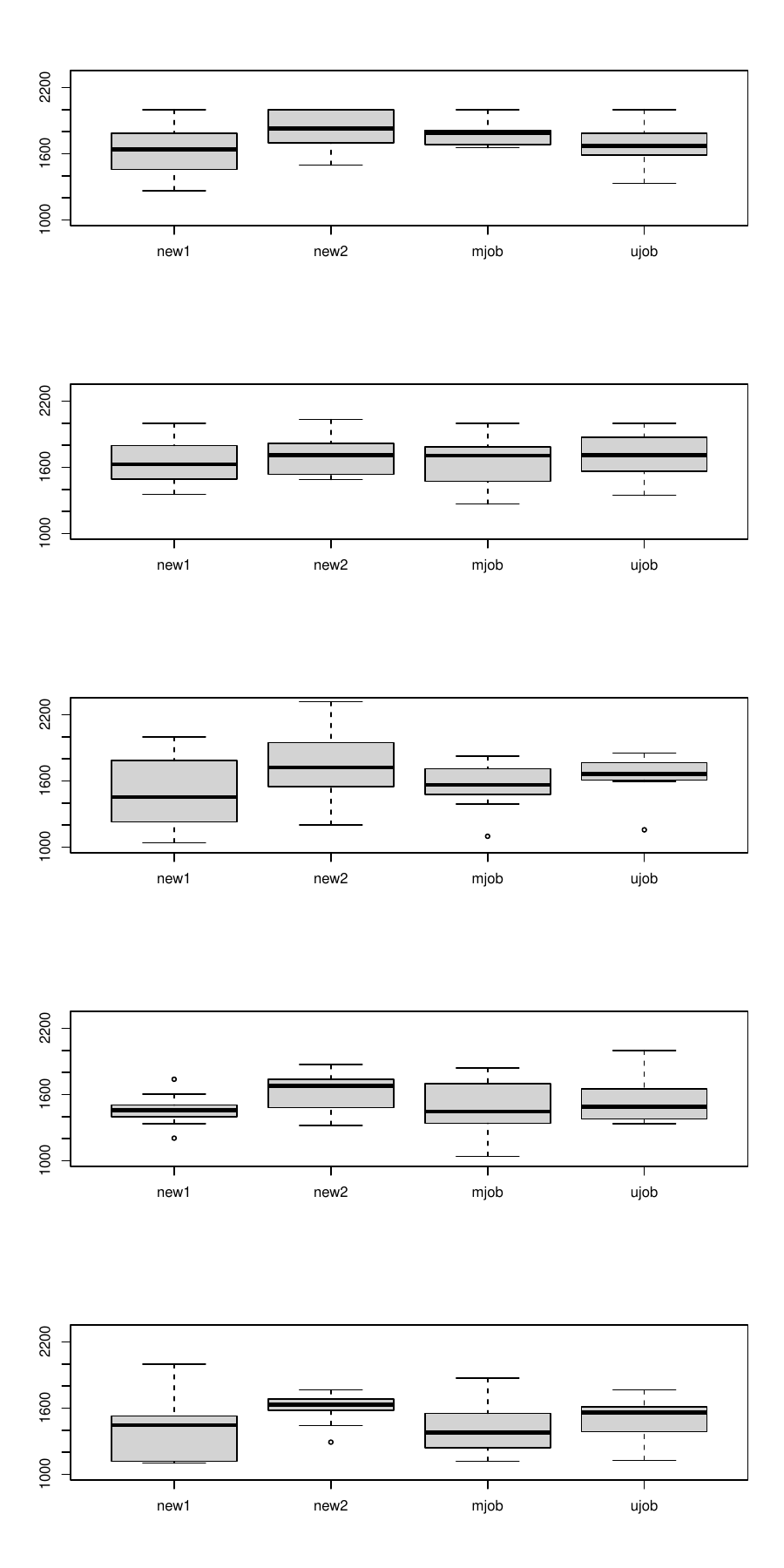}
\caption{Boxplots of effective sample sizes of $\be _1 , \dots , \be _{10}$ for the methods of Section \ref{sec:horseshoe} (new1), Section \ref{subsec:improved} (new2), \cite{bkp2022} (mjob), and \cite{job2020} (ujob) based on five datasets with $(n, p) = (100, 75)$. }
\label{fig:be_large} 
\end{figure}

Boxplots of effective sample sizes of $\be _1 , \dots , \be _{10}$ are shown in Figures \ref{fig:be_small} and \ref{fig:be_large}. 
The results for the four methods are similar when $p = 25$. 
When $p = 75$, all four do a good job but new2 is better than the other methods.

Effective sample sizes of $\sis $ and $\tas $ are reported in Tables \ref{table:sis} and \ref{table:tas}, respectively. 
The JOB methods are better when we consider $\sis $. 
(It would be interesting to investigate if our methods can be further improved by introducing a working parameter related to $\bbet $ and $\sis $). 
Meanwhile, our methods are better for efficiently sampling $\tas $ when $p = 25$. 
When $p = 75$, ujob is good, new1 is bad (the omitted value is $2,000$ probably because of some numerical error), and new2 is somewhat better than mjob.

\small
\begin{table}[!thb]
\caption{Effective sample sizes of $\sis $ for the methods of Section \ref{sec:horseshoe} (new1), Section \ref{subsec:improved} (new2), \cite{bkp2022} (mjob), and \cite{job2020} (ujob) based on five datasets with $(n, p) = (100, 25)$ (left) and five datasets with $(n, p) = (100, 25)$ (right). 
Each row of the left and right blocks corresponds to one dataset. }
\begin{center}
$
{\renewcommand\arraystretch{1.1}\small
\begin{array}{c@{\hspace{2mm}}
              r@{\hspace{2mm}}
              r@{\hspace{2mm}}
              r@{\hspace{2mm}}
              r@{\hspace{2mm}}
              r@{\hspace{2mm}}
              r@{\hspace{2mm}}
              r@{\hspace{2mm}}
              r@{\hspace{2mm}}
              r@{\hspace{2mm}}
              r@{\hspace{2mm}}
              r@{\hspace{2mm}}
              r
             }
\hline
\multicolumn{4}{c}{p = 25} &  &  & \multicolumn{4}{c}{p = 75} \\

\hline 
\text{new1} & \text{new2} & \text{mjob} & \text{ujob} &  &  & \text{new1} & \text{new2} & \text{mjob} & \text{ujob} \\

\hline

1392 & 1549 & 1752 & 2000 &  &  & 620 & 1014 & 1791 & 2000 \\
973 & 1231 & 2000 & 2000 &  &  & 496 & 981 & 1015 & 2000 \\
1416 & 1442 & 2000 & 2000 &  &  & 784 & 577 & 1433 & 1739 \\
1315 & 1294 & 2000 & 2000 &  &  & 612 & 519 & 1156 & 1426 \\
1356 & 1365 & 2138 & 2000 &  &  & 334 & 581 & 1507 & 1442 \\

\hline
\end{array}
}
$
\end{center}
\label{table:sis} 
\end{table}
\normalsize

\small
\begin{table}[!thb]
\caption{Effective sample sizes of $\tas $ for the methods of Section \ref{sec:horseshoe} (new1), Section \ref{subsec:improved} (new2), \cite{bkp2022} (mjob), and \cite{job2020} (ujob) based on five datasets with $(n, p) = (100, 25)$ (left) and five datasets with $(n, p) = (100, 25)$ (right). 
Each row of the left and right blocks corresponds to one dataset. }
\begin{center}
$
{\renewcommand\arraystretch{1.1}\small
\begin{array}{c@{\hspace{2mm}}
              r@{\hspace{2mm}}
              r@{\hspace{2mm}}
              r@{\hspace{2mm}}
              r@{\hspace{2mm}}
              r@{\hspace{2mm}}
              r@{\hspace{2mm}}
              r@{\hspace{2mm}}
              r@{\hspace{2mm}}
              r@{\hspace{2mm}}
              r@{\hspace{2mm}}
              r@{\hspace{2mm}}
              r
             }
\hline
\multicolumn{4}{c}{p = 25} &  &  & \multicolumn{4}{c}{p = 75} \\

\hline 
\text{new1} & \text{new2} & \text{mjob} & \text{ujob} &  &  & \text{new1} & \text{new2} & \text{mjob} & \text{ujob} \\

\hline

552 & 726 & 109 & 145 &  &  & - & 250 & 178 & 255 \\
528 & 888 & 185 & 173 &  &  & 179 & 349 & 240 & 247 \\
419 & 969 & 127 & 189 &  &  & 117 & 153 & 188 & 238 \\
642 & 788 & 175 & 154 &  &  & 117 & 273 & 218 & 303 \\
613 & 746 & 149 & 213 &  &  & 115 & 256 & 208 & 238 \\

\hline
\end{array}
}
$
\end{center}
\label{table:tas} 
\end{table}
\normalsize

Finally, we remark again that it is difficult to compare the computation times of the four methods. 
Also, different methods are good for the efficient sampling of different variables. 
In interpreting the results for $\bbe $, we recall that $\bbet $ is sampled twice in each scan when we use new2.

\section{Discussion}
\label{sec:discussion} 
In the first half of this study, we considered a Gibbs algorithm for a horseshoe model under the usual $( \bbe , \sis )$ parameterization, in which the global shrinkage parameter $\tas $ is easy to sample. 
We showed that in proving geometric ergodicity of the algorithm, we do not have to assume that $\tas $ has a finite approximately $(p / 2)$-th negative moment, which is a novel result. 
In the second half of this study, we saw that if we are willing to use reparameterization and rejection sampling for one-dimensional distributions, we can establish geometric ergodicity %
for a normal linear regression model with virtually any global-local shrinkage prior. 
We also constructed a practical sampling method. 
As in many studies, we did not address the problem of obtaining good quantitative bounds on the rate of %
convergence, which is a relevant issue.

\section*{Acknowledgments}
Research of the author was supported in part by JSPS KAKENHI Grant Numbers JP25K21163, JP22K20132 from Japan Society for the Promotion of Science.

\newpage
\setcounter{page}{1}
\setcounter{section}{0}
\renewcommand{\thesection}{S\arabic{section}}
\setcounter{table}{0}
\renewcommand{\thetable}{S\arabic{table}}
\setcounter{figure}{0}
\renewcommand{\thefigure}{S\arabic{figure}}

\begin{center}
{\LARGE\bf Supplementary Materials for ``Geometric Ergodicity of Gibbs Algorithms for a Normal Model With a Global-Local Shrinkage Prior''}
\end{center}

\medskip

\begin{center}
\Large
Yasuyuki Hamura\footnote{Graduate School of Economics, Kyoto University, 
Yoshida-Honmachi, Sakyo-ku, Kyoto, 606-8501, JAPAN. 
\newline{
E-Mail: yasu.stat@gmail.com}} \
\end{center}

\bigskip

In Section \ref{sec:proof_general_prp}, we prove Proposition \ref{prp:normal-cauchy}. 
In Section \ref{sec:details_rejection_proposed}, we give details of an accept-reject algorithm corresponding to Proposition \ref{prp:new}, which is proved in Section \ref{sec:proof_general_prp_new}. 
In Section \ref{sec:naive_rejection}, we describe an approach to sampling $\tas $. 
In Section \ref{sec:rejecton_expint}, we rewrite $p( \bla ^2 | \bbe , \sis , \tas , \y )$. 
In Section \ref{sec:proof_thm_first}, we prove Theorem \ref{thm:lat}. 
In Section \ref{sec:failure}, we demonstrate that establishing geometric ergodicity for Algorithm \ref{algo:lat} may be difficult when $\pi _{\ta }$ has a heavier left tail. 
In Section \ref{sec:proof_thm_second}, we prove Theorem \ref{thm:general_ergodicity}.

\section{Proof of Proposition \ref{prp:normal-cauchy}}
\label{sec:proof_general_prp} 
If $A$ is a subset of $\mathbb{R}$, then ${\rm{TN}} (0, 1; A)$ denotes the standard normal distribution truncated to $A$ and ${\rm{TC}} (r, s; A)$ denotes the truncated Cauchy distribution with density proportional to $1(x \in A) / \{ r^2 + (x - s)^2 \} $, $x \in \mathbb{R}$, for $(r, s) \in (0, \infty ) \times \mathbb{R}$. 

Without loss of generality, we assume that $s \ge 0$. 
Fix $x \in \mathbb{R}$. 
To prove Proposition \ref{prp:normal-cauchy}, the basic idea is that for any density $f \colon \mathbb{R} \to (0, \infty )$, constants $C_1 , C_2 , C_3 > 0$, intervals $A_1 , A_2 , A_3 \subset \mathbb{R}$ satisfying $A_1 \cup A_2 \cup A_3 = \mathbb{R}$, %
and normalized densities $h_1 , h_2 , h_3 \colon \mathbb{R} \to [0, \infty )$, the condition that 
\begin{align}
f(x) 1(x \in A_j ) &\le C_j h_j (x) \non 
\end{align}
for all $j = 1, 2, 3$ implies the following: 
\begin{align}
f(x) &\le ( C_1 + C_2 + C_3 ) \sum_{j = 1}^{3} {C_j \over C_1 + C_2 + C_3} h_j (x) \text{.} \non 
\end{align}

\subsection{The case of $s > 1$}
\subsubsection{The case of $s > 1$ and $r > 1$}
\label{subsubsec:pp} 
First, 
\begin{align}
f(x | r, s) 1(x \in (- \infty , s / 2)) &\le {e^{- x^2 / 2} 1(x \in (- \infty , s / 2)) \over r^2 + (x - s)^2} / \int_{0}^{s / 2} {e^{- \xt ^2 / 2} \over r^2 + ( \xt - s)^2} d\xt \non \\
&\le {e^{- x^2 / 2} 1(x \in (- \infty , s / 2)) \over r^2 + s^2 / 4} / \int_{0}^{s / 2} {e^{- \xt ^2 / 2} \over r^2 + s^2} d\xt \non \\
&= {r^2 + s^2 \over r^2 + s^2 / 4} {\int_{- \infty }^{s / 2} e^{- \xt ^2 / 2} d\xt \over \int_{0}^{s / 2} e^{- \xt ^2 / 2} d\xt} {e^{- x^2 / 2} 1(x \in (- \infty , s / 2)) \over \int_{- \infty }^{s / 2} e^{- \xt ^2 / 2} d\xt } \non \\
&\le (1 + 4) \Big( \sup_{\st \in [1, \infty )} {\int_{- \infty }^{\st / 2} e^{- \xt ^2 / 2} d\xt \over \int_{0}^{\st / 2} e^{- \xt ^2 / 2} d\xt} \Big) {\rm{TN}} (x | 0, 1; (- \infty , s / 2)) \text{.} \non 
\end{align}
Second, 
\begin{align}
&f(x | r, s) 1(x \in [s / 2, \infty )) \non \\
&\le {e^{- x^2 / 2} 1(x \in [s / 2, \infty )) \over r^2} / \int_{- \infty }^{\infty } {e^{- \xt ^2 / 2} \over r^2 + ( \xt - s)^2} d\xt \non \\
&\le e^{- x^2 / 2} 1(x \in [s / 2, \infty )) / \int_{- \infty }^{\infty } {e^{- \xt ^2 / 2} \over 1 + ( \xt - s)^2} d\xt \non \\
&= \Big\{ (1 + s^2) \int_{s / 2}^{\infty } e^{- \xt ^2 / 2} d\xt / \int_{- \infty }^{\infty } {(1 + s^2) e^{- \xt ^2 / 2} \over 1 + ( \xt - s)^2} d\xt \Big\} {e^{- x^2 / 2} 1(x \in [s / 2, \infty )) \over \int_{s / 2}^{\infty } e^{- \xt ^2 / 2} d\xt } \non \\
&\le \Big[ \sup_{\st \in [1, \infty )} \Big\{ (1 + \st ^2) \int_{\st / 2}^{\infty } e^{- \xt ^2 / 2} d\xt / \int_{- \infty }^{\infty } {(1 + \st ^2) e^{- \xt ^2 / 2} \over 1 + 2 ( \xt^2 + \st ^2 )} d\xt \Big\} \Big] {\rm{TN}} (x | 0, 1; [s / 2, \infty )) \text{.} \non 
\end{align}

\subsubsection{The case of $s > 1$ and $r \le 1$}
First, as in Section \ref{subsubsec:pp}, we have 
\begin{align}
f(x | r, s) 1(x \in (- \infty , s / 2)) &\le (1 + 4) \Big( \sup_{\st \in [1, \infty )} {\int_{- \infty }^{\st / 2} e^{- \xt ^2 / 2} d\xt \over \int_{0}^{\st / 2} e^{- \xt ^2 / 2} d\xt} \Big) {\rm{TN}} (x | 0, 1; (- \infty , s / 2)) \text{.} \non 
\end{align}
Second, 
\begin{align}
&f(x | r, s) 1(x \in [s / 2, \infty ) \setminus [s - 1 / s, s + 1 / s]) \non \\
&\le {e^{- x^2 / 2} 1(x \in [s / 2, \infty ) \setminus [s - 1 / s, s + 1 / s]) \over 1 / s^2} / \int_{- \infty }^{\infty } {e^{- \xt ^2 / 2} \over 1 + ( \xt - s)^2} d\xt \non \\
&\le \Big\{ s^2 \int_{[s / 2, \infty )} e^{- \xt ^2 / 2} d\xt / \int_{- \infty }^{\infty } {e^{- \xt ^2 / 2} \over 1 + ( \xt - s)^2} d\xt \Big\} {e^{- x^2 / 2} 1(x \in [s / 2, \infty ) \setminus [s - 1 / s, s + 1 / s]) \over \int_{[s / 2, \infty ) \setminus [s - 1 / s, s + 1 / s]} e^{- \xt ^2 / 2} d\xt } \non \\
&\le \Big[ \sup_{\st \in [1, \infty )} \Big\{ \st ^2 (1 + \st ^2 ) \int_{[ \st  / 2, \infty )} e^{- \xt ^2 / 2} d\xt / \int_{- \infty }^{\infty } {(1 + \st ^2 ) e^{- \xt ^2 / 2} \over 1 + 2 ( \xt ^2 + \st ^2 )} d\xt \Big\} \Big] \non \\
&\quad \times {\rm{TN}} (x | 0, 1; [s / 2, \infty ) \setminus [s - 1 / s, s + 1 / s]) \text{.} \non 
\end{align}
Third, 
\begin{align}
&f(x | r, s) 1(x \in [s / 2, \infty ) \cap [s - 1 / s, s + 1 / s]) \non \\
&\le {e^{- x^2 / 2} 1(x \in [s - 1 / s, s + 1 / s]) \over r^2 + (x - s)^2} / \int_{s - 1 / s}^{s + 1 / s} {e^{- \xt ^2 / 2} \over r^2 + ( \xt - s)^2} d\xt \non \\
&\le {e^{- (s - 1 / s)^2 / 2} \over e^{- (s + 1 / s)^2 / 2}} {1(x \in [s - 1 / s, s + 1 / s]) \over r^2 + (x - s)^2} / \int_{s - 1 / s}^{s + 1 / s} {1 \over r^2 + ( \xt - s)^2} d\xt \non \\
&= e^2 {\rm{TC}} (x | r, s; [s - 1 / s, s + 1 / s]) \text{.} \non 
\end{align}

\subsection{The case of $s \le 1$}
\subsubsection{The case of $s \le 1$ and $r > 1$}
We have 
\begin{align}
f(x | r, s) &\le {e^{- x^2 / 2} \over r^2} / \int_{- \infty }^{\infty } {e^{- \xt ^2 / 2} \over r^2 + 2 ( \xt^2 + s^2 )} d\xt \non \\
&\le e^{- x^2 / 2} / \int_{- \infty }^{\infty } {e^{- \xt ^2 / 2} \over 1 + 2 ( \xt^2 + 1)} d\xt \non \\
&= \Big\{ (2 \pi )^{1 / 2} / \int_{- \infty }^{\infty } {e^{- \xt ^2 / 2} \over 1 + 2 ( \xt^2 + 1)} d\xt \Big\} {\rm{N}} (x | 0, 1) \text{.} \non 
\end{align}

\subsubsection{The case of $s \le 1$ and $r \le 1$}
First, 
\begin{align}
&f(x | r, s) 1(x \in \mathbb{R} \setminus [s - 1 , s + 1]) \non \\
&\le e^{- x^2 / 2} 1(x \in \mathbb{R} \setminus [s - 1 , s + 1]) / \int_{- \infty }^{\infty } {e^{- \xt ^2 / 2} \over 1 + 2 ( \xt ^2 + 1)} d\xt \non \\
&\le \Big\{ \int_{- \infty }^{\infty } e^{- \xt ^2 / 2} d\xt / \int_{- \infty }^{\infty } {e^{- \xt ^2 / 2} \over 1 + 2 ( \xt ^2 + 1)} d\xt \Big\} {\rm{TN}} (x | 0, 1; \mathbb{R} \setminus [s - 1 , s + 1]) \text{.} \non 
\end{align}
Second, 
\begin{align}
f(x | r, s) 1(x \in [s - 1 , s + 1]) &\le {1(x \in [s - 1 , s + 1]) \over r^2 + (x - s)^2} / \int_{s - 1}^{s + 1} {e^{- \xt ^2 / 2} \over r^2 + ( \xt - s)^2} d\xt \non \\
&\le e^{2^2 / 2} {\rm{TC}} (x | r, s; [s - 1 , s + 1]) \text{.} \non 
\end{align}

\section{Details of an Accept-Reject Algorithm Corresponding to Proposition \ref{prp:new}}
\label{sec:details_rejection_proposed} 
Here, we give details of the accept-reject algorithm for updating $\la _1 , \dots , \la _p$ which is used for implementing the new2 method of Section \ref{sec:sim}. 
Fix $(r, s) \in (0, \infty )^2$. 
We have 
\begin{align}
&f(x | r, s) = \frac{ \displaystyle {1 \over 1 + x^2} {1 \over (r + x^2 )^{1 / 2}} \exp \Big( s {x^2 \over r + x^2} \Big) }{ \displaystyle \int_{- \infty }^{\infty } {1 \over 1 + \xt ^2} {1 \over (r + \xt ^2 )^{1 / 2}} \exp \Big( s {\xt ^2 \over r + \xt ^2} \Big) d\xt } \non 
\end{align}
for all $x \in \mathbb{R}$. 
We take the same approach as that of Section \ref{sec:proof_general_prp}.

\subsection{The case of $s \le 1$}
Fix $x \in \mathbb{R}$. 

\subsubsection{The case of $s \le 1$ and $r > 1$}
We have 
\begin{align}
f(x | r, s) &\le \frac{ \displaystyle {1 \over 1 + x^2} {e \over r^{1 / 2}} }{ \displaystyle \int_{- \infty }^{\infty } {1 \over 1 + \xt ^2} d\xt } \frac{ \displaystyle \int_{- \infty }^{\infty } {1 \over 1 + \xt ^2} d\xt }{ \displaystyle \int_{- \infty }^{\infty } {1 \over 1 + \xt ^2} {1 \over (r + \xt ^2 )^{1 / 2}} \exp \Big( s {\xt ^2 \over r + \xt ^2} \Big) d\xt } \text{.} \non 
\end{align}

\begin{remark}
Fix $q_0 \in \mathbb{R}$ and $p_0 \in (0, 1)$. 
Then 
\begin{align}
&\int_{- \infty }^{q_0} {1 \over 1 + \xt ^2} d\xt = \arctan q_0 - \Big( - {\pi \over 2} \Big) \text{.} \non 
\end{align}
Therefore, 
\begin{align}
&p_0 = \int_{- \infty }^{q_0} {1 \over 1 + \xt ^2} d\xt / \int_{- \infty }^{\infty } {1 \over 1 + \xt ^2} d\xt \non 
\end{align}
if and only if 
\begin{align}
q_0 &= \tan ( \pi p_0 - \pi / 2) \text{.} \non 
\end{align}
\end{remark}

\subsubsection{The case of $s \le 1$ and $r \le 1$}
We have 
\begin{align}
1(|x| > 1) f(x | r, s) &\le \frac{ \displaystyle 1(|x| > 1) {1 \over 1 + x^2} e }{ \displaystyle \int_{\mathbb{R} \setminus [- 1, 1]} {1 \over 1 + \xt ^2} d\xt } \frac{ \displaystyle \int_{\mathbb{R} \setminus [- 1, 1]} {1 \over 1 + \xt ^2} d\xt }{ \displaystyle \int_{- \infty }^{\infty } {1 \over 1 + \xt ^2} {1 \over (r + \xt ^2 )^{1 / 2}} \exp \Big( s {\xt ^2 \over r + \xt ^2} \Big) d\xt } \non 
\end{align}
and 
\begin{align}
&1(|x| \le 1) f(x | r, s) \le \frac{ \displaystyle 1(|x| \le 1) {1 \over (r + x^2 )^{1 / 2}} e }{ \displaystyle \int_{- 1}^{1} {1 \over (r + \xt ^2 )^{1 / 2}} d\xt } \frac{ \displaystyle \int_{- 1}^{1} {1 \over (r + \xt ^2 )^{1 / 2}} d\xt }{ \displaystyle \int_{- \infty }^{\infty } {1 \over 1 + \xt ^2} {1 \over (r + \xt ^2 )^{1 / 2}} \exp \Big( s {\xt ^2 \over r + \xt ^2} \Big) d\xt } \text{.} \non 
\end{align}

\begin{remark}
The first component is a symmetric density. 
Fix $q_0 \in (1, \infty )$ and $p_0 \in (0, 1)$. 
Then 
\begin{align}
&\int_{1}^{q_0} {1 \over 1 + \xt ^2} d\xt = \arctan q_0 - {\pi \over 4} \non 
\end{align}
Therefore, 
\begin{align}
&p_0 = \int_{1}^{q_0} {1 \over 1 + \xt ^2} d\xt / \int_{1}^{\infty } {1 \over 1 + \xt ^2} d\xt \non 
\end{align}
if and only if 
\begin{align}
q_0 &= \tan \{ ( \pi / 4) p_0 + \pi / 4 \} \text{.} \non 
\end{align}
\end{remark}

\begin{remark}
Fix $q_0 \in (- 1, 1)$ and $p_0 \in (0, 1)$. 
Then 
\begin{align}
\int_{- 1}^{q_0} {1 \over (r + \xt ^2 )^{1 / 2}} d\xt &= \log {q_0 / r^{1 / 2} + (1 + {q_0}^2 / r)^{1 / 2} \over - 1 / r^{1 / 2} + (1 + 1 / r)^{1 / 2}} \text{.} \non 
\end{align}
Therefore, 
\begin{align}
&p_0 = \int_{- 1}^{q_0} {1 \over (r + \xt ^2 )^{1 / 2}} d\xt / \int_{- 1}^{1} {1 \over (r + \xt ^2 )^{1 / 2}} d\xt \non 
\end{align}
if and only if 
\begin{align}
q_0 &= {r^{1 / 2} \over 2} \Big( \exp \Big[ p_0 \log {1 / r^{1 / 2} + (1 + 1 / r)^{1 / 2} \over - 1 / r^{1 / 2} + (1 + 1 / r)^{1 / 2}} + \log \{ - 1 / r^{1 / 2} + (1 + 1 / r)^{1 / 2} \} \Big] \non \\
&\quad - 1 / \exp \Big[ p_0 \log {1 / r^{1 / 2} + (1 + 1 / r)^{1 / 2} \over - 1 / r^{1 / 2} + (1 + 1 / r)^{1 / 2}} + \log \{ - 1 / r^{1 / 2} + (1 + 1 / r)^{1 / 2} \} \Big] \Big) \text{.} \non 
\end{align}
\end{remark}

\subsection{The case of $s > 1$}
Note that if $X$ is distributed with density $f(x | r, s)$, $x \in \mathbb{R}$, then $Y = s |X|^2 / (r + |X|^2 )$ is distributed with density 
\begin{align}
g(y | r, s) &= \frac{ \displaystyle 1(0 < y < s) {y^{1 / 2 - 1} e^y \over 1 - y / s + r y / s} }{ \displaystyle \int_{0}^{s} {\yt ^{1 / 2 - 1} e^{\yt } \over 1 - \yt / s + r \yt / s} d\yt } \text{,} \quad y \in (0, \infty ) \text{.} \label{eq:prop_new_2} 
\end{align}
Fix $y \in (0, \infty )$. 

\subsubsection{The case of $s > 1$ and $r > 1$}
We have 
\begin{align}
&g(y | r, s) = \frac{ \displaystyle 1(0 < y < s) {y^{1 / 2 - 1} e^y \over 1 + (r - 1) y / s} }{ \displaystyle \int_{0}^{s} {\yt ^{1 / 2 - 1} e^{\yt } \over 1 + (r - 1) \yt / s} d\yt } \text{.} \non 
\end{align}
It follows that 
\begin{align}
&1(0 < y < 1) g(y | r, s) \le \frac{ \displaystyle 1(0 < y < 1) {y^{1 / 2 - 1} \over 1 + (r - 1) y / s} e }{ \displaystyle \int_{0}^{1} {\yt ^{1 / 2 - 1} \over 1 + (r - 1) \yt / s} d\yt } \frac{ \displaystyle \int_{0}^{1} {\yt ^{1 / 2 - 1} \over 1 + (r - 1) \yt / s} d\yt }{ \displaystyle \int_{0}^{s} {\yt ^{1 / 2 - 1} e^{\yt } \over 1 + (r - 1) \yt / s} d\yt } \text{,} \non 
\end{align}
that 
\begin{align}
1(s / 2 < y < s) g(y | r, s) &\le \frac{ \displaystyle 1(s / 2 < y < s) e^y {(s / 2)^{1 / 2 - 1} \over 1 + (r - 1) / 2} }{ \displaystyle \int_{s / 2}^{s} e^{\yt } d\yt } \frac{ \displaystyle \int_{s / 2}^{s} e^{\yt } d\yt }{ \displaystyle \int_{0}^{s} {\yt ^{1 / 2 - 1} e^{\yt } \over 1 + (r - 1) \yt / s} d\yt } \text{,} \non 
\end{align}
and that 
\begin{align}
1(1 < y < s / 2) g(y | r, s) &\le \frac{ \displaystyle 1(1 < y < s / 2) {y^{1 / 2 - 1} e^{s / 2} \over 1 + (r - 1) / s} }{ \displaystyle \int_{1}^{s / 2} \yt ^{1 / 2 - 1} d\yt } \frac{ \displaystyle \int_{1}^{s / 2} \yt ^{1 / 2 - 1} d\yt }{ \displaystyle \int_{0}^{s} {\yt ^{1 / 2 - 1} e^{\yt } \over 1 + (r - 1) \yt / s} d\yt } \text{.} \non 
\end{align}

\begin{remark}
Fix $q_0 \in (0, 1)$ and $p_0 \in (0, 1)$. 
Then 
\begin{align}
&\int_{0}^{q_0} {\yt ^{1 / 2 - 1} \over 1 + (r - 1) \yt / s} d\yt = {2 s^{1 / 2} \over (r - 1)^{1 / 2}} \arctan [ \{ (r - 1) / s \} ^{1 / 2} {q_0}^{1 / 2}] \text{.} \non 
\end{align}
Therefore, 
\begin{align}
&p_0 = \int_{0}^{q_0} {\yt ^{1 / 2 - 1} \over 1 + (r - 1) \yt / s} d\yt / \int_{0}^{1} {\yt ^{1 / 2 - 1} \over 1 + (r - 1) \yt / s} d\yt \non 
\end{align}
if and only if 
\begin{align}
&q_0 = {\{ \tan ( p_0 \arctan [ \{ (r - 1) / s \} ^{1 / 2} ]) \} ^2 \over (r - 1) / s} \non 
\end{align}
\end{remark}

\begin{remark}
Fix $q_0 \in (s / 2, s)$ and $p_0 \in (0, 1)$. 
Then 
\begin{align}
&\int_{s / 2}^{q_0} e^{\yt } d\yt = e^{q_0} - e^{s / 2} \text{.} \non 
\end{align}
Therefore, 
\begin{align}
&p_0 = \int_{s / 2}^{q_0} e^{\yt } d\yt / \int_{s / 2}^{s} e^{\yt } d\yt \non 
\end{align}
if and only if 
\begin{align}
&q_0 = \log \{ ( e^s - e^{s / 2} ) p_0 + e^{s / 2} \} \text{.} \non 
\end{align}
\end{remark}

\begin{remark}
Suppose that $s / 2 > 1$. 
Fix $q_0 \in (1, s / 2)$ and $p_0 \in (0, 1)$. 
Then 
\begin{align}
&\int_{1}^{q_0} \yt ^{1 / 2 - 1} d\yt = 2 ( {q_0}^{1 / 2} - 1) \text{.} \non 
\end{align}
Therefore, 
\begin{align}
&p_0 = \int_{1}^{q_0} \yt ^{1 / 2 - 1} d\yt / \int_{1}^{s / 2} \yt ^{1 / 2 - 1} d\yt \non 
\end{align}
if and only if 
\begin{align}
&q_0 = [ \{ (s / 2)^{1 / 2} - 1 \} p_0 + 1]^2 \text{.} \non 
\end{align}
\end{remark}

\subsubsection{The case of $s > 1$ and $1 / 2 < r \le 1$}
We have 
\begin{align}
g(y | r, s) &= \frac{ \displaystyle 1(0 < y < s) {y^{1 / 2 - 1} e^y \over 1 - (1 - r) y / s} }{ \displaystyle \int_{0}^{s} {\yt ^{1 / 2 - 1} e^{\yt } \over 1 - (1 - r) \yt / s} d\yt } \text{.} \non 
\end{align}
We have 
\begin{align}
&1(0 < y < s / 2) g(y | r, s) \le \frac{ \displaystyle 1(0 < y < s / 2) {y^{1 / 2 - 1} e^{s / 2} \over 1 - (1 - r) / 2} }{ \displaystyle \int_{0}^{s / 2} \yt ^{1 / 2 - 1} d\yt } \frac{ \displaystyle \int_{0}^{s / 2} \yt ^{1 / 2 - 1} d\yt }{ \displaystyle \int_{0}^{s} {\yt ^{1 / 2 - 1} e^{\yt } \over 1 - (1 - r) \yt / s} d\yt } \text{.} \non 
\end{align}
and 
\begin{align}
&1(s / 2 < y < s) g(y | r, s) \le \frac{ \displaystyle 1(s / 2 < y < s) {(s / 2)^{1 / 2 - 1} e^y \over r} }{ \displaystyle \int_{s / 2}^{s} e^{\yt } d\yt } \frac{ \displaystyle \int_{s / 2}^{s} e^{\yt } d\yt }{ \displaystyle \int_{0}^{s} {\yt ^{1 / 2 - 1} e^{\yt } \over 1 - (1 - r) \yt / s} d\yt } \text{.} \non 
\end{align}

\begin{remark}
The second component is considered in the previous case. 
For the first component, fix $q_0 \in (0, s / 2)$ and $p_0 \in (0, 1)$. 
Then 
\begin{align}
&\int_{0}^{q_0} \yt ^{1 / 2 - 1} d\yt = 2 {q_0}^{1 / 2} \text{.} \non 
\end{align}
Therefore, 
\begin{align}
&p_0 = \int_{0}^{q_0} \yt ^{1 / 2 - 1} d\yt / \int_{0}^{s / 2} \yt ^{1 / 2 - 1} d\yt \non 
\end{align}
if and only if 
\begin{align}
&q_0 = \{ (s / 2)^{1 / 2} p_0 \} ^2 \text{.} \non 
\end{align}
\end{remark}

\subsubsection{The case of $s > 1$ and $r \le 1 / 2$}
In this case, we have to deal with something like the exponential integral when considering the right tail of (\ref{eq:prop_new_2}). 
However, we can use results of \cite{d2021}. 
In order to use them, we make a change of varaiables at the beginning. 
If $Y$ is distributed with density $g(y | r, s)$, $y \in (0, \infty )$, then $Y^{\dag } = \log \{ 1 - (1 - r) Y / s \} - \log r$ is distributed with density 
\begin{align}
g^{\dag } ( \yd | r, s) &= \frac{ \displaystyle {1(0 < \yd < \log (1 / r)) \over (1 - r e^{\yd } )^{1 / 2}} \exp \Big( - {r s e^{\yd } \over 1 - r} \Big) }{ \displaystyle \int_{0}^{\log (1 / r)} {1 \over (1 - r e^{\ytd } )^{1 / 2}} \exp \Big( - {r s e^{\ytd } \over 1 - r} \Big) d\ytd } \text{,} \quad \ytd \in (0, \infty ) \text{.} \non 
\end{align}
Fix $\ytd \in (0, \infty )$. 

First, 
\begin{align}
&1 \Big( \log {1 + r \over 2 r} < \yd < \log {1 \over r} \Big) g^{\dag } ( \yd | r, s) \non \\
&\le \frac{ \displaystyle 1 \Big( \log {1 + r \over 2 r} < \yd < \log {1 \over r} \Big) {1 \over (1 - r e^{\yd } )^{1 / 2}} \exp \Big( - {s \over 2} {1 + r \over 1 - r} \Big) }{ \displaystyle \int_{\log \{ (1 + r) / (2 r) \} }^{\log (1 / r)} {1 \over (1 - r e^{\ytd } )^{1 / 2}} d\ytd } \frac{ \displaystyle \int_{\log \{ (1 + r) / (2 r) \} }^{\log (1 / r)} {1 \over (1 - r e^{\ytd } )^{1 / 2}} d\ytd }{ \displaystyle \int_{0}^{\log (1 / r)} {1 \over (1 - r e^{\ytd } )^{1 / 2}} \exp \Big( - {r s e^{\ytd } \over 1 - r} \Big) d\ytd } \text{.} \non 
\end{align}
Second, 
\begin{align}
1 \Big( 0 < \yd < \log {1 + r \over 2 r} \Big) g^{\dag } ( \yd | r, s) &\le \frac{ \displaystyle 1 \Big( 0 < \yd < \log {1 + r \over 2 r} \Big) {2^{1 / 2} \over (1 - r)^{1 / 2}} \exp \Big( - {r s e^{\yd } \over 1 - r} \Big) }{ \displaystyle \int_{0}^{\log (1 / r)} {1 \over (1 - r e^{\ytd } )^{1 / 2}} \exp \Big( - {r s e^{\ytd } \over 1 - r} \Big) d\ytd } \non \\
&= \frac{ \displaystyle 1 \Big( 0 < \yd < \log {T \over S} \Big) {2^{1 / 2} \over (1 - r)^{1 / 2}} \exp (- S e^{\yd } ) e^S e^{- S} }{ \displaystyle \int_{0}^{\log (1 / r)} {1 \over (1 - r e^{\ytd } )^{1 / 2}} \exp \Big( - {r s e^{\ytd } \over 1 - r} \Big) d\ytd } \non \\
&= \frac{ \displaystyle 1 \Big( 0 < \yd < \log {T \over S} \Big) {2^{1 / 2} \over (1 - r)^{1 / 2}} \exp \{ H(y; S, T) \} e^{- S} }{ \displaystyle \int_{0}^{\log (1 / r)} {1 \over (1 - r e^{\ytd } )^{1 / 2}} \exp \Big( - {r s e^{\ytd } \over 1 - r} \Big) d\ytd } \text{,} \non 
\end{align}
where 
\begin{align}
&S = {r s \over 1 - r} \quad \text{and} \quad T = {r s \over 1 - r} {1 + r \over 2 r} \non 
\end{align}
and where 
\begin{align}
&H(y; S, T) = - S ( e^y - 1) \non 
\end{align}
for $y \in (0, \log (T / S))$. 
Let 
\begin{align}
&Z_0 = \log (T / S) \text{,} \quad Z_1 = \log \{ 1 + 1 / (2 S) \} \text{,} \quad Z_2 = \infty \text{,} \non 
\end{align}
let 
\begin{align}
&Z = \min \{ Z_0 , Z_1 , Z_2 \} \text{,} \non 
\end{align}
and let 
\begin{align}
&A = S e^Z \text{.} \non 
\end{align}
Then, by an inequality in \cite{d2021}, 
\begin{align}
&\exp \{ H(y; S, T) \}\le \exp \{ G(y; S, T) \} \non 
\end{align}
for all $y \in (0, \log (T / S))$, where 
\begin{align}
&G(y; S, T) = \begin{cases} 0 \text{,} & \text{if $y < Z$} \text{,} \\ H(Z; S, T) - A (y - Z) \text{,} & \text{if $y \ge Z$} \text{,} \end{cases} \non 
\end{align}
for $y \in (0, \log (T / S))$. 
Therefore, 
\begin{align}
&1 \Big( 0 < \yd < \log {1 + r \over 2 r} \Big) g^{\dag } ( \yd | r, s) \non \\
&\le \frac{ \displaystyle 1 \Big( 0 < \yd < \log {T \over S} \Big) {2^{1 / 2} \over (1 - r)^{1 / 2}} \exp \{ G(y; S, T) \} e^{- S} }{ \displaystyle \int_{0}^{\log (T / S)} \exp \{ G( \ytd ; S, T) \} d\ytd } \frac{ \displaystyle \int_{0}^{\log (T / S)} \exp \{ G( \ytd ; S, T) \} d\ytd }{ \displaystyle \int_{0}^{\log (1 / r)} {1 \over (1 - r e^{\ytd } )^{1 / 2}} \exp \Big( - {r s e^{\ytd } \over 1 - r} \Big) d\ytd } \text{,} \non 
\end{align}

\begin{remark}
Fix $q_0 \in ( \log \{ (1 + r) / (2 r) \} , \log (1 / r))$ and $p_0 \in (0, 1)$. 
Then 
\begin{align}
&\int_{\log \{ (1 + r) / (2 r) \} }^{q_0} {1 \over (1 - r e^{\ytd } )^{1 / 2}} d\ytd = 2 \log {\{ (1 - r) / (1 + r) \} ^{1 / 2} + \{ 2 / (1 + r) \} ^{1 / 2} \over \{ 1 / (r e^{q_0} ) - 1 \} ^{1 / 2} + \{ 1 / (r e^{q_0} ) \} ^{1 / 2}} \text{.} \non 
\end{align}
Therefore, 
\begin{align}
&p_0 = \int_{\log \{ (1 + r) / (2 r) \} }^{q_0} {1 \over (1 - r e^{\ytd } )^{1 / 2}} d\ytd / \int_{\log \{ (1 + r) / (2 r) \} }^{\log (1 / r)} {1 \over (1 - r e^{\ytd } )^{1 / 2}} d\ytd \non 
\end{align}
if and only if 
\begin{align}
&q_0 = \log {4 / r \over ([ \{ (1 - r)^{1 / 2} + 2^{1 / 2} \} / (1 + r)^{1 / 2} ]^{1 - p_0} + 1 / [ \{ (1 - r)^{1 / 2} + 2^{1 / 2} \} / (1 + r)^{1 / 2} ]^{1 - p_0} )^2} \text{.} \non 
\end{align}
\end{remark}

\begin{remark}
Fix $q_0 \in (0, \log (T / S))$ and $p_0 \in (0, 1)$. 
Then 
\begin{align}
&\int_{0}^{q_0} \exp \{ G( \ytd ; S, T) \} d\ytd = q_0 \wedge Z + 1( q_0 > Z) \int_{Z}^{q_0} \exp \{ - S ( e^Z - 1) - A (y - Z) \} d\ytd \text{.} \non 
\end{align}
Therefore, 
\begin{align}
&p_0 = \int_{0}^{q_0} \exp \{ G( \ytd ; S, T) \} d\ytd / \int_{0}^{\log (T / S)} \exp \{ G( \ytd ; S, T) \} d\ytd \non 
\end{align}
if and only if 
\begin{align}
&\Big[ p_0 \le \frac{ Z }{ \displaystyle Z + [ \exp \{ - S ( e^Z - 1) \} / A] (1 - \exp [- A \{ \log (T / S) - Z \} ]) } \quad \text{and} \non \\
&\quad q_0 = p_0 \{ Z + [ \exp \{ - S ( e^Z - 1) \} / A] (1 - \exp [- A \{ \log (T / S) - Z \} ]) \} \Big] \non \\
&\quad \text{or} \quad \Big( p_0 > \frac{ Z }{ \displaystyle Z + [ \exp \{ - S ( e^Z - 1) \} / A] (1 - \exp [- A \{ \log (T / S) - Z \} ]) } \quad \text{and} \non \\
&\quad q_0 = Z - {1 \over A} \log \Big[ 1 - {p_0 \{ Z + [ \exp \{ - S ( e^Z - 1) \} / A] (1 - \exp [- A \{ \log (T / S) - Z \} ]) \} - Z \over \exp \{ - S ( e^Z - 1) \} / A} \Big] \Big) \text{.} \non 
\end{align}
\end{remark}

\section{Proof of Proposition \ref{prp:new}}
\label{sec:proof_general_prp_new} 
We show that the rejection constants of the previous section are bounded. 
\subsection{The case of $s \le 1$}
\subsubsection{The case of $s \le 1$ and $r > 1$}
We have 
\begin{align}
&{e \over r^{1 / 2}} \frac{ \displaystyle \int_{- \infty }^{\infty } {1 \over 1 + \xt ^2} d\xt }{ \displaystyle \int_{- \infty }^{\infty } {1 \over 1 + \xt ^2} {1 \over (r + \xt ^2 )^{1 / 2}} \exp \Big( s {\xt ^2 \over r + \xt ^2} \Big) d\xt } \le {e \over r^{1 / 2}} \frac{ \displaystyle \int_{- \infty }^{\infty } {1 \over 1 + \xt ^2} d\xt }{ \displaystyle \int_{- 1}^{1} {1 \over 1 + \xt ^2} {1 \over (r + 1)^{1 / 2}} d\xt } \le e 2^{3 / 2} \text{.} \non 
\end{align}

\subsubsection{The case of $s \le 1$ and $r \le 1$}
We have 
\begin{align}
&e \frac{ \displaystyle \int_{\mathbb{R} \setminus [- 1, 1]} {1 \over 1 + \xt ^2} d\xt }{ \displaystyle \int_{- \infty }^{\infty } {1 \over 1 + \xt ^2} {1 \over (r + \xt ^2 )^{1 / 2}} \exp \Big( s {\xt ^2 \over r + \xt ^2} \Big) d\xt } \le e \frac{ \displaystyle \int_{\mathbb{R} \setminus [- 1, 1]} {1 \over 1 + \xt ^2} d\xt }{ \displaystyle \int_{- \infty }^{\infty } {1 \over (1 + \xt ^2 )^{3 / 2}} d\xt } = {e \pi \over 4} \non 
\end{align}
and 
\begin{align}
&e \frac{ \displaystyle \int_{- 1}^{1} {1 \over (r + \xt ^2 )^{1 / 2}} d\xt }{ \displaystyle \int_{- \infty }^{\infty } {1 \over 1 + \xt ^2} {1 \over (r + \xt ^2 )^{1 / 2}} \exp \Big( s {\xt ^2 \over r + \xt ^2} \Big) d\xt } \le e \frac{ \displaystyle \int_{- 1}^{1} {1 \over (r + \xt ^2 )^{1 / 2}} d\xt }{ \displaystyle \int_{- 1}^{1} {1 \over 2} {1 \over (r + \xt ^2 )^{1 / 2}} d\xt } = 2 e \text{.} \non 
\end{align}

\subsection{The case of $s > 1$}
\subsubsection{The case of $s > 1$ and $r > 1$}
We have that 
\begin{align}
&e \frac{ \displaystyle \int_{0}^{1} {\yt ^{1 / 2 - 1} \over 1 + (r - 1) \yt / s} d\yt }{ \displaystyle \int_{0}^{s} {\yt ^{1 / 2 - 1} e^{\yt } \over 1 + (r - 1) \yt / s} d\yt } \le e \frac{ \displaystyle \int_{0}^{1} {\yt ^{1 / 2 - 1} \over 1 + (r - 1) \yt / s} d\yt }{ \displaystyle \int_{0}^{1} {\yt ^{1 / 2 - 1} \over 1 + (r - 1) \yt / s} d\yt } = e \text{,} \non 
\end{align}
that 
\begin{align}
&{(s / 2)^{1 / 2 - 1} \over 1 + (r - 1) / 2} \frac{ \displaystyle \int_{s / 2}^{s} e^{\yt } d\yt }{ \displaystyle \int_{0}^{s} {\yt ^{1 / 2 - 1} e^{\yt } \over 1 + (r - 1) \yt / s} d\yt } \le {(s / 2)^{1 / 2 - 1} \over 1 + (r - 1) / 2} \frac{ \displaystyle \int_{s / 2}^{s} e^{\yt } d\yt }{ \displaystyle \int_{0}^{s} {s^{1 / 2 - 1} e^{\yt } \over r} d\yt } \le {2^{1 / 2} r \over 1 + (r - 1) / 2} \text{,} \non 
\end{align}
and that 
\begin{align}
{e^{s / 2} \over 1 + (r - 1) / s} \frac{ \displaystyle \int_{1}^{s / 2} \yt ^{1 / 2 - 1} d\yt }{ \displaystyle \int_{0}^{s} {\yt ^{1 / 2 - 1} e^{\yt } \over 1 + (r - 1) \yt / s} d\yt } &\le {e^{s / 2} \over 1 / s + (r - 1) / s} \frac{ \displaystyle \int_{1}^{s / 2} \yt ^{1 / 2 - 1} d\yt }{ \displaystyle \int_{(2 / 3) s}^{s} {\yt ^{1 / 2 - 1} e^{(2 / 3) s} \over r} d\yt } \non \\
&= {s \over e^{s / 6}} \frac{ 2 \{ (s / 2)^{1 / 2} - 1 \} }{ 2 \{ s^{1 / 2} - (2 / 3)^{1 / 2} s^{1 / 2} \} } \le {s \over e^{s / 6}} \frac{ (1 / 2)^{1 / 2} }{ 1 - (2 / 3)^{1 / 2} } \text{.} \non 
\end{align}

\subsubsection{The case of $s > 1$ and $1 / 2 < r \le 1$}
We have 
\begin{align}
&{e^{s / 2} \over 1 - (1 - r) / 2} \frac{ \displaystyle \int_{0}^{s / 2} \yt ^{1 / 2 - 1} d\yt }{ \displaystyle \int_{0}^{s} {\yt ^{1 / 2 - 1} e^{\yt } \over 1 - (1 - r) \yt / s} d\yt } \le {e^{s / 2} \over 1 - (1 - r) / 2} \frac{ \displaystyle \int_{0}^{s / 2} \yt ^{1 / 2 - 1} d\yt }{ \displaystyle \int_{s / 2}^{s} {\yt ^{1 / 2 - 1} e^{s / 2} \over 1 - (1 - r) / 2} d\yt } = {(s / 2)^{1 / 2} \over s^{1 / 2} - (s / 2)^{1 / 2}} \non 
\end{align}
and 
\begin{align}
&{(s / 2)^{1 / 2 - 1} \over r} \frac{ \displaystyle \int_{s / 2}^{s} e^{\yt } d\yt }{ \displaystyle \int_{0}^{s} {\yt ^{1 / 2 - 1} e^{\yt } \over 1 - (1 - r) \yt / s} d\yt } \le {(s / 2)^{1 / 2 - 1} \over r} \frac{ \displaystyle \int_{s / 2}^{s} e^{\yt } d\yt }{ \displaystyle \int_{s / 2}^{s} {s^{1 / 2 - 1} e^{\yt } \over 1 - (1 - r) / 2} d\yt } = {(1 / 2)^{1 / 2 - 1} \over r} {1 + r \over 2} \text{.} \non 
\end{align}

\subsubsection{The case of $s > 1$ and $r \le 1 / 2$}
First, let $\phi = (1 + 1 / r) / 2$. 
Then $1 < \phi < 1 / r$ and 
\begin{align}
&\exp \Big( - {s \over 2} {1 + r \over 1 - r} \Big) \frac{ \displaystyle \int_{\log \{ (1 + r) / (2 r) \} }^{\log (1 / r)} {1 \over (1 - r e^{\ytd } )^{1 / 2}} d\ytd }{ \displaystyle \int_{0}^{\log (1 / r)} {1 \over (1 - r e^{\ytd } )^{1 / 2}} \exp \Big( - {r s e^{\ytd } \over 1 - r} \Big) d\ytd } \non \\
&\le \exp \Big( - {s \over 2} {1 + r \over 1 - r} \Big) \exp \Big( {r s \phi \over 1 - r} \Big) \frac{ \displaystyle \int_{\log \{ (1 + r) / (2 r) \} }^{\log (1 / r)} {1 \over (1 - r e^{\ytd } )^{1 / 2}} d\ytd }{ \displaystyle \int_{0}^{\log \phi } {1 \over (1 - r e^{\ytd } )^{1 / 2}} d\ytd } \text{.} \non 
\end{align}
Note that for all $\eta > 0$, we have that 
\begin{align}
\int_{0}^{\eta } {1 \over (1 - r e^{\ytd } )^{1 / 2}} d\ytd %
&= \int_{r}^{r e^{\eta }} {1 \over (1 - \tilde{\yt } ^{\dag } )^{1 / 2} \tilde{\yt } ^{\dag }} d{\tilde{\yt } ^{\dag }} %
= 2 \Big[ \log \{ \yh + (1 + \yh ^2 )^{1 / 2} \} \Big] _{\yh = \{ 1 / (r e^{\eta } ) - 1 \} ^{1 / 2}}^{\yh = (1 / r - 1)^{1 / 2}} \non 
\end{align}
if $\eta < \log (1 / r)$. 
Then 
\begin{align}
&\int_{\log \{ (1 + r) / (2 r) \} }^{\log (1 / r)} {1 \over (1 - r e^{\ytd } )^{1 / 2}} d\ytd / \int_{0}^{\log \phi } {1 \over (1 - r e^{\ytd } )^{1 / 2}} d\ytd \non \\
&= 2 \Big[ \log \{ \yh + (1 + \yh ^2 )^{1 / 2} \} \Big] _{\yh = 0}^{\yh = \{ (1 - r) / (1 + r) \} ^{1 / 2}} / \int_{r}^{r \phi } {1 \over (1 - \tilde{\yt } ^{\dag } )^{1 / 2} \tilde{\yt } ^{\dag }} d{\tilde{\yt } ^{\dag }} \non \\
&\le 2 \log (1 + 2^{1 / 2} ) / \int_{r}^{r \phi } {1 \over (1 - r)^{1 / 2} r \phi } d{\tilde{\yt } ^{\dag }} = 2 \log (1 + 2^{1 / 2} ) / {\phi - 1 \over \phi } \text{.} \non 
\end{align}
Therefore, 
\begin{align}
&\exp \Big( - {s \over 2} {1 + r \over 1 - r} \Big) \frac{ \displaystyle \int_{\log \{ (1 + r) / (2 r) \} }^{\log (1 / r)} {1 \over (1 - r e^{\ytd } )^{1 / 2}} d\ytd }{ \displaystyle \int_{0}^{\log (1 / r)} {1 \over (1 - r e^{\ytd } )^{1 / 2}} \exp \Big( - {r s e^{\ytd } \over 1 - r} \Big) d\ytd } \non \\
&\le \exp \Big( - {s \over 2} {1 + r \over 1 - r} \Big) \exp \Big( {r s \phi \over 1 - r} \Big) 2 \log (1 + 2^{1 / 2} ) / {\phi - 1 \over \phi } \non \\
&\le 2 \{ \log (1 + 2^{1 / 2} ) \} 3 \text{.} \non 
\end{align}

Second, 
\begin{align}
&{2^{1 / 2} \over (1 - r)^{1 / 2}} e^{- S} \frac{ \displaystyle \int_{0}^{\log (T / S)} \exp \{ G( \ytd ; S, T) \} d\ytd }{ \displaystyle \int_{0}^{\log (1 / r)} {1 \over (1 - r e^{\ytd } )^{1 / 2}} \exp \Big( - {r s e^{\ytd } \over 1 - r} \Big) d\ytd } \non \\
&\le {2^{1 / 2} \over (1 - r)^{1 / 2}} e^{- S} \frac{ \displaystyle \int_{0}^{\log (T / S)} \exp \{ G( \ytd ; S, T) \} d\ytd }{ \displaystyle \int_{0}^{\log (T / S)} {1 \over (1 - r)^{1 / 2}} \exp (- S e^{\ytd } ) d\ytd } \non \\
&= 2^{1 / 2} \frac{ \displaystyle \int_{0}^{\log (T / S)} \exp \{ G( \ytd ; S, T) \} d\ytd }{ \displaystyle \int_{0}^{\log (T / S)} \exp \{ H( \ytd ; S, T) \} d\ytd } \text{.} \non 
\end{align}
As in the proof of Theorem 1 of \cite{d2021}, it can be seen that 
\begin{align}
&\int_{0}^{\log (T / S)} \exp \{ G( \ytd ; S, T) \} d\ytd / \int_{0}^{\log (T / S)} \exp \{ H( \ytd ; S, T) \} d\ytd \non \\
&\le %
\exp \{ - H(Z; S, T) \} + 1( Z_1 < Z_0 ) 2 / [(1 + 2 S) \log \{ 1 + 1 / (2 S) \} ] \non \\
&\le e^{1 / 2} + 2 \text{.} \non 
\end{align}
Thus, 
\begin{align}
&{2^{1 / 2} \over (1 - r)^{1 / 2}} e^{- S} \frac{ \displaystyle \int_{0}^{\log (T / S)} \exp \{ G( \ytd ; S, T) \} d\ytd }{ \displaystyle \int_{0}^{\log (1 / r)} {1 \over (1 - r e^{\ytd } )^{1 / 2}} \exp \Big( - {r s e^{\ytd } \over 1 - r} \Big) d\ytd } \le 2^{1 / 2} ( e^{1 / 2} + 2) \text{.} \non 
\end{align}

\section{A Naive Approach to Sampling $\tas $}
\label{sec:naive_rejection} 
The conditional distribution of $\tas $ in Section 2.1 of \cite{bkp2022} is 
\begin{align}
p( \tas | \bla , \y ) &\propto {1 / \sqrt{| \I ^{(p)} + \tas \X ^{\top } \X \bLa ^2 |} \over ( \y ^{\top } [ \I ^{(n)} - \X \{ \X ^{\top } \X + ( \tas \bLa ^2 )^{- 1} \} ^{- 1} \X ^{\top } ] \y / 2 + b' )^{a' + n / 2}} {( \tas )^{1 / 2 - 1} \over 1 + \tas } \non 
\end{align}
in the horseshoe case. 
By making the change of variables $\ka = ( \tas )^{1 / 2} / \{ 1 + ( \tas )^{1 / 2} \} \in (0, 1)$, 
\begin{align}
&p( \ka | \bla , \y ) \non \\
&\propto %
\Big\{ {1 / \sqrt{| \I ^{(p)} + \tas \X ^{\top } \X \bLa ^2 |} \over ( \y ^{\top } [ \I ^{(n)} - \X \{ \X ^{\top } \X + ( \tas \bLa ^2 )^{- 1} \} ^{- 1} \X ^{\top } ] \y / 2 + b' )^{a' + n / 2}} \Big| _{\ta = \ka / (1 - \ka )} \Big\} {1 \over (1 - \ka )^2 + \ka ^2} \text{,} \non 
\end{align}
which is bounded. 
Furthermore, 
\begin{align}
&{1 / \sqrt{| \I ^{(p)} + \tas \X ^{\top } \X \bLa ^2 |}} \non 
\end{align}
is decreasing in $\ta $, 
\begin{align}
&{1 \over ( \y ^{\top } [ \I ^{(n)} - \X \{ \X ^{\top } \X + ( \tas \bLa ^2 )^{- 1} \} ^{- 1} \X ^{\top } ] \y / 2 + b' )^{a' + n / 2}} \non 
\end{align}
is increasing in $\ta $, and 
\begin{align}
&{1 \over (1 - \ka )^2 + \ka ^2} \non 
\end{align}
is increasing in $\ka $ for $\ka < 1 / 2$ and decreasing in $\ka $ for $\ka > 1 / 2$. 
Therefore, we can sample $\ka $ using rejection sampling based on the following inequality. 

\begin{lem}
\label{lem:naive} 
Let $h \colon [0, 1] \to [0, \infty )$ be a bounded function. 
Then for all $G \in \mathbb{N}$, we have that 
\begin{align}
\frac{ \displaystyle h( \ka ) }{ \displaystyle \int_{0}^{\infty } h( \tilde{\ka } ) d\tilde{\ka } } &\le %
\frac{ \displaystyle \sum_{g = 1}^{G} \overline{h} _g }{ \displaystyle G \int_{0}^{\infty } h( \tilde{\ka } ) d\tilde{\ka } } \sum_{g = 1}^{G} \frac{ \displaystyle \overline{h} _g }{ \displaystyle \sum_{\tilde{g} = 1}^{G} \overline{h} _{\tilde{g}} } {1((g - 1) / G \le \ka \le g / G) \over 1 / G} \non 
\end{align}
for all $\ka \in [0, 1]$, where 
\begin{align}
\overline{h} _g &\in [ \sup_{\tilde{\ka } \in [(g - 1) / G, g / G]} h( \tilde{\ka } ), \infty ) \non 
\end{align}
for $g = 1, \dots , G$. 
\end{lem}

\section{Rejection Sampling for the Conditional Distribution of $\bla ^2$ in the Unmodified JOB Sampling Scheme}
\label{sec:rejecton_expint} 
In the sampling scheme of \cite{job2020}, $\bla ^2 = ( \la _{k}^{2} )_{k = 1}^{p}$ is sampled from its conditional distribution without introducing the latent variables $\bnu $, as described on page 7 of \cite{bkp2022}. 
For this purpose, we use the method of \cite{d2021} instead of that of Appendix S1 of \cite{job2020}. 

Specifically, the conditional distribution of $\bla ^2$ is 
\begin{align}
p( \bla ^2 | \bbe , \sis , \tas , \y ) &\propto %
\prod_{k = 1}^{p} \Big\{ {1 \over ( \la _{k}^{2} )^2} {1 \over 1 + 1 / \lat _{k}^{2}} \exp \Big( - {{\be _k}^2 \over 2 \sis \tas \la _{k}^{2}} \Big) \Big\} \text{.} \non 
\end{align}
By making the change of variables $\bpsi ^{\dag } = ( \psi _{k}^{\dag } )_{k = 1}^{p} = %
\{ {\be _k}^2 / (2 \sis \tas ) \} (1 + 1 / \la _{k}^{2} )_{k = 1}^{p}$, 
\begin{align}
p( \bpsi ^{\dag } | \bbe , \sis , \tas , \y ) &\propto \prod_{k = 1}^{p} \Big\{ 1( \psi _{k}^{\dag } > {\be _k}^2 / (2 \sis \tas )) {\exp (- \psi _{k}^{\dag } ) \over \psi _{k}^{\dag }} \Big\} \text{.} \non 
\end{align}
Thus, we can use the method of Section 4 of \cite{d2021}.

\section{Proof of Theorem \ref{thm:lat}}
\label{sec:proof_thm_first} 
We use the drift and minorization technique to prove Theorem \ref{thm:lat}; see %
\cite{jh2001} for this technique. 
We suppress the dependence on $\y $ in the remainder of this Supplementary Material. 

For $V \colon (0, \infty )^p \to [0, \infty )$ and $\v \in (0, \infty )^p$, we write 
\begin{align}
(P V) ( \v ) &= E^{( \bbe , \sis , \bnu , \tas ) | {\blat {}^2}} [ E^{{\blat {}^2} | ( \bbe , \sis , \bnu , \tas )} [ V( {\blat {}^2} ) | \bbe , \sis , \bnu , \tas ] | {\blat {}^2} = \v ] \text{.} \non 
\end{align}

\begin{lem}
\label{lem:minorization} 
Let $L > 0$. 
Suppose that ${\blat {}^2} \in [1 / L, L]^p$. 
Then there exists a normalized probability density $q \colon \mathbb{R} ^p \times (0, \infty ) \times (0, \infty )^p \times (0, \infty ) \to (0, \infty )$ and a constant $\ep > 0$ such that 
\begin{align}
&p( \bbe , \sis , \bnu , \tas | {\blat {}^2} ) \ge \ep q( \bbe , \sis , \bnu , \tas ) \non 
\end{align}
for all $( \bbe , \sis , \bnu , \tas  ) \in \mathbb{R} ^p \times (0, \infty ) \times (0, \infty )^p \times (0, \infty )$. 
\end{lem}

\noindent
{\bf Proof%
.} \ \ We have 
\begin{align}
&p( \bbe , \sis , \bnu , \tas | {\blat {}^2} ) \non \\
&= {\rm{IG}} ( \sis | n / 2 + a' , \{ \| \y \| ^2 - \y ^{\top } \X ( \X ^{\top } \X + \bPsi )^{- 1} \X ^{\top } \y \} / 2 + b' ) \non \\
&\quad \times {\rm{N}}_p ( \bbe | ( \X ^{\top } \X + \bPsi )^{- 1} \X ^{\top } \y , \sis ( \X ^{\top } \X + \bPsi )^{- 1} ) \non \\
&\quad \times \Big( \prod_{k = 1}^{p} \Big[ \Big( 1 + {\tas \over \lat _{k}^{2}} \Big) {1 \over {\nu _k}^2} \exp \Big\{ - {1 \over \nu _k} \Big( 1 + {\tas \over \lat _{k}^{2}} \Big) \Big\} \Big] \Big) \non \\
&\quad \times \pi _{\ta } ( \tas ) \Big\{ \prod_{k = 1}^{p} {( \tas )^{1 / 2} \over \lat _{k}^{2} + \tas } \Big\} / \int_{0}^{\infty } \pi _{\ta } (t) \Big( \prod_{k = 1}^{p} {t^{1 / 2} \over \lat _{k}^{2} + t} \Big) dt \non \\
&= {[ \{ \| \y \| ^2 - \y ^{\top } \X ( \X ^{\top } \X + \bPsi )^{- 1} \X ^{\top } \y \} / 2 + b' ]^{n / 2 + a'} \over \Ga (n / 2 + a' )} \non \\
&\quad \times {1 \over ( \sis )^{1 + n / 2 + a'}} \exp \Big( - {1 \over \sis } [ \{ \| \y \| ^2 - \y ^{\top } \X ( \X ^{\top } \X + \bPsi )^{- 1} \X ^{\top } \y \} / 2 + b' ] \Big) \non \\
&\quad \times {| \X ^{\top } \X + \bPsi |^{1 / 2} \over (2 \pi )^{p / 2} ( \sis )^{p / 2}} \exp \Big[ - {\{ \bbe - ( \X ^{\top } \X + \bPsi )^{- 1} \X ^{\top } \y \} ^{\top } ( \X ^{\top } \X + \bPsi ) \{ \bbe - ( \X ^{\top } \X + \bPsi )^{- 1} \X ^{\top } \y \} \over 2 \sis } \Big] \non \\
&\quad \times \pi _{\ta } ( \tas ) \Big( \prod_{k = 1}^{p} \Big[ {( \tas )^{1 / 2} \over {\nu _k}^2} \exp \Big\{ - {1 \over \nu _k} \Big( 1 + {\tas \over \lat _{k}^{2}} \Big) \Big\} \Big] \Big) / \int_{0}^{\infty } \pi _{\ta } (t) \Big( \prod_{k = 1}^{p} {t^{1 / 2} \over 1 + t / \lat _{k}^{2}} \Big) dt \non \\
&\ge {( b' )^{n / 2 + a'} \over \Ga (n / 2 + a' )} {1 \over ( \sis )^{1 + n / 2 + a'}} \exp \Big\{ - {1 \over \sis } ( \| \y \| ^2 / 2 + b' ) \Big\} \non \\
&\quad \times {| \X ^{\top } \X |^{1 / 2} \over (2 \pi )^{p / 2}} {1 \over ( \sis )^{p / 2}} \exp \Big[ - {2 \{ \bbe ^{\top } ( \X ^{\top } \X + \bPsi ) \bbe + \y ^{\top } \X ( \X ^{\top } \X + \bPsi )^{- 1} \X ^{\top } \y \} \over 2 \sis } \Big] \non \\
&\quad \times \pi _{\ta } ( \tas ) \Big( \prod_{k = 1}^{p} \Big[ {( \tas )^{1 / 2} \over {\nu _k}^2} \exp \Big\{ - {1 \over \nu _k} \Big( 1 + {\tas \over \lat _{k}^{2}} \Big) \Big\} \Big] \Big) / \int_{0}^{\infty } \pi _{\ta } (t) \Big( \prod_{k = 1}^{p} {t^{1 / 2} \over 1 + t / \lat _{k}^{2}} \Big) dt \text{.} \non 
\end{align}
Since $\bPsi \le L \I ^{(p)}$ and $\lat _{k}^{2} \le L$ for all $k = 1, \dots , p$ by assumption, 
\begin{align}
&p( \bbe , \sis , \bnu , \tas | {\blat {}^2} ) \non \\
&\ge {( b' )^{n / 2 + a'} \over \Ga (n / 2 + a' )} {1 \over ( \sis )^{1 + n / 2 + a'}} \exp \Big\{ - {1 \over \sis } ( \| \y \| ^2 / 2 + b' ) \Big\} \non \\
&\quad \times {| \X ^{\top } \X |^{1 / 2} \over (2 \pi )^{p / 2}} {1 \over ( \sis )^{p / 2}} \exp \Big[ - {2 \{ \bbe ^{\top } ( \X ^{\top } \X + L \I ^{(p)} ) \bbe + \y ^{\top } \X ( \X ^{\top } \X )^{- 1} \X ^{\top } \y \} \over 2 \sis } \Big] \non \\
&\quad \times \pi _{\ta } ( \tas ) \Big( \prod_{k = 1}^{p} \Big[ {( \tas )^{1 / 2} \over {\nu _k}^2} \exp \Big\{ - {1 \over \nu _k} (1 + L \tas ) \Big\} \Big] \Big) / \int_{0}^{\infty } \pi _{\ta } (t) \Big( \prod_{k = 1}^{p} {t^{1 / 2} \over 1 + t / L} \Big) dt \text{.} \non 
\end{align}
The right-hand side is an integrable function of $( \bbe , \sis , \bnu , \tas )$ which does not depend on ${\blat {}^2}$ and this completes the proof. 
\hfill$\Box$

\bigskip

For $\ep > 0$, define the function $V_{1, \ep } \colon (0, \infty )^p \to [0, \infty )$ by 
\begin{align}
V_{1, \ep } (( v_k )_{k = 1}^{p} ) &= \sum_{k = 1}^{p} {v_k}^{\ep } \text{,} \quad ( v_k )_{k = 1}^{p} \in (0, \infty )^p \text{.} \non 
\end{align}

\begin{lem}
\label{lem:large} 
There exists $\ep > 0$ such that 
\begin{align}
(P V_{1, \ep } ) ( \v ) \le \De + \de V_{1, \ep } ( \v ) \non 
\end{align}
for all $\v \in (0, \infty )^p$ for some $0 < \de < 1$ and $\De > 0$. 
\end{lem}

\noindent
{\bf Proof%
.} \ \ Fix $\ep > 0$. 
For all $k = 1, \dots , p$, by making the change of variables $\psi = 1 / \lat _{k}^{2}$, 
\begin{align}
E[ ( \lat _{k}^{2} )^{\ep } | \bbe , \sis , \bnu , \tas ] &= \int_{0}^{\infty } \Big( {{\be _k}^2 \over 2 \sis } + {\tas \over \nu _k} \Big) \psi ^{1 - \ep - 1} \exp \Big\{ - \psi \Big({{\be _k}^2 \over 2 \sis }  + {\tas \over \nu _k} \Big) \Big\} d\psi \non \\
&= \Ga (1 - \ep ) \Big( {{\be _k}^2 \over 2 \sis } + {\tas \over \nu _k} \Big) ^{\ep } \non \\
&\le \Ga (1 - \ep ) \Big( {{\be _k}^2 \over 2 \sis } \Big) ^{\ep } + \Ga (1 - \ep ) \Big( {\tas \over \nu _k} \Big) ^{\ep } \text{.} \label{llargep1} 
\end{align}
This is as in the case considered on page 12 of %
\cite{bkp2022}. 

Fix $k = 1, \dots , p$. 
Note that 
\begin{align}
&E \Big[ {{\be _k}^2 \over 2 \sis } \Big| {\blat {}^2} \Big] \non \\
&= E \Big[ {1 \over 2 \sis } \{ \sis ( \X ^{\top } \X + \bPsi )^{k, k} + \y ^{\top } \X ( \X ^{\top } \X + \bPsi )^{- 1} \e _{k}^{(p)} ( \e _{k}^{(p)} )^{\top } ( \X ^{\top } \X + \bPsi )^{- 1} \X ^{\top } \y \} \Big| {\blat {}^2} \Big] \non \\
&\le E \Big[ {1 \over 2 \sis } \{ \sis ( \X ^{\top } \X + \bPsi )^{k, k} + \y ^{\top } \X ( \X ^{\top } \X + \bPsi )^{- 1} ( M_1 \X ^{\top } \X ) ( \X ^{\top } \X + \bPsi )^{- 1} \X ^{\top } \y \} \Big| {\blat {}^2} \Big] \non \\
&\le E \Big[ {1 \over 2 \sis } \{ \sis ( \X ^{\top } \X + \bPsi )^{k, k} + M_1 \y ^{\top } \X ( \X ^{\top } \X + \bPsi )^{- 1} \X ^{\top } \y \} \Big| {\blat {}^2} \Big] \non \\
&\le E \Big[ {1 \over 2 \sis } \{ \sis ( \X ^{\top } \X )^{k, k} + M_1 \y ^{\top } \X ( \X ^{\top } \X )^{- 1} \X ^{\top } \y \} \Big| {\blat {}^2} \Big] \non 
\end{align}
for some $M_1 > 0$. 
Then 
\begin{align}
&E \Big[ {{\be _k}^2 \over 2 \sis } \Big| {\blat {}^2} \Big] \non \\
&\le {1 \over 2} \Big[ ( \X ^{\top } \X )^{k, k} + M_1 \y ^{\top } \X ( \X ^{\top } \X )^{- 1} \X ^{\top } \y {n / 2 + a' \over \{ \| \y \| ^2 - \y ^{\top } \X ( \X ^{\top } \X + \bPsi )^{- 1} \X ^{\top } \y \} / 2 + b'} \Big] \non \\
&\le {1 \over 2} \Big\{ ( \X ^{\top } \X )^{k, k} + M_1 \y ^{\top } \X ( \X ^{\top } \X )^{- 1} \X ^{\top } \y {n / 2 + a' \over b'} \Big\} \text{.} \non 
\end{align}
Therefore, by Jensen's inequality, 
\begin{align}
&E \Big[ \Big( {{\be _k}^2 \over 2 \sis } \Big) ^{\ep } \Big| {\blat {}^2} \Big] \le M_2 \label{llargep2} 
\end{align}
for some constant $M_2 > 0$. 

Meanwhile, by making the change of variables $\mu _k = 1 / \nu _k$ for $k = 1, \dots , p$, 
\begin{align}
\sum_{k = 1}^{p} E \Big[ \Big( {\tas \over \nu _k} \Big) ^{\ep } \Big| {\blat {}^2} \Big] %
&= \sum_{k = 1}^{p} E \Big[ ( \tas )^{\ep } \int_{0}^{\infty } {\mu _k}^{\ep } \Big( 1 + {\tas \over \lat _{k}^{2}} \Big) \exp \Big\{ - \mu _k \Big( 1 + {\tas \over \lat _{k}^{2}} \Big) \Big\} d{\mu _k} \Big| {\blat {}^2} \Big] \non \\
&= \sum_{k = 1}^{p} E \Big[ ( \lat _{k}^{2} )^{\ep } \Ga (1 + \ep ) {( \tas )^{\ep } \over ( \lat _{k}^{2} + \tas )^{\ep }} \Big| {\blat {}^2} \Big]  \non \\
&\le \sum_{\substack{k = 1 \\ k \notin S( {\blat {}^2} )}}^{p} c^{\ep } \Ga (1 + \ep ) + \sum_{k \in S( {\blat {}^2} )} E \Big[ ( \lat _{k}^{2} )^{\ep } \Ga (1 + \ep ) {( \tas )^{\ep } \over ( \lat _{k}^{2} + \tas )^{\ep }} \Big| {\blat {}^2} \Big] \text{,} \non 
\end{align}
where $S( {\blat {}^2} ) = \{ k = 1, \dots , p | \lat _{k}^{2} \ge c \} $. 
By the covariance inequality %
(see Lemma 5.6.6 of \cite{lc1998}) 
and Jensen's inequality, 
\begin{align}
&\sum_{k = 1}^{p} E \Big[ \Big( {\tas \over \nu _k} \Big) ^{\ep } \Big| {\blat {}^2} \Big] - \sum_{\substack{k = 1 \\ k \notin S( {\blat {}^2} )}}^{p} c^{\ep } \Ga (1 + \ep ) \non \\
&\le |S( {\blat {}^2} )| E \Big[ \Big\{ {1 \over |S( {\blat {}^2} )|} \sum_{k \in S( {\blat {}^2} )} ( \lat _{k}^{2} )^{\ep } \Ga (1 + \ep ) \Big\} {1 \over |S( {\blat {}^2} )|} \sum_{k \in S( {\blat {}^2} )} {( \tas )^{\ep } \over ( \lat _{k}^{2} + \tas )^{\ep }} \Big| {\blat {}^2} \Big] \non \\
&\le {\Ga (1 + \ep ) \over |S( {\blat {}^2} )|} \Big\{ \sum_{k \in S( {\blat {}^2} )} ( \lat _{k}^{2} )^{\ep } \Big\} \sum_{k \in S( {\blat {}^2} )} \Big( E \Big[ {\tas \over \lat _{k}^{2} + \tas } \Big| {\blat {}^2} \Big] \Big) ^{\ep } \non \\
&\le \Ga (1 + \ep ) \Big\{ \sum_{k \in S( {\blat {}^2} )} ( \lat _{k}^{2} )^{\ep } \Big\} \Big\{ {1 \over |S( {\blat {}^2} )|} \sum_{k \in S( {\blat {}^2} )} E \Big[ {\tas \over \lat _{k}^{2} + \tas } \Big| {\blat {}^2} \Big] \Big\} ^{\ep } \text{.} \non 
\end{align}
Since by %
Lemma \ref{lem:K} 
\begin{align}
{p \over 2} + a &\ge E \Big[ \{ a + b + p - |S( {\blat {}^2} )| \} {\tas \over c + \tas } + \sum_{k \in S( {\blat {}^2} )} {\tas \over \lat _{k}^{2} + \tas } \Big| {\blat {}^2} \Big] \non \\
&\ge E \Big[ \sum_{k \in S( {\blat {}^2} )} {a + b + p \over |S( {\blat {}^2} )|} {\tas \over \lat _{k}^{2} + \tas } \Big| {\blat {}^2} \Big] \text{,} \non 
\end{align}
it follows that 
\begin{align}
&\sum_{k = 1}^{p} E \Big[ \Big( {\tas \over \nu _k} \Big) ^{\ep } \Big| {\blat {}^2} \Big] - \sum_{\substack{k = 1 \\ k \notin S( {\blat {}^2} )}}^{p} c^{\ep } \Ga (1 + \ep ) \le \Ga (1 + \ep ) \Big\{ \sum_{k \in S( {\blat {}^2} )} ( \lat _{k}^{2} )^{\ep } \Big\} \Big( {p / 2 + a \over a + b + p} \Big) ^{\ep } \text{.} \label{llargep3} 
\end{align}

By (\ref{llargep1}), (\ref{llargep2}), and (\ref{llargep3}), 
\begin{align}
(P V_{1, \ep } ) ( {\blat {}^2} ) &\le p \Ga (1 - \ep ) M_2 + \Ga (1 - \ep ) \Ga (1 + \ep ) \Big[ \{ p - |S( {\blat {}^2} )| \} c^{\ep } + \Big\{ \sum_{k \in S( {\blat {}^2} )} ( \lat _{k}^{2} )^{\ep } \Big\} \Big( {p / 2 + a \over a + b + p} \Big) ^{\ep } \Big] \text{.} \non 
\end{align}
Note that 
\begin{align}
&\Ga (1 - \ep ) \Ga (1 + \ep ) \Big( {p / 2 + a \over a + b + p} \Big) ^{\ep } = \Ga (1 - \ep ) \Ga ( \ep ) \ep \Big( {p / 2 + a \over a + b + p} \Big) ^{\ep } = {\pi \over \sin ( \pi \ep )} \ep \Big( {p / 2 + a \over a + b + p} \Big) ^{\ep } \gtreqless 1 \non 
\end{align}
if and only if 
\begin{align}
&z \rho ^z \gtreqless \sin z \text{,} \non 
\end{align}
where $z = \pi \ep $ and $\rho = \{ (p / 2 + a) / (a + b + p) \} ^{1 / \pi }$, and that 
\begin{align}
&{\pd \over \pd z} (z \rho ^z - \sin z) = \rho ^z + z \rho ^z \log \rho - \cos z \le \rho ^z - \cos z < 0 \non 
\end{align}
if $z > 0$ is sufficiently small. 
Then 
\begin{align}
&\Ga (1 - \ep ) \Ga (1 + \ep ) \Big( {p / 2 + a \over a + b + p} \Big) ^{\ep } < 1 \non 
\end{align}
for some sufficiently small $0 < \ep < 1$. 
This completes the proof. 
\hfill$\Box$

\bigskip

For $\ep > 0$, define the function $V_{2, \ep } \colon (0, \infty )^p \to [0, \infty )$ by 
\begin{align}
V_{2, \ep } (( v_k )_{k = 1}^{p} ) &= \sum_{k = 1}^{p} {1 \over {v_k}^{\ep }} \text{,} \quad ( v_k )_{k = 1}^{p} \in (0, \infty )^p \text{.} \non 
\end{align}

\begin{lem}
\label{lem:small} 
For all $0 < \ep < 1 / 2$, there exists $M > 0$ such that 
\begin{align}
(P V_{2, \ep } ) ( \v ) &\le M + \sum_{k = 1}^{p} {1 \over {v_k}^{\ep }} {\Ga (1 + \ep ) \over \Ga ( \ep )} \int_{0}^{\infty } {w^{\ep - 1} \over (1 + w)^{1 / 2}} E \Big[ {v_k + \tas \over v_k + (1 + w) \tas } \Big| {\blat {}^2} = \v \Big] dw \non 
\end{align}
for all $\v = ( v_k )_{k = 1}^{p} \in (0, \infty )^p$. 
\end{lem}

\noindent
{\bf Proof%
.} \ \ Fix $0 < \ep < 1 / 2$. 
Fix $k = 1, \dots , p$. 
Then, by making the change of variables $\psi = 1 / \lat _{k}^{2}$, 
\begin{align}
E \Big[ {1 \over ( \lat _{k}^{2} )^{\ep }} \Big| \bbe , \sis , \bnu , \tas \Big] &= \int_{0}^{\infty } \Big( {{\be _k}^2 \over 2 \sis } + {\tas \over \nu _k} \Big) \psi ^{1 + \ep - 1} \exp \Big\{ - \psi \Big( {{\be _k}^2 \over 2 \sis }  + {\tas \over \nu _k} \Big) \Big\} d\psi \non \\
&= \Ga (1 + \ep ) / \Big( {{\be _k}^2 \over 2 \sis } + {\tas \over \nu _k} \Big) ^{\ep } \non \\
&= \int_{0}^{\infty } {\Ga (1 + \ep ) \over \Ga ( \ep )} u^{\ep - 1} \exp \Big( - u {{\be _k}^2 \over 2 \sis } \Big) \exp \Big( - u {\tas \over \nu _k} \Big) du \text{.} \label{lsmallp1} 
\end{align}

Fix $u \in (0, \infty )$. 
Then 
\begin{align}
E \Big[ \exp \Big( - u {{\be _k}^2 \over 2 \sis } \Big) \Big| \sis , {\blat {}^2} \Big] %
&= \int_{- \infty }^{\infty } \Big\{ {1 \over (2 \pi )^{1 / 2}} {1 \over ( \sis )^{1 / 2}} {1 \over \{ ( \X ^{\top } \X + \bPsi )^{k, k} \} ^{1 / 2}} \non \\
&\quad \times \exp \Big( - {1 \over 2 \sis } \Big[ {\{ \be _k - ( \e _{k}^{(p)} )^{\top } ( \X ^{\top } \X + \bPsi )^{- 1} \X ^{\top } \y \} ^2 \over ( \X ^{\top } \X + \bPsi )^{k, k}} + u {\be _k}^2 \Big] \Big) \Big\} d{\be _k} \text{,} \non 
\end{align}
where 
\begin{align}
&{\{ \be _k - ( \e _{k}^{(p)} )^{\top } ( \X ^{\top } \X + \bPsi )^{- 1} \X ^{\top } \y \} ^2 \over ( \X ^{\top } \X + \bPsi )^{k, k}} + u {\be _k}^2 \non \\
&= \Big\{ {1 \over ( \X ^{\top } \X + \bPsi )^{k, k}} + u \Big\} \Big\{ \be _k - {( \e _{k}^{(p)} )^{\top } ( \X ^{\top } \X + \bPsi )^{- 1} \X ^{\top } \y / ( \X ^{\top } \X + \bPsi )^{k, k} \over 1 / ( \X ^{\top } \X + \bPsi )^{k, k} + u} \Big\} ^2 \non \\
&\quad + {u \{ ( \e _{k}^{(p)} )^{\top } ( \X ^{\top } \X + \bPsi )^{- 1} \X ^{\top } \y \} ^2 \over 1 + u ( \X ^{\top } \X + \bPsi )^{k, k}} \text{.} \non 
\end{align}
Therefore, 
\begin{align}
&E \Big[ \exp \Big( - u {{\be _k}^2 \over 2 \sis } \Big) \Big| {\blat {}^2} \Big] \non \\
&= E \Big[ {1 \over \{ ( \X ^{\top } \X + \bPsi )^{k, k} \} ^{1 / 2}} {1 \over \{ 1 / ( \X ^{\top } \X + \bPsi )^{k, k} + u \} ^{1 / 2}} \exp \Big[ - {1 \over 2 \sis } {u \{ ( \e _{k}^{(p)} )^{\top } ( \X ^{\top } \X + \bPsi )^{- 1} \X ^{\top } \y \} ^2 \over 1 + u ( \X ^{\top } \X + \bPsi )^{k, k}} \Big] \Big| {\blat {}^2} \Big] \non \\
&= \int_{0}^{\infty } \Big( {[ \{ \| \y \| ^2 - \y ^{\top } \X ( \X ^{\top } \X + \bPsi )^{- 1} \X ^{\top } \y \} / 2 + b' ]^{n / 2 + a'} \over \Ga (n / 2 + a' )} \non \\
&\quad \times {1 \over ( \sis )^{1 + n / 2 + a'}} \exp \Big[ - {\{ \| \y \| ^2 - \y ^{\top } \X ( \X ^{\top } \X + \bPsi )^{- 1} \X ^{\top } \y \} / 2 + b' \over \sis} \Big] \non \\
&\quad \times {1 \over \{ 1 + u ( \X ^{\top } \X + \bPsi )^{k, k} \} ^{1 / 2}} \exp \Big[ - {1 \over 2 \sis } {u \{ ( \e _{k}^{(p)} )^{\top } ( \X ^{\top } \X + \bPsi )^{- 1} \X ^{\top } \y \} ^2 \over 1 + u ( \X ^{\top } \X + \bPsi )^{k, k}} \Big] \Big) d\sis \non \\
&= {1 \over \{ 1 + u ( \X ^{\top } \X + \bPsi )^{k, k} \} ^{1 / 2}} \Big\{ {\| \y \| ^2 - \y ^{\top } \X ( \X ^{\top } \X + \bPsi )^{- 1} \X ^{\top } \y \over 2} + b' \Big\} ^{n / 2 + a'} \non \\
&\quad / \Big[ {u \{ ( \e _{k}^{(p)} )^{\top } ( \X ^{\top } \X + \bPsi )^{- 1} \X ^{\top } \y \} ^2 \over 2 \{ 1 + u ( \X ^{\top } \X + \bPsi )^{k, k} \} } + {\| \y \| ^2 - \y ^{\top } \X ( \X ^{\top } \X + \bPsi )^{- 1} \X ^{\top } \y \over 2} + b' \Big] ^{n / 2 + a'} \non \\
&\le {1 \over \{ 1 + u ( \X ^{\top } \X + \bPsi )^{k, k} \} ^{1 / 2}} \le {1 \over \{ 1 + u ( M_1 \I ^{(p)} + \bPsi )^{k, k} \} ^{1 / 2}} = {1 \over \{ 1 + u / ( M_1 + 1 / \lat _{k}^{2} ) \} ^{1 / 2}} \label{lsmallp2} 
\end{align}
for some $M_1 > 0$. 
Meanwhile, by making the change of variables $\mu = 1 / \nu _k$, 
\begin{align}
E \Big[ \exp \Big( - u {\tas \over \nu _k} \Big) \Big| {\blat {}^2} \Big] &= E \Big[ \int_{0}^{\infty } \Big( 1 + {\tas \over \lat _{k}^{2}} \Big) \exp \Big\{ - \mu \Big( 1 + {\tas \over \lat _{k}^{2}} \Big) \Big\} \exp (- \mu u \tas ) d\mu \Big| {\blat {}^2} \Big] \non \\
&= E \Big[ {1 + \tas / \lat _{k}^{2} \over 1 + \tas / \lat _{k}^{2} + u \tas } \Big| {\blat {}^2} \Big] \text{.} \label{lsmallp3} 
\end{align}

By (\ref{lsmallp1}), (\ref{lsmallp2}), and (\ref{lsmallp3}), 
\begin{align}
&E \Big[ E \Big[ {1 \over ( \lat _{k}^{2} )^{\ep }} \Big| \bbe , \sis , \bnu , \tas \Big] \Big| {\blat {}^2} \Big] \non \\
&\le \int_{0}^{\infty } {\Ga (1 + \ep ) \over \Ga ( \ep )} u^{\ep - 1} {1 \over \{ 1 + u / ( M_1 + 1 / \lat _{k}^{2} ) \} ^{1 / 2}} E \Big[ {1 + \tas / \lat _{k}^{2} \over 1 + \tas / \lat _{k}^{2} + u \tas } \Big| {\blat {}^2} \Big] du \non \\
&= \int_{0}^{\infty } {\Ga (1 + \ep ) \over \Ga ( \ep )} w^{\ep - 1} {( M_1 + 1 / \lat _{k}^{2} )^{\ep } \over (1 + w)^{1 / 2}} E \Big[ {1 + \tas / \lat _{k}^{2} \over 1 + \tas / \lat _{k}^{2} + w \tas ( M_1 + 1 / \lat _{k}^{2} )} \Big| {\blat {}^2} \Big] dw \non \\
&\le \int_{0}^{\infty } {\Ga (1 + \ep ) \over \Ga ( \ep )} w^{\ep  - 1} {{M_1}^{\ep } + (1 / \lat _{k}^{2} )^{\ep } \over (1 + w)^{1 / 2}} E \Big[ {1 + \tas / \lat _{k}^{2} \over 1 + \tas / \lat _{k}^{2} + w \tas / \lat _{k}^{2}} \Big| {\blat {}^2} \Big] dw \non \\
&\le {\Ga (1 + \ep ) \over \Ga ( \ep )} B( \ep , 1 / 2 - \ep ) {M_1}^{\ep } + {1 \over ( \lat _{k}^{2} )^{\ep }} {\Ga (1 + \ep ) \over \Ga ( \ep )} \int_{0}^{\infty } {w^{\ep - 1} \over (1 + w)^{1 / 2}} E \Big[ {1 + \tas / \lat _{k}^{2} \over 1 + (1 + w) \tas / \lat _{k}^{2}} \Big| {\blat {}^2} \Big] dw \text{.} \non 
\end{align}
The result follows. 
\hfill$\Box$

\bigskip

\noindent
{\bf Proof of Theorem \ref{thm:lat}.} \ \ By Corollary 9.14 of %
\cite{rc2004}, the result follows if the Markov chain corresponding to ${\blat {}^2}$ is geometrically ergodic. 
By Lemmas \ref{lem:minorization} and \ref{lem:large}, it is sufficient to show that $V_{2, \ep }$ satisfies the drift condition or that there exists $\ep > 0$ such that 
\begin{align}
(P V_{2, \ep } ) ( \v ) \le \De + \de V_{2, \ep } ( \v ) \non 
\end{align}
for all $\v \in (0, \infty )^p$ for some $0 < \de < 1$ and $\De > 0$. 
Fix $0 < \ep < 1 / 2$. 
Then, by Lemma \ref{lem:small}, there exists $M > 0$ such that 
\begin{align}
(P V_{2, \ep } ) ( \v ) &\le M + \sum_{k = 1}^{p} {1 \over {v_k}^{\ep }} {\Ga (1 + \ep ) \over \Ga ( \ep )} \int_{0}^{\infty } {w^{\ep - 1} \over (1 + w)^{1 / 2}} E \Big[ {v_k + \tas \over v_k + (1 + w) \tas } \Big| {\blat {}^2} = \v \Big] dw \non 
\end{align}
for all $\v = ( v_k )_{k = 1}^{p} \in (0, \infty )^p$. 

For part (i), 
\begin{align}
&\sum_{k = 1}^{p} {1 \over ( \lat _{k}^{2} )^{\ep }} {\Ga (1 + \ep ) \over \Ga ( \ep )} \int_{0}^{\infty } {w^{\ep  - 1} \over (1 + w)^{1 / 2}} E \Big[ {\lat _{k}^{2} + \tas \over \lat _{k}^{2} + (1 + w) \tas } \Big| {\blat {}^2} \Big] dw \non \\
&\le \sum_{k = 1}^{p} {1 \over ( \lat _{k}^{2} )^{\ep }} {\Ga (1 + \ep ) \over \Ga ( \ep )} \int_{0}^{\infty } {w^{\ep  - 1} \over (1 + w)^{1 / 2}} \Big[ E \Big[ \Big\{ {\lat _{k}^{2} \over \lat _{k}^{2} + (1 + w) \tas } \Big\} ^{\ep } \Big| {\blat {}^2} \Big] + E \Big[ {1 \over1 + w} \Big| {\blat {}^2} \Big] \Big] dw \non \\
&\le \sum_{k = 1}^{p} {\Ga (1 + \ep ) \over \Ga ( \ep )} \int_{0}^{\infty } {w^{\ep  - 1} \over (1 + w)^{\ep + 1 / 2}} E \Big[{1 \over ( \tas )^{\ep }} \Big| {\blat {}^2}  \Big] dw + \sum_{k = 1}^{p} {1 \over ( \lat _{k}^{2} )^{\ep }} {\Ga (1 + \ep ) \over \Ga ( \ep )} \int_{0}^{\infty } {w^{\ep  - 1} \over (1 + w)^{3 / 2}} dw \text{.} \non 
\end{align}
If $\ep < a - p / 2$, then $E[ \pi _{\ta } ( \tas ) / ( \tas )^{\ep + p / 2} ] < \infty $ and we have 
\begin{align}
E \Big[ {1 \over ( \tas )^{\ep }} \Big| {\blat {}^2} \Big] &= \int_{0}^{\infty } {1 \over ( \tas )^{\ep }} \pi _{\ta } ( \tas ) \Big\{ \prod_{k = 1}^{p} {( \tas )^{1 / 2} \over \lat _{k}^{2} + \tas } \Big\} d\tas / \int_{0}^{\infty } \pi _{\ta } ( \tas ) \Big\{ \prod_{k = 1}^{p} {( \tas )^{1 / 2} \over \lat _{k}^{2} + \tas } \Big\} d\tas \non \\
&= \int_{0}^{\infty } {1 \over ( \tas )^{\ep }} {\pi _{\ta } ( \tas ) \over ( \tas )^{p / 2}} \Big( \prod_{k = 1}^{p} {\tas \over \lat _{k}^{2} + \tas } \Big) d\tas / \int_{0}^{\infty } {\pi _{\ta } ( \tas ) \over ( \tas )^{p / 2}} \Big( \prod_{k = 1}^{p} {\tas \over \lat _{k}^{2} + \tas } \Big) d\tas \non \\
&\le \int_{0}^{\infty } {1 \over ( \tas )^{\ep }} {\pi _{\ta } ( \tas ) \over ( \tas )^{p / 2}} d\tas / \int_{0}^{\infty } {\pi _{\ta } ( \tas ) \over ( \tas )^{p / 2}} d\tas < \infty \non 
\end{align}
by the covariance inequality. 
Meanwhile, since 
\begin{align}
&{\Ga (1 + \ep ) \over \Ga ( \ep )} \int_{0}^{\infty } {w^{\ep  - 1} \over (1 + w)^{3 / 2}} dw = {\Ga (1 + \ep ) \over \Ga ( \ep )} {\Ga ( \ep ) \Ga (3 / 2 -  \ep ) \over \Ga (3 / 2)} = {\Ga (1 + \ep ) \Ga (3 / 2 - \ep ) \over \Ga (3 / 2)} \non 
\end{align}
and since 
\begin{align}
&{\pd \over \pd \tilde{\ep }} \log {\Ga (1 + \tilde{\ep } ) \Ga (3 / 2 - \tilde{\ep } ) \over \Ga (3 / 2)} = \psi (1 + \tilde{\ep } ) - \psi (3 / 2 - \tilde{\ep } ) < 0 \non 
\end{align}
for all $- 1 < \tilde{\ep } < 1 / 4$, it follows that 
\begin{align}
&{\Ga (1 + \ep ) \over \Ga ( \ep )} \int_{0}^{\infty } {w^{\ep  - 1} \over (1 + w)^{3 / 2}} dw < 1 \non 
\end{align}
if $\ep > 0$ is sufficiently small. 
This proves part (i). 

For part (ii), 
\begin{align}
&\sum_{k = 1}^{p} {1 \over ( \lat _{k}^{2} )^{\ep }} {\Ga (1 + \ep ) \over \Ga ( \ep )} \int_{0}^{\infty } {w^{\ep - 1} \over (1 + w)^{1 / 2}} E \Big[ {\lat _{k}^{2} + \tas \over \lat _{k}^{2} + (1 + w) \tas } \Big| {\blat {}^2} \Big] dw \non \\
&\le \Big\{ \sum_{k = 1}^{p} {1 \over ( \lat _{k}^{2} )^{\ep }} \Big\} {1 \over p} \sum_{k = 1}^{p} {\Ga (1 + \ep ) \over \Ga ( \ep )} \int_{0}^{\infty } {w^{\ep - 1} \over (1 + w)^{1 / 2}} E \Big[ {\lat _{k}^{2} + \tas \over \lat _{k}^{2} + (1 + w) \tas } \Big| {\blat {}^2} \Big] dw \non 
\end{align}
by the covariance inequality. 
Note that 
\begin{align}
&E \Big[ {\lat _{k}^{2} + \tas \over \lat _{k}^{2} + (1 + w) \tas } \Big| {\blat {}^2} \Big] \non \\
&= E \Big[ {1 \over 1 + w} + {w \over 1 + w} {\lat _{k}^{2} \over \lat _{k}^{2} + (1 + w) \tas } \Big| {\blat {}^2} \Big] \non \\
&= {1 \over 1 + w} + {w \over 1 + w} \Big[ E \Big[ {\lat _{k}^{2} \over \lat _{k}^{2} + \tas } \Big| {\blat {}^2} \Big] - E \Big[ {w \lat _{k}^{2} \tas \over ( \lat _{k}^{2} + \tas ) \{ \lat _{k}^{2} + (1 + w) \tas \} } \Big| {\blat {}^2} \Big] \Big] \non 
\end{align}
for all $k = 1, \dots , p$. 
Then 
\begin{align}
&\sum_{k = 1}^{p} {1 \over ( \lat _{k}^{2} )^{\ep }} {\Ga (1 + \ep ) \over \Ga ( \ep )} \int_{0}^{\infty } {w^{\ep - 1} \over (1 + w)^{1 / 2}} E \Big[ {\lat _{k}^{2} + \tas \over \lat _{k}^{2} + (1 + w) \tas } \Big| {\blat {}^2} \Big] dw \non \\
&\le \Big\{ \sum_{k = 1}^{p} {1 \over ( \lat _{k}^{2} )^{\ep }} \Big\} {1 \over p} \sum_{k = 1}^{p} {\Ga (1 + \ep ) \over \Ga ( \ep )} \int_{0}^{\infty } {w^{\ep - 1} \over (1 + w)^{1 / 2}} \Big( {1 \over 1 + w} + {w \over 1 + w} E \Big[ {\lat _{k}^{2} \over \lat _{k}^{2} + \tas } \Big| {\blat {}^2} \Big] \Big) dw \non \\
&= \Big\{ \sum_{k = 1}^{p} {1 \over ( \lat _{k}^{2} )^{\ep }} \Big\} {1 \over p} \sum_{k = 1}^{p} {\Ga (1 + \ep ) \over \Ga ( \ep )} \Big\{ B( \ep , 3 / 2 - \ep ) + B( \ep + 1, 1 / 2 - \ep ) E \Big[ {\lat _{k}^{2} \over \lat _{k}^{2} + \tas } \Big| {\blat {}^2} \Big] \Big\} \text{.} \non 
\end{align}
By %
Lemma \ref{lem:K}, 
\begin{align}
&\sum_{k = 1}^{p} {1 \over ( \lat _{k}^{2} )^{\ep }} {\Ga (1 + \ep ) \over \Ga ( \ep )} \int_{0}^{\infty } {w^{\ep - 1} \over (1 + w)^{1 / 2}} E \Big[ {\lat _{k}^{2} + \tas \over \lat _{k}^{2} + (1 + w) \tas } \Big| {\blat {}^2} \Big] dw \non \\
&\le \Big\{ \sum_{k = 1}^{p} {1 \over ( \lat _{k}^{2} )^{\ep }} \Big\} {\Ga (1 + \ep ) \over \Ga ( \ep )} \Big\{ B( \ep , 3 / 2 - \ep ) + B( \ep + 1, 1 / 2 - \ep ) {1 \over p} E \Big[ {p \over 2} + b - (a + b) {c \over c + \tas } \Big| {\blat {}^2} \Big] \Big\} \non \\
&= \Big\{ \sum_{k = 1}^{p} {1 \over ( \lat _{k}^{2} )^{\ep }} \Big\} {\Ga (1 + \ep ) \over \Ga ( \ep )} B( \ep , 3 / 2 - \ep ) \Big\{ 1 + {B( \ep + 1, 1 / 2 - \ep ) \over B( \ep , 3 / 2 - \ep )} {1 \over p} E \Big[ {p \over 2} - a + (a + b) {\tas \over c + \tas } \Big| {\blat {}^2} \Big] \Big\} \text{.} \non 
\end{align}
Now, suppose that $0 < \ep < A_{*}$, where $A_{*} > 0$ is defined as in Lemma \ref{lem:median}, and fix $\eta > 0$. 
Let $\lat _{*}^{2}$ be as in Lemma \ref{lem:median}. 
Then, by part (iii) of Lemma \ref{lem:median}, there exists $\de > 0$ such that for all ${\blat {}^2} \in (0, \infty )^p$ satisfying $\lat _{*}^{2} < \de $, we have 
\begin{align}
&E \Big[ {\tas \over c + \tas } \Big| {\blat {}^2} \Big] < \eta \text{.} \non 
\end{align}
Meanwhile, for all ${\blat {}^2} \in (0, \infty )^p$ satisfying $\lat _{*}^{2} \ge \de $, we have, by part (i) of Lemma \ref{lem:median}, that 
\begin{align}
&E \Big[ {1 \over ( \tas )^{\ep }} \Big| {\blat {}^2} \Big] \non \\
&\le \int_{0}^{\infty } {1 \over ( \tas )^{\ep }} {( \tas )^{A_{*} - 1} \over (c + \tas )^{a + b}} {1 \over ( \lat _{*}^{2} + \tas )^{A_{*} + p / 2 - a}} d\tas / \int_{0}^{\infty } {( \tas )^{A_{*} - 1} \over (c + \tas )^{a + b}} {1 \over ( \lat _{*}^{2} + \tas )^{A_{*} + p / 2 - a}} d\tas \non \\
&\le \int_{0}^{\infty } {1 \over ( \tas )^{\ep }} {( \tas )^{A_{*} - 1} \over (c + \tas )^{a + b}} {1 \over ( \de + \tas )^{A_{*} + p / 2 - a}} d\tas / \int_{0}^{\infty } {( \tas )^{A_{*} - 1} \over (c + \tas )^{a + b}} {1 \over ( \de + \tas )^{A_{*} + p / 2 - a}} d\tas < \infty \text{,} \non 
\end{align}
where the second inequality follows from the covariance inequailty. 
Therefore, if $\lat _{*}^{2} \ge \de $, then as in the proof of part (i), 
\begin{align}
&\sum_{k = 1}^{p} {1 \over ( \lat _{k}^{2} )^{\ep }} {\Ga (1 + \ep ) \over \Ga ( \ep )} \int_{0}^{\infty } {w^{\ep - 1} \over (1 + w)^{1 / 2}} E \Big[ {\lat _{k}^{2} + \tas \over \lat _{k}^{2} + (1 + w) \tas } \Big| {\blat {}^2} \Big] dw \non \\
&\le \sum_{k = 1}^{p} {\Ga (1 + \ep ) \over \Ga ( \ep )} \int_{0}^{\infty } {w^{\ep - 1} \over (1 + w)^{\ep + 1 / 2}} E \Big[ {1 \over ( \tas )^{\ep }} \Big| {\blat {}^2} \Big] dw + \sum_{k = 1}^{p} {1 \over ( \lat _{k}^{2} )^{\ep }} {\Ga (1 + \ep ) \over \Ga ( \ep )} \int_{0}^{\infty } {w^{\ep - 1} \over (1 + w)^{3 / 2}} dw \non \\
&\le M_1 + \Big\{ \sum_{k = 1}^{p} {1 \over ( \lat _{k}^{2} )^{\ep }} \Big\} {\Ga (1 + \ep ) \over \Ga ( \ep )} \int_{0}^{\infty } {w^{\ep - 1} \over (1 + w)^{3 / 2}} dw \non 
\end{align}
for some $M_1 > 0$, where 
\begin{align}
&{\Ga (1 + \ep ) \over \Ga ( \ep )} \int_{0}^{\infty } {w^{\ep - 1} \over (1 + w)^{3 / 2}} dw < 1 \non 
\end{align}
provided $\ep > 0$ is sufficiently small, whereas if $\lat _{*}^{2} < \de $, then 
\begin{align}
&\sum_{k = 1}^{p} {1 \over ( \lat _{k}^{2} )^{\ep }} {\Ga (1 + \ep ) \over \Ga ( \ep )} \int_{0}^{\infty } {w^{\ep - 1} \over (1 + w)^{1 / 2}} E \Big[ {\lat _{k}^{2} + \tas \over \lat _{k}^{2} + (1 + w) \tas } \Big| {\blat {}^2} \Big] dw \non \\
&\le \Big\{ \sum_{k = 1}^{p} {1 \over ( \lat _{k}^{2} )^{\ep }} \Big\} {\Ga (1 + \ep ) \over \Ga ( \ep )} B( \ep , 3 / 2 - \ep ) \Big[ 1 + {B( \ep + 1, 1 / 2 - \ep ) \over B( \ep , 3 / 2 - \ep )} {1 \over p} \Big\{ {p \over 2} - a + (a + b) \eta \Big\} \Big] \non \\
&= \Big\{ \sum_{k = 1}^{p} {1 \over ( \lat _{k}^{2} )^{\ep }} \Big\} {\Ga (1 + \ep ) \Ga (3 / 2 - \ep ) \over \Ga (3 / 2)} \Big( 1 + {\ep \over 1 / 2 - \ep } \rho \Big) \text{,} \non 
\end{align}
where 
\begin{align}
\rho &= {1 \over p} \Big\{ {p \over 2} - a + (a + b) \eta \Big\} \text{.} \non 
\end{align}
We have 
\begin{align}
&\lim_{\tilde{\ep } \to 0} {\pd \over \pd \tilde{\ep }} \log \Big\{ {\Ga (1 + \tilde{\ep } ) \Ga (3 / 2 - \tilde{\ep } ) \over \Ga (3 / 2)} \Big( 1 + {\tilde{\ep } \over 1 / 2 - \tilde{\ep }} \rho \Big) \Big\} \non \\
&= \lim_{\tilde{\ep } \to 0} \Big[ \psi (1 + \tilde{\ep } ) - \psi (3 / 2 - \tilde{\ep } ) + {\{ 1 / (1 / 2 - \tilde{\ep } ) \} \rho + \{ \tilde{\ep } / (1 / 2 - \tilde{\ep } )^2 \} \rho \over 1 + \{ \tilde{\ep } / (1 / 2 - \tilde{\ep } ) \} \rho } \Big] \non \\
&= \psi (1) - \psi (3 / 2) + 2 \rho %
= 2 (- 1 + \log 2 + \rho ) \text{,} \non 
\end{align}
where the last equality follows from properties of the digamma function $\psi $ %
(see \cite{as1972}). 
If $\eta > 0$ is sufficiently small, it follows from the assumption that the right-hand side is negative, or that 
\begin{align}
\{ a - (a + b) \eta \} / p > \log 2 - 1 / 2 \text{,} \non 
\end{align}
and we conclude that 
\begin{align}
&{\Ga (1 + \ep ) \Ga (3 / 2 - \ep ) \over \Ga (3 / 2)} \Big( 1 + {\ep \over 1 / 2 - \ep } \rho \Big) < 1 \non 
\end{align}
provided that $\ep > 0$ is sufficiently small. 
This completes the proof. 
\hfill$\Box$

\begin{lem}
\label{lem:K} 
We have 
\begin{align}
{p \over 2} + a &= E \Big[ (a + b) {\tas \over c + \tas } + \sum_{k = 1}^{p} {\tas \over \lat _{k}^{2} + \tas } \Big| {\blat {}^2} \Big] \non 
\end{align}
and 
\begin{align}
{p \over 2} + b &= E \Big[ (a + b) {c \over c + \tas } + \sum_{k = 1}^{p} {\lat _{k}^{2} \over \lat _{k}^{2} + \tas } \Big| {\blat {}^2} \Big] \text{.} \non 
\end{align}
\end{lem}

\noindent
{\bf Proof%
.} \ \ %
By integration by parts, 
\begin{align}
&\int_{0}^{\infty } {( \tas )^{a - 1} \over (c + \tas )^{a + b}} \Big\{ \prod_{k = 1}^{p} {( \tas )^{1 / 2} \over \lat _{k}^{2} + \tas } \Big\} d\tas \non \\
&= \int_{0}^{\infty } {( \tas )^{p / 2 + a - 1} \over (c + \tas )^{a + b}} \Big( \prod_{k = 1}^{p} {1 \over \lat _{k}^{2} + \tas } \Big) d\tas \non \\
&= \Big[ {1 \over p / 2 + a} {( \tas )^{p / 2 + a} \over (c + \tas )^{a + b}} \prod_{k = 1}^{p} {1 \over \lat _{k}^{2} + \tas } \Big] _{0}^{\infty } \non \\
&\quad - \int_{0}^{\infty } {1 \over p / 2 + a} {( \tas )^{p / 2 + a} \over (c + \tas )^{a + b}} \Big( \prod_{k = 1}^{p} {1 \over \lat _{k}^{2} + \tas } \Big) \Big( - {a + b \over c + \tas } - \sum_{k = 1}^{p} {1 \over \lat _{k}^{2} + \tas } \Big) d\tas \non \\
&= \int_{0}^{\infty } {( \tas )^{a - 1} \over (c + \tas )^{a + b}} \Big\{ \prod_{k = 1}^{p} {( \tas )^{1 / 2} \over \lat _{k}^{2} + \tas } \Big\} {\tas \over p / 2 + a} \Big( {a + b \over c + \tas } + \sum_{k = 1}^{p} {1 \over \lat _{k}^{2} + \tas } \Big) d\tas \non %
\end{align}
as in the proof of Lemma 5 of \cite{k2015}. 
\hfill$\Box$

\bigskip

If $\alt $ is a real number, $[ \alt ]$ denotes the integer part of $\alt $. 

\begin{lem}
\label{lem:median} 
Assume that $a \le p / 2$. 
Then 
\begin{align}
\{ 1, \dots , p \} &\ni \begin{cases} p / 2 - [- a] \text{,} & \text{if $p \in 2 \mathbb{N}$} \text{,} \\ (p + 1) / 2 + [a] \text{,} & \text{if $p \in 2 \mathbb{N} - 1$ and $a - [a] \le 1 / 2$} \text{,} \\ (p + 1) / 2 + [a] + 1 \text{,} & \text{if $p \in 2 \mathbb{N} - 1$ and $a - [a] > 1 / 2$} \text{.} \end{cases} \non 
\end{align}
Let $\lat _{(1)}^{2} \le \dots \le \lat _{(p)}^{2}$ denote the order statistics of $\lat _{1}^{2} , \dots , \lat _{p}^{2}$ and let 
\begin{align}
\lat _{*}^{2} & = \begin{cases} \lat _{(p / 2 - [- a])}^{2} \text{,} & \text{if $p \in 2 \mathbb{N}$} \text{,} \\ \lat _{((p + 1) / 2 + [a])}^{2} \text{,} & \text{if $p \in 2 \mathbb{N} - 1$ and $a - [a] \le 1 / 2$} \text{,} \\ \lat _{((p + 1) / 2 + [a] + 1)}^{2} \text{,} & \text{if $p \in 2 \mathbb{N} - 1$ and $a - [a] > 1 / 2$} \text{.} \end{cases} \non 
\end{align}
\begin{itemize}
\item[{\rm{(i)}}]
For all $\ep > 0$, we have 
\begin{align}
&E \Big[ {1 \over ( \tas )^{\ep }} \Big| {\blat {}^2} \Big] \non \\
&\le \int_{0}^{\infty } {1 \over ( \tas )^{\ep }} {( \tas )^{A_{*} - 1} \over (c + \tas )^{a + b}} {1 \over ( \lat _{*}^{2} + \tas )^{A_{*} + p / 2 - a}} d\tas / \int_{0}^{\infty } {( \tas )^{A_{*} - 1} \over (c + \tas )^{a + b}} {1 \over ( \lat _{*}^{2} + \tas )^{A_{*} + p / 2 - a}} d\tas \text{,} \non 
\end{align}
where 
\begin{align}
A_{*} &= \begin{cases} a + [- a] + 1 \text{,} & \text{if $p \in 2 \mathbb{N}$} \text{,} \\ 1 / 2 + a - [a] \text{,} & \text{if $p \in 2 \mathbb{N} - 1$ and $a - [a] \le 1 / 2$} \text{,} \\ a - [a] - 1 / 2 \text{,} & \text{if $p \in 2 \mathbb{N} - 1$ and $a - [a] > 1 / 2$} \text{.} \end{cases} \non 
\end{align}
\item[{\rm{(ii)}}]
For all $\ga > 0$, we have 
\begin{align}
&E \Big[ {\tas \over \ga + \tas } \Big| {\blat {}^2} \Big] \le \int_{0}^{\infty } {\tas \over \ga + \tas } {( \tas )^{p / 2 + a - 1} \over (c + \tas )^{a + b}} {1 \over ( \lat _{*}^{2} + \tas )^{B_{*}}} d\tas / \int_{0}^{\infty } {( \tas )^{p / 2 + a - 1} \over (c + \tas )^{a + b}} {1 \over ( \lat _{*}^{2} + \tas )^{B_{*}}} d\tas \text{,} \non 
\end{align}
where 
\begin{align}
B_{*} &= \begin{cases} p / 2 - [- a] \text{,} & \text{if $p \in 2 \mathbb{N}$} \text{,} \\ (p + 1) / 2 + [a] \text{,} & \text{if $p \in 2 \mathbb{N} - 1$ and $a - [a] \le 1 / 2$} \text{,} \\ (p + 1) / 2 + [a] + 1 \text{,} & \text{if $p \in 2 \mathbb{N} - 1$ and $a - [a] > 1 / 2$} \text{.} \end{cases} \non 
\end{align}
\item[{\rm{(iii)}}]
Assume that $a \notin \mathbb{N}$ if $p \in 2 \mathbb{N}$. 
Assume that $a \notin \mathbb{N} - 1 / 2$ if $p \in 2 \mathbb{N} - 1$. 
Then %
for all $\ga > 0$, we have 
\begin{align}
&E \Big[ {\tas \over \ga + \tas } \Big| {\blat {}^2} \Big] \to 0 \non 
\end{align}
as $\lat _{*}^{2} \to 0$; %
that is, for all $\ga > 0$ and all $\eta > 0$, there exists $\de > 0$ such that for all ${\blat {}^2} \in (0, \infty )^p$ satisfying $\lat _{*}^{2} < \de $, we have 
\begin{align}
0 &\le E \Big[ {\tas \over \ga + \tas } \Big| {\blat {}^2} \Big] < \eta \text{.} \non 
\end{align}
\end{itemize}
\end{lem}

\noindent
{\bf Proof%
.} \ \ For part (i), fix $\ep > 0$. 
Then 
\begin{align}
&E \Big[ {1 \over ( \tas )^{\ep }} \Big| {\blat {}^2} \Big] = \int_{0}^{\infty } {1 \over ( \tas )^{\ep }} {( \tas )^{p / 2 + a - 1} \over (c + \tas )^{a + b}} \Big( \prod_{k = 1}^{p} {1 \over \lat _{k}^{2} + \tas } \Big) d\tas / \int_{0}^{\infty } {( \tas )^{p / 2 + a - 1} \over (c + \tas )^{a + b}} \Big( \prod_{k = 1}^{p} {1 \over \lat _{k}^{2} + \tas } \Big) d\tas \text{.} \non 
\end{align}
If $p \in 2 \mathbb{N}$, then by the covariance inequality, 
\begin{align}
E \Big[ {1 \over ( \tas )^{\ep }} \Big| {\blat {}^2} \Big] &\le \int_{0}^{\infty } {1 \over ( \tas )^{\ep }} {( \tas )^{a + [- a] + 1 - 1} \over (c + \tas )^{a + b}} \Big( \prod_{k = p / 2 - [- a]}^{p} {1 \over \lat _{(k)}^{2} + \tas } \Big) d\tas \non \\
&\quad / \int_{0}^{\infty } {( \tas )^{a + [- a] + 1 - 1} \over (c + \tas )^{a + b}} \Big( \prod_{k = p / 2 - [- a]}^{p} {1 \over \lat _{(k)}^{2} + \tas } \Big) d\tas \non \\
&\le \int_{0}^{\infty } {1 \over ( \tas )^{\ep }} {( \tas )^{a + [- a] + 1 - 1} \over (c + \tas )^{a + b}} {1 \over ( \lat _{(p / 2 - [- a])}^{2} + \tas )^{p / 2 + [- a] + 1}} d\tas \non \\
&\quad / \int_{0}^{\infty } {( \tas )^{a + [- a] + 1 - 1} \over (c + \tas )^{a + b}} {1 \over ( \lat _{(p / 2 - [- a])}^{2} + \tas )^{p / 2 + [- a] + 1}} d\tas \text{.} \non 
\end{align}
Similarly, if $p \in 2 \mathbb{N} - 1$ and $a - [a] \le 1 / 2$, then 
\begin{align}
E \Big[ {1 \over ( \tas )^{\ep }} \Big| {\blat {}^2} \Big] &\le \int_{0}^{\infty } {1 \over ( \tas )^{\ep }} {( \tas )^{1 / 2 + a - [a] - 1} \over (c + \tas )^{a + b}} \Big\{ \prod_{k = (p + 1) / 2 + [a]}^{p} {1 \over \lat _{(k)}^{2} + \tas } \Big\} d\tas \non \\
&\quad / \int_{0}^{\infty } {( \tas )^{1 / 2 + a - [a] - 1} \over (c + \tas )^{a + b}} \Big\{ \prod_{k = (p + 1) / 2 + [a]}^{p} {1 \over \lat _{(k)}^{2} + \tas } \Big\} d\tas \non \\
&\le \int_{0}^{\infty } {1 \over ( \tas )^{\ep }} {( \tas )^{1 / 2 + a - [a] - 1} \over (c + \tas )^{a + b}} {1 \over \{ \lat _{((p + 1) / 2 + [a])}^{2} + \tas \} ^{(p + 1) / 2 - [a]}} d\tas \non \\
&\quad / \int_{0}^{\infty } {( \tas )^{1 / 2 + a - [a] - 1} \over (c + \tas )^{a + b}} {1 \over \{ \lat _{((p + 1) / 2 + [a])}^{2} + \tas \} ^{(p + 1) / 2 - [a]}} d\tas \text{.} \non 
\end{align}
Furthermore, if $p \in 2 \mathbb{N} - 1$ and if $a - [a] > 1 / 2$, then 
\begin{align}
E \Big[ {1 \over ( \tas )^{\ep }} \Big| {\blat {}^2} \Big] &\le \int_{0}^{\infty } {1 \over ( \tas )^{\ep }} {( \tas )^{a - [a] - 1 / 2 - 1} \over (c + \tas )^{a + b}} \Big\{ \prod_{k = (p + 1) / 2 + [a] + 1}^{p} {1 \over \lat _{(k)}^{2} + \tas } \Big\} d\tas \non \\
&\quad / \int_{0}^{\infty } {( \tas )^{a - [a] - 1 / 2 - 1} \over (c + \tas )^{a + b}} \Big\{ \prod_{k = (p + 1) / 2 + [a] + 1}^{p} {1 \over \lat _{(k)}^{2} + \tas } \Big\} d\tas \non \\
&\le \int_{0}^{\infty } {1 \over ( \tas )^{\ep }} {( \tas )^{a - [a] - 1 / 2 - 1} \over (c + \tas )^{a + b}} {1 \over \{ \lat _{((p + 1) / 2 + [a] + 1)}^{2} + \tas \} ^{(p - 1) / 2 - [a]}} d\tas \non \\
&\quad / \int_{0}^{\infty } {( \tas )^{a - [a] - 1 / 2 - 1} \over (c + \tas )^{a + b}} {1 \over \{ \lat _{((p + 1) / 2 + [a] + 1)}^{2} + \tas \} ^{(p - 1) / 2 - [a]}} d\tas \text{.} \non 
\end{align}
Thus, part (i) follows. 

For part (ii), fix $\ga > 0$. 
Then 
\begin{align}
&E \Big[ {\tas \over \ga + \tas } \Big| {\blat {}^2} \Big] = \int_{0}^{\infty } {\tas \over \ga + \tas } {( \tas )^{p / 2 + a - 1} \over (c + \tas )^{a + b}} \Big( \prod_{k = 1}^{p} {1 \over \lat _{k}^{2} + \tas } \Big) d\tas / \int_{0}^{\infty } {( \tas )^{p / 2 + a - 1} \over (c + \tas )^{a + b}} \Big( \prod_{k = 1}^{p} {1 \over \lat _{k}^{2} + \tas } \Big) d\tas \text{.} \non 
\end{align}
As in the proof of part (i), it follows from the covariance inequality that 
\begin{align}
&E \Big[ {\tas \over \ga + \tas } \Big| {\blat {}^2} \Big] %
\le \int_{0}^{\infty } {\tas \over \ga + \tas } {( \tas )^{p / 2 + a - 1} \over (c + \tas )^{a + b}} \Big( \prod_{k = 1}^{B_{*}} {1 \over \lat _{*}^{2} + \tas } \Big) d\tas / \int_{0}^{\infty } {( \tas )^{p / 2 + a - 1} \over (c + \tas )^{a + b}} \Big( \prod_{k = 1}^{B_{*}} {1 \over \lat _{*}^{2} + \tas } \Big) d\tas \text{,} \non 
\end{align}
which proves part (ii). 

For part (iii), assume that 
\begin{align}
&\begin{cases} a \notin \mathbb{N} \text{,} & \text{if $p \in 2 \mathbb{N}$} \text{,} \\ a \notin \mathbb{N} - 1 / 2 \text{,} & \text{if $p \in 2 \mathbb{N} - 1$} \text{.} \end{cases} \non 
\end{align}
Let $A_{*}$ and $B_{*}$ be as in parts (ii) and (iii). 
Then $A_{*} > 0$ and $B_{*} > p / 2 + a$. 
Let $C_{*} = \min \{ 1, B_{*} - (p / 2 + a) \} / 2$. 
Then $C_{*} \in (0, 1)$. 
Fix $\ga > 0$. 
Then by part (ii) and the monotone convergence theorem, 
\begin{align}
&E \Big[ {\tas \over \ga + \tas } \Big| {\blat {}^2} \Big] \non \\
&\le \int_{0}^{\infty } \Big( {\tas \over \ga + \tas } \Big) ^{C_{*}} {( \tas )^{p / 2 + a - 1} \over (c + \tas )^{a + b}} {1 \over ( \lat _{*}^{2} + \tas )^{B_{*}}} d\tas / \int_{0}^{\infty } {( \tas )^{p / 2 + a - 1} \over (c + \tas )^{a + b}} {1 \over ( \lat _{*}^{2} + \tas )^{B_{*}}} d\tas \non \\
&= ( \lat _{*}^{2} )^{C_{*}} \int_{0}^{\infty } \Big( {\tilde{\ta } ^2 \over \ga + \lat _{*}^{2} \tilde{\ta } ^2} \Big) ^{C_{*}} {( \tilde{\ta } ^2 )^{p / 2 + a - 1} \over (c + \lat _{*}^{2} \tilde{\ta } ^2 )^{a + b}} {1 \over (1 + \tilde{\ta } ^2 )^{B_{*}}} d{\tilde{\ta } ^2} / \int_{0}^{\infty } {( \tilde{\ta } ^2 )^{p / 2 + a - 1} \over (c + \lat _{*}^{2} \tilde{\ta } ^2 )^{a + b}} {1 \over (1 + \tilde{\ta } ^2 )^{B_{*}}} d{\tilde{\ta } ^2} \non \\
&\sim {( \lat _{*}^{2} )^{C_{*}} \over \ga ^{C_{*}}} \int_{0}^{\infty } {( \tilde{\ta } ^2 )^{C_{*} + p / 2 + a - 1} \over (1 + \tilde{\ta } ^2 )^{B_{*}}} d{\tilde{\ta } ^2} / \int_{0}^{\infty } {( \tilde{\ta } ^2 )^{p / 2 + a - 1} \over (1 + \tilde{\ta } ^2 )^{B_{*}}} d{\tilde{\ta } ^2} \non \\
&= {( \lat _{*}^{2} )^{C_{*}} \over \ga ^{C_{*}}} {B( C_{*} + p / 2 + a, B_{*} - ( C_{*} + p / 2 + a)) \over B(p / 2 + a, B_{*} - (p / 2 + a))} \text{,} \non 
\end{align}
where 
\begin{align}
&B_{*} - ( C_{*} + p / 2 + a) \ge B_{*} - [ \{ B_{*} - (p / 2 + a) \} / 2 + p / 2 + a] = \{ B_{*} - (p / 2 + a) \} / 2 > 0 \text{.} \non 
\end{align}
Thus, $E[ \tas / ( \ga + \tas ) | {\blat {}^2} ] \to 0$ as $\lat _{*}^{2} \to 0$. 
This completes the proof of part (iii). 
\hfill$\Box$

\begin{remark}
\label{rem:small_a} 
As shown in the next section, the drift condition does not hold for the energy function $V_{2, \ep }$ for all $\ep > 0$ near zero when $a$ is smaller. 
We have several possibilities for proving geometric ergodicity. 
We may consider trying 
\begin{itemize}
\item
$V_{2, \ep }$ for some $\ep > 0$ not very close to zero, 
\item
a different form of energy function, or 
\item
a different algorithm or parameterization. 
\end{itemize}
However, it is not standard practice to try the first or second possibilities, and such approaches may not be fruitful. 
In Section \ref{sec:general}, we saw that we can dispense with the assumption (\ref{eq:condition_ii}) if we consider different algorithms based on a reparameterized model. 
\end{remark}

\section{Failure of a Usual Energy Function When $a$ is Not Large}
\label{sec:failure} 
The following proposition shows that when $a$ is smaller, the drift condition does not hold for the energy function $V_{2, \ep }$ for all $\ep > 0$ near zero. 

\begin{prp}
\label{prp:bad_drift_function} 
Assume that $a \notin \mathbb{N}$ if $p \in 2 \mathbb{N}$. 
Assume that $a \notin \mathbb{N} - 1 / 2$ if $p \in 2 \mathbb{N} - 1$. 
Suppose that $a + 1 < \{ (2 \log 2 - 1) / 2 \} p \le p / 2$. 
Then there exists $0 < \ep _0 < 1 / 2$ such that there is no $0 < \ep _1 < \ep _0$ satisfying the condition that for some $0 < \de < 1$ and $\De > 0$, we have for all $\v \in (0, \infty )^p$ that 
\begin{align}
&(P V_{2, \ep _1} ) ( \v ) \le \De + \de V_{2, \ep _1} ( \v ) \text{.} \label{eq:condition_bad} 
\end{align}
\end{prp}

\noindent
{\bf Proof%
.} \ \ Fix $0 < \ep < 1 / 2$. 
Fix $k = 1, \dots , p$ and $u \in (0, \infty )$. 
By (\ref{lsmallp2}), 
\begin{align}
&E \Big[ \exp \Big( - u {{\be _k}^2 \over 2 \sis } \Big) \Big| {\blat {}^2} \Big] %
\ge {1 \over (1 + u \bPsi ^{k, k} )^{1 / 2}} ( b' )^{n / 2 + a'} / \Big[ {u \{ ( \e _{k}^{(p)} )^{\top } ( \X ^{\top } \X + \bPsi )^{- 1} \X ^{\top } \y \} ^2 \over 2 \{ 1 + u ( \X ^{\top } \X + \bPsi )^{k, k} \} } + b' \Big] ^{n / 2 + a'} \text{.} \non 
\end{align}
Note that 
\begin{align}
\{ ( \e _{k}^{(p)} )^{\top } ( \X ^{\top } \X + \bPsi )^{- 1} \X ^{\top } \y \} ^2 &\le ( \e _{k}^{(p)} )^{\top } ( \X ^{\top } \X + \bPsi )^{- 1} ( M_1 \X ^{\top } \X ) ( \X ^{\top } \X + \bPsi )^{- 1} \e _{k}^{(p)} \non \\
&\le M_1 ( \e _{k}^{(p)} )^{\top } \{ \bPsi ( \X ^{\top } \X )^{- 1} \bPsi \} ^{- 1} \e _{k}^{(p)} \non \\
&\le M_2 ( \e _{k}^{(p)} )^{\top } \bPsi ^{- 1} \I ^{(p)} \bPsi ^{- 1} \e _{k}^{(p)} = M_2 ( \lat _{k}^{2} )^2 \non 
\end{align}
for some $M_1 , M_2 > 0$. 
Then 
\begin{align}
E \Big[ \exp \Big( - u {{\be _k}^2 \over 2 \sis } \Big) \Big| {\blat {}^2} \Big] &\ge {1 \over (1 + u \bPsi ^{k, k} )^{1 / 2}} ( b' )^{n / 2 + a'} / \Big[ {u M_2 ( \lat _{k}^{2} )^2 \over 2 \{ 1 + u ( M_3 \I ^{(p)} + \bPsi )^{k, k} \} } + b' \Big] ^{n / 2 + a'} \non \\
&= {1 \over (1 + u \lat _{k}^{2} )^{1 / 2}} / \Big[ {u M_2 ( \lat _{k}^{2} )^2 / b' \over 2 \{ 1 + u / ( M_3 + 1 / \lat _{k}^{2} ) \} } + 1 \Big] ^{n / 2 + a'} \label{eq:lower_1} 
\end{align}
for some $M_3 > 0$. 

By (\ref{lsmallp1}), (\ref{lsmallp3}), and (\ref{eq:lower_1}), 
\begin{align}
(P V_{2, \ep } ) ( {\blat {}^2} ) &\ge \sum_{k = 1}^{p} \int_{0}^{\infty } \Big\{ {\Ga (1 + \ep ) \over \Ga ( \ep )} u^{\ep - 1} \non \\
&\quad \times \Big( {1 \over (1 + u \lat _{k}^{2} )^{1 / 2}} / \Big[ {u M_2 ( \lat _{k}^{2} )^2 / b' \over 2 \{ 1 + u / ( M_3 + 1 / \lat _{k}^{2} ) \} } + 1 \Big] ^{n / 2 + a'} \Big) E \Big[ {1 + \tas / \lat _{k}^{2} \over 1 + \tas / \lat _{k}^{2} + u \tas } \Big| {\blat {}^2} \Big] \Big\} du \text{.} \non 
\end{align}
Therefore, 
\begin{align}
(P V_{2, \ep } ) ( {\blat {}^2} ) &\ge \sum_{k = 1}^{p} \Big[ {1 \over ( \lat _{k}^{2} )^{\ep }} \int_{0}^{\infty } \Big\{ {\Ga (1 + \ep ) \over \Ga ( \ep )} w^{\ep - 1} \non \\
&\quad \times \Big( {1 \over (1 + w)^{1 / 2}} / \Big[ {w M_2 \lat _{k}^{2} / b' \over 2 \{ 1 + w / ( M_3 \lat _{k}^{2} + 1) \} } + 1 \Big] ^{n / 2 + a'} \Big) E \Big[ {\lat _{k}^{2} + \tas \over \lat _{k}^{2} + (1 + w) \tas } \Big| {\blat {}^2} \Big] \Big\} dw \Big] \text{.} \label{eq:v3_0} 
\end{align}

Fix $0 < K < 1$. 
Then by (\ref{eq:v3_0}), 
\begin{align}
&(P V_{2, \ep } ) (K {\blat {}^2} ) \non \\
&\ge {1 \over K^{\ep }} \sum_{k = 1}^{p} \Big[ {1 \over ( \lat _{k}^{2} )^{\ep }} \int_{0}^{\infty } \Big\{ {\Ga (1 + \ep ) \over \Ga ( \ep )} w^{\ep - 1} \non \\
&\quad \times \Big( {1 \over (1 + w)^{1 / 2}} / \Big[ {w M_2 K \lat _{k}^{2} / b' \over 2 \{ 1 + w / ( M_3 K \lat _{k}^{2} + 1) \} } + 1 \Big] ^{n / 2 + a'} \Big) E \Big[ {K \lat _{k}^{2} + \tas \over K \lat _{k}^{2} + (1 + w) \tas } \Big| K {\blat {}^2} \Big] \Big\} dw \Big] \non \\
&\ge {1 \over K^{\ep }} \sum_{k = 1}^{p} \Big[ {1 \over ( \lat _{k}^{2} )^{\ep }} \int_{0}^{\infty } \Big\{ {\Ga (1 + \ep ) \over \Ga ( \ep )} w^{\ep - 1} \non \\
&\quad \times \Big( {1 \over (1 + w)^{1 / 2}} / \Big[ {w M_2 \lat _{k}^{2} / b' \over 2 \{ 1 + w / ( M_3 \lat _{k}^{2} + 1) \} } + 1 \Big] ^{n / 2 + a'} \Big) E \Big[ {K \lat _{k}^{2} + \tas \over K \lat _{k}^{2} + (1 + w) \tas } \Big| K {\blat {}^2} \Big] \Big\} dw \Big] \text{.} \non 
\end{align}
If $g \colon (0, \infty ) \to [0, \infty )$, we write 
\begin{align}
\widehat{E} [ g( \tas ) | \blat ] &= \int_{0}^{\infty } g( \hat{\ta } ^2 ) ( \hat{\ta } ^2 )^{p / 2 + a - 1} \Big( \prod_{k = 1}^{p} {1 \over \lat _{k}^{2} + \hat{\ta } ^2} \Big) d{\hat{\ta } ^2} / \int_{0}^{\infty } ( \hat{\ta } ^2 )^{p / 2 + a - 1} \Big( \prod_{k = 1}^{p} {1 \over \lat _{k}^{2} + \hat{\ta } ^2} \Big) d{\hat{\ta } ^2} \text{.} \non 
\end{align}
For all $k = 1, \dots , p$, we have 
\begin{align}
&E \Big[ {K \lat _{k}^{2} + \tas \over K \lat _{k}^{2} + (1 + w) \tas } \Big| K {\blat {}^2} \Big] \non \\
&= \int_{0}^{\infty } {K \lat _{k}^{2} + \tas \over K \lat _{k}^{2} + (1 + w) \tas } {( \tas )^{p / 2 + a - 1} \over (c + \tas )^{a + b}} \Big( \prod_{k = 1}^{p} {1 \over K \lat _{k}^{2} + \tas } \Big) d\tas / \int_{0}^{\infty } {( \tas )^{p / 2 + a - 1} \over (c + \tas )^{a + b}} \Big( \prod_{k = 1}^{p} {1 \over K \lat _{k}^{2} + \tas } \Big) d\tas \non \\
&= \int_{0}^{\infty } {\lat _{k}^{2} + \hat{\ta } ^2 \over \lat _{k}^{2} + (1 + w) \hat{\ta } ^2} {( \hat{\ta } ^2 )^{p / 2 + a - 1} \over (c / K + \hat{\ta } ^2 )^{a + b}} \Big( \prod_{k = 1}^{p} {1 \over \lat _{k}^{2} + \hat{\ta } ^2} \Big) d{\hat{\ta } ^2} / \int_{0}^{\infty } {( \hat{\ta } ^2 )^{p / 2 + a - 1} \over (c / K + \hat{\ta } ^2 )^{a + b}} \Big( \prod_{k = 1}^{p} {1 \over \lat _{k}^{2} + \hat{\ta } ^2} \Big) d{\hat{\ta } ^2} \non \\
&\ge \int_{0}^{\infty } {\lat _{k}^{2} + \hat{\ta } ^2 \over \lat _{k}^{2} + (1 + w) \hat{\ta } ^2} ( \hat{\ta } ^2 )^{p / 2 + a - 1} \Big( \prod_{k = 1}^{p} {1 \over \lat _{k}^{2} + \hat{\ta } ^2} \Big) d{\hat{\ta } ^2} / \int_{0}^{\infty } ( \hat{\ta } ^2 )^{p / 2 + a - 1} \Big( \prod_{k = 1}^{p} {1 \over \lat _{k}^{2} + \hat{\ta } ^2} \Big) d{\hat{\ta } ^2} \non \\
&= \widehat{E} \Big[ {\lat _{k}^{2} + \tas \over \lat _{k}^{2} + (1 + w) \tas } \Big| {\blat {}^2} \Big] \text{,} \non 
\end{align}
where the inequality follows from the covariance inequality. 
Therefore, 
\begin{align}
&(P V_{2, \ep } ) (K {\blat {}^2} ) \ge {1 \over K^{\ep }} \sum_{k = 1}^{p} {W_k ( {\blat {}^2} , \ep ) \over ( \lat _{k}^{2} )^{\ep }} \text{,} \non 
\end{align}
where 
\begin{align}
&W_k ( {\blat {}^2} , \ep ) \non \\
&= \int_{0}^{\infty } \Big\{ {\Ga (1 + \ep ) \over \Ga ( \ep )} w^{\ep - 1} \Big( {1 \over (1 + w)^{1 / 2}} / \Big[ {w M_2 \lat _{k}^{2} / b' \over 2 \{ 1 + w / ( M_3 \lat _{k}^{2} + 1) \} } + 1 \Big] ^{n / 2 + a'} \Big) \widehat{E} \Big[ {\lat _{k}^{2} + \tas \over \lat _{k}^{2} + (1 + w) \tas } \Big| {\blat {}^2} \Big] \Big\} dw \non 
\end{align}
for $k = 1, \dots , p$. 
Thus, 
\begin{align}
{(P V_{2, \ep } ) (K {\blat {}^2} ) \over V_{2, \ep } (K {\blat {}^2} )} &\ge {1 \over K^{\ep }} \Big\{ \sum_{k = 1}^{p} {W_k ( {\blat {}^2} , \ep ) \over ( \lat _{k}^{2} )^{\ep }} \Big\} / \Big\{ {1 \over K^{\ep }} \sum_{k = 1}^{p} {1 \over ( \lat _{k}^{2} )^{\ep }} \Big\} = \Big\{ \sum_{k = 1}^{p} {W_k ( {\blat {}^2} , \ep ) \over ( \lat _{k}^{2} )^{\ep }} \Big\} / \sum_{k = 1}^{p} {1 \over ( \lat _{k}^{2} )^{\ep }} \text{.} \label{eq:v3_0.5} 
\end{align}

For all $k = 1, \dots , p$, by integration by parts, 
\begin{align}
&W_k ( {\blat {}^2} , \ep ) \non \\
&= \widehat{E} \Big[ \int_{0}^{\infty } \Big\{ w^{\ep } \Big( {1 \over (1 + w)^{1 / 2}} / \Big[ {w M_2 \lat _{k}^{2} / b' \over 2 \{ 1 + w / ( M_3 \lat _{k}^{2} + 1) \} } + 1 \Big] ^{n / 2 + a'} \Big) {\lat _{k}^{2} + \tas \over \lat _{k}^{2} + (1 + w) \tas } \non \\
&\quad \times \Big[ {1 / 2 \over 1 + w} + {\tas \over \lat _{k}^{2} + (1 + w) \tas } + {M_2 \lat _{k}^{2} / b' \over 2 + w \{ M_2 \lat _{k}^{2} / b' + 2 / ( M_3 \lat _{k}^{2} + 1) \} } {n / 2 + a' \over 1 + w / ( M_3 \lat _{k}^{2} + 1)} \Big] \Big\} dw \Big| {\blat {}^2} \Big] \non 
\end{align}
since 
\begin{align}
&{\pd \over \pd w} \log \Big[ {w M_2 \lat _{k}^{2} / b' \over 2 \{ 1 + w / ( M_3 \lat _{k}^{2} + 1) \} } + 1 \Big] = {M_2 \lat _{k}^{2} / b' \over 2 + w \{ M_2 \lat _{k}^{2} / b' + 2 / ( M_3 \lat _{k}^{2} + 1) \} } {1 \over 1 + w / ( M_3 \lat _{k}^{2} + 1)} \non 
\end{align}
for all $w \in (0, \infty )$. 
Therefore, for all $k = 1, \dots , p$, 
\begin{align}
&\lim_{\ep \to 0} W_k ( {\blat {}^2} , \ep ) \non \\
&= \widehat{E} \Big[ \int_{0}^{\infty } \Big\{ \Big( {1 \over (1 + w)^{1 / 2}} / \Big[ {w M_2 \lat _{k}^{2} / b' \over 2 \{ 1 + w / ( M_3 \lat _{k}^{2} + 1) \} } + 1 \Big] ^{n / 2 + a'} \Big) {\lat _{k}^{2} + \tas \over \lat _{k}^{2} + (1 + w) \tas } \non \\
&\quad \times \Big[ {1 / 2 \over 1 + w} + {\tas \over \lat _{k}^{2} + (1 + w) \tas } + {M_2 \lat _{k}^{2} / b' \over 2 + w \{ M_2 \lat _{k}^{2} / b' + 2 / ( M_3 \lat _{k}^{2} + 1) \} } {n / 2 + a' \over 1 + w / ( M_3 \lat _{k}^{2} + 1)} \Big] \Big\} dw \Big| {\blat {}^2} \Big] = 1 \non 
\end{align}
and 
\begin{align}
&\lim_{\ep \to 0} {\pd \over \pd \ep } W_k ( {\blat {}^2} , \ep ) = W_{k}^{\dag } ( {\blat {}^2} ) \text{,} \non 
\end{align}
where 
\begin{align}
W_{k}^{\dag } ( {\blat {}^2} ) &= \int_{0}^{\infty } \Big\{ ( \log w) \Big( {1 \over (1 + w)^{1 / 2}} / \Big[ {w M_2 \lat _{k}^{2} / b' \over 2 \{ 1 + w / ( M_3 \lat _{k}^{2} + 1) \} } + 1 \Big] ^{n / 2 + a'} \Big) \widehat{E} \Big[ {\lat _{k}^{2} + \tas \over \lat _{k}^{2} + (1 + w) \tas } \times \Big[ \non \\
&\quad {1 / 2 \over 1 + w} + {\tas \over \lat _{k}^{2} + (1 + w) \tas } + {M_2 \lat _{k}^{2} / b' \over 2 + w \{ M_2 \lat _{k}^{2} / b' + 2 / ( M_3 \lat _{k}^{2} + 1) \} } {n / 2 + a' \over 1 + w / ( M_3 \lat _{k}^{2} + 1)} \Big] \Big| {\blat {}^2} \Big] \Big\} dw \non 
\end{align}
for $k = 1, \dots , p$. 
Thus, 
\begin{align}
&\lim_{\ep \to 0} \Big\{ \sum_{k = 1}^{p} {W_k ( {\blat {}^2} , \ep ) \over ( \lat _{k}^{2} )^{\ep }} \Big\} / \sum_{k = 1}^{p} {1 \over ( \lat _{k}^{2} )^{\ep }} = 1 \label{eq:v3_0.75} 
\end{align}
and 
\begin{align}
&\lim_{\ep \to 0} {\pd \over \pd \ep } \Big[ \Big\{ \sum_{k = 1}^{p} {W_k ( {\blat {}^2}, \ep ) \over ( \lat _{k}^{2} )^{\ep }} \Big\} / \sum_{k = 1}^{p} {1 \over ( \lat _{k}^{2} )^{\ep }} \Big] \non \\
&= \lim_{\ep \to 0} \Big[ \Big\{ \sum_{k = 1}^{p} {1 \over ( \lat _{k}^{2} )^{\ep }} {\pd \over \pd \ep } W_k ( {\blat {}^2}, \ep ) + \sum_{k = 1}^{p} {W_k ( {\blat {}^2}, \ep ) \over ( \lat _{k}^{2} )^{\ep }} \log {1 \over \lat _{k}^{2}} \Big\} / \sum_{k = 1}^{p} {1 \over ( \lat _{k}^{2} )^{\ep }} \non \\
&\quad - \Big\{ \sum_{k = 1}^{p} {W_k ( {\blat {}^2}, \ep ) \over ( \lat _{k}^{2} )^{\ep }} \Big\} \Big\{ \sum_{k = 1}^{p} {1 \over ( \lat _{k}^{2} )^{\ep }} \log {1 \over \lat _{k}^{2}} \Big\} / \Big\{ \sum_{k = 1}^{p} {1 \over ( \lat _{k}^{2} )^{\ep }} \Big\} ^2 \Big] \non \\
&= {1 \over p} \sum_{k = 1}^{p} \lim_{\ep \to 0} {\pd \over \pd \ep } W_k ( {\blat {}^2}, \ep ) = {1 \over p} \sum_{k = 1}^{p} W_{k}^{\dag } ( {\blat {}^2} ) \text{.} \label{eq:v3_1} 
\end{align}

Now, assume $a \le p / 2$ and let 
\begin{align}
k_{*} &= \begin{cases} p / 2 - [- a] \text{,} & \text{if $p \in 2 \mathbb{N}$} \text{,} \\ (p + 1) / 2 + [a] \text{,} & \text{if $p \in 2 \mathbb{N} - 1$ and $a - [a] \le 1 / 2$} \text{,} \\ (p + 1) / 2 + [a] + 1 \text{,} & \text{if $p \in 2 \mathbb{N} - 1$ and $a - [a] > 1 / 2$} \text{.} \end{cases} \non 
\end{align}
Then $k_{*} \in \{ 1, \dots , p \} $ by Lemma \ref{lem:median}. 
Suppose that $\lat _{1}^{2} = \dots = \lat _{k_{*} - 1}^{2} = \latsu $ and that $\lat _{k_{*} + 1}^{2} = \dots = \lat _{p}^{2} = \latso $ for $0 < \latsu \ll \latsm \ll \latso \ll 1$ with $\latso + \latsm / \latso + \latsu / \latsm \to 0$. 
Then for all $g \colon (0, \infty ) \to [0, \infty )$, we have 
\begin{align}
\widehat{E} [ g( \tas ) | \blat ] &= \int_{0}^{\infty } g( \hat{\ta } ^2 ) ( \hat{\ta } ^2 )^{p / 2 + a - 1} {1 \over ( \latsu + \hat{\ta } ^2 )^{k_{*} - 1}} {1 \over \latsm + \hat{\ta } ^2} {1 \over ( \latso + \hat{\ta } ^2 )^{p - k_{*}}} d{\hat{\ta } ^2} \non \\
&\quad / \int_{0}^{\infty } ( \hat{\ta } ^2 )^{p / 2 + a - 1} {1 \over ( \latsu + \hat{\ta } ^2 )^{k_{*} - 1}} {1 \over \latsm + \hat{\ta } ^2} {1 \over ( \latso + \hat{\ta } ^2 )^{p - k_{*}}} d{\hat{\ta } ^2} \non \\
&= \int_{0}^{\infty } g( \latsm \tats ) ( \tats )^{p / 2 + a - 1} {1 \over ( \latsu / \latsm + \tats )^{k_{*} - 1}} {1 \over 1 + \tats } {1 \over ( \latso / \latsm + \tats )^{p - k_{*}}} d\tats \non \\
&\quad / \int_{0}^{\infty } ( \tats )^{p / 2 + a - 1} {1 \over ( \latsu / \latsm + \tats )^{k_{*} - 1}} {1 \over 1 + \tats } {1 \over ( \latso / \latsm + \tats )^{p - k_{*}}} d\tats \text{.} \non 
\end{align}
Fix $k = 1, \dots , p$. 
Then as $\latso + \latsm / \latso + \latsu / \latsm \to 0$, 
\begin{align}
&\widehat{E} \Big[ {\lat _{k}^{2} + \tas \over \lat _{k}^{2} + (1 + w) \tas } \Big| \blat \Big] \non \\
&\to \begin{cases} \displaystyle \int_{0}^{\infty } {1 \over 1 + w} {( \tats )^{p / 2 + a - k_{*} + 1 - 1} \over 1 + \tats } d\tas / \int_{0}^{\infty } {( \tats )^{p / 2 + a - k_{*} + 1 - 1} \over 1 + \tats } d\tats \text{,} & \text{if $k < k_{*}$} \text{,} \\ \displaystyle \int_{0}^{\infty } {1 + \tats \over 1 + (1 + w) \tats } {( \tats )^{p / 2 + a - k_{*} + 1 - 1} \over 1 + \tats } d\tas / \int_{0}^{\infty } {( \tats )^{p / 2 + a - k_{*} + 1 - 1} \over 1 + \tats } d\tats \text{,} & \text{if $k = k_{*}$} \text{,} \\ \displaystyle \int_{0}^{\infty } {( \tats )^{p / 2 + a - k_{*} + 1 - 1} \over 1 + \tats } d\tas / \int_{0}^{\infty } {( \tats )^{p / 2 + a - k_{*} + 1 - 1} \over 1 + \tats } d\tats \text{,} & \text{if $k > k_{*}$} \text{,} \end{cases} \non \\
&= \begin{cases} \displaystyle {1 \over 1 + w} \text{,} & \text{if $k < k_{*}$} \text{,} \\ \displaystyle {1 \over (1 + w)^{p / 2 + a - k_{*} + 1}} \text{,} & \text{if $k = k_{*}$} \text{,} \\ \displaystyle 1 \text{,} & \text{if $k > k_{*}$} \text{,} \end{cases} \non 
\end{align}
where 
\begin{align}
&\int_{0}^{\infty } {( \tats )^{p / 2 + a - k_{*} + 1 - 1} \over 1 + \tats } d\tats < \infty \non 
\end{align}
under the assumption that 
\begin{align}
&\begin{cases} a \notin \mathbb{N} \text{,} & \text{if $p \in 2 \mathbb{N}$} \text{,} \\ a \notin \mathbb{N} - 1 / 2 \text{,} & \text{if $p \in 2 \mathbb{N} - 1$} \text{.} \end{cases} \non 
\end{align}
Also, 
\begin{align}
&\widehat{E} \Big[ {\lat _{k}^{2} + \tas \over \lat _{k}^{2} + (1 + w) \tas } {\tas \over \lat _{k}^{2} + (1 + w) \tas } \Big| \blat \Big] \non \\
&\to \begin{cases} \displaystyle \int_{0}^{\infty } {1 \over (1 + w)^2} {( \tats )^{p / 2 + a - k_{*} + 1 - 1} \over 1 + \tats } d\tas / \int_{0}^{\infty } {( \tats )^{p / 2 + a - k_{*} + 1 - 1} \over 1 + \tats } d\tats \text{,} & \text{if $k < k_{*}$} \text{,} \\ \displaystyle \int_{0}^{\infty } {(1 + \tats ) \tats \over \{ 1 + (1 + w) \tats \} ^2} {( \tats )^{p / 2 + a - k_{*} + 1 - 1} \over 1 + \tats } d\tas / \int_{0}^{\infty } {( \tats )^{p / 2 + a - k_{*} + 1 - 1} \over 1 + \tats } d\tats \text{,} & \text{if $k = k_{*}$} \text{,} \\ \displaystyle \int_{0}^{\infty } 0 {( \tats )^{p / 2 + a - k_{*} + 1 - 1} \over 1 + \tats } d\tas / \int_{0}^{\infty } {( \tats )^{p / 2 + a - k_{*} + 1 - 1} \over 1 + \tats } d\tats \text{,} & \text{if $k > k_{*}$} \text{,} \end{cases} \non \\
&= \begin{cases} \displaystyle {1 \over (1 + w)^2} \text{,} & \text{if $k < k_{*}$} \text{,} \\ \displaystyle {p / 2 + a - k_{*} + 1 \over (1 + w)^{p / 2 + a - k_{*} + 2}} \text{,} & \text{if $k = k_{*}$} \text{,} \\ \displaystyle 0 \text{,} & \text{if $k > k_{*}$} \text{,} \end{cases} \non 
\end{align}
as $\latso + \latsm / \latso + \latsu / \latsm \to 0$. 
Therefore, 
\begin{align}
&\lim_{\latso + \latsm / \latso + \latsu / \latsm \to 0} W_{k}^{\dag } ( {\blat {}^2} ) \non \\
&= \begin{cases} \displaystyle \int_{0}^{\infty } {\log w \over (1 + w)^{1 / 2}} \Big\{ {1 / 2 \over (1 + w)^2} + {1 \over (1 + w)^2} \Big\} dw \text{,} & \text{if $k < k_{*}$} \text{,} \\ \displaystyle \int_{0}^{\infty } {\log w \over (1 + w)^{1 / 2}} \Big\{ {1 / 2 \over (1 + w)^{p / 2 + a - k_{*} + 2}} + {p / 2 + a - k_{*} + 1 \over (1 + w)^{p / 2 + a - k_{*} + 2}} \Big\} dw \text{,} & \text{if $k = k_{*}$} \text{,} \\ \displaystyle \int_{0}^{\infty } {\log w \over (1 + w)^{1 / 2}} \Big\{ {1 / 2 \over 1 + w} + 0 \Big\} dw \text{,} & \text{if $k > k_{*}$} \text{,} \end{cases} \non \\
&= \begin{cases} \displaystyle {3 \over 2} \int_{0}^{\infty } {\log w \over (1 + w)^{5 / 2}} dw \text{,} & \text{if $k < k_{*}$} \text{,} \\ \displaystyle \Big( {p + 1 \over 2} + a - k_{*} + 1 \Big) \int_{0}^{\infty } {\log w \over (1 + w)^{(p + 1) / 2 + a - k_{*} + 2}} dw \text{,} & \text{if $k = k_{*}$} \text{,} \\ \displaystyle {1 \over 2} \int_{0}^{\infty } {\log w \over (1 + w)^{3 / 2}} dw \text{,} & \text{if $k > k_{*}$} \text{.} \end{cases} \non 
\end{align}
Note that 
\begin{align}
&{3 \over 2} \int_{0}^{\infty } {\log w \over (1 + w)^{3 / 2}} dw - {3 \over 2} \int_{0}^{\infty } {\log w \over (1 + w)^{5 / 2}} dw %
= {3 \over 2} \int_{0}^{\infty } {w \log w \over (1 + w)^{5 / 2}} dw %
= \int_{0}^{\infty } {\log w + 1 \over (1 + w)^{3 / 2}} dw \text{,} \non 
\end{align}
where the second equality follows from integration by parts, and that 
\begin{align}
&\Big( {p + 1 \over 2} + a - k_{*} + 1 \Big) \int_{0}^{\infty } {\log w \over (1 + w)^{(p + 1) / 2 + a - k_{*} + 2}} dw \non \\
&\ge \int_{0}^{\infty } {\log w \over (1 + w)^{5 / 2}} dw / \int_{0}^{\infty } {1 \over (1 + w)^{5 / 2}} dw %
= {3 \over 2} \int_{0}^{\infty } {\log w \over (1 + w)^{5 / 2}} dw \text{,} \non 
\end{align}
where the inequality follows from the covariance inequality. 
Then 
\begin{align}
&\lim_{\latso + \latsm / \latso + \latsu / \latsm \to 0} W_{k}^{\dag } ( {\blat {}^2} ) \begin{cases} \displaystyle = W - 2 \text{,} & \text{if $k < k_{*}$} \text{,} \\ \displaystyle \ge W - 2 \text{,} & \text{if $k = k_{*}$} \text{,} \\ \displaystyle = W \text{,} & \text{if $k > k_{*}$} \text{,} \end{cases} \label{eq:v3_2} 
\end{align}
where 
\begin{align}
W &= {1 \over 2} \int_{0}^{\infty } {\log w \over (1 + w)^{3 / 2}} dw = - \ga - \psi (1 / 2) = - \ga - (- \ga - 2 \log 2) = 2 \log 2 \label{eq:v3_3} 
\end{align}
by a result on page 544 of %
\cite{gr2014}. 

By (\ref{eq:v3_1}), (\ref{eq:v3_2}), and (\ref{eq:v3_3}) and by assumption, 
\begin{align}
&\lim_{\latso + \latsm / \latso + \latsu / \latsm \to 0} \lim_{\ep \to 0} {\pd \over \pd \ep } \Big[ \Big\{ \sum_{k = 1}^{p} {W_k ( {\blat {}^2}, \ep ) \over ( \lat _{k}^{2} )^{\ep }} \Big\} / \sum_{k = 1}^{p} {1 \over ( \lat _{k}^{2} )^{\ep }} \Big] %
\non \\
&\ge W - {2 k_{*} \over p} > W - 1 - {2 (a + 1) \over p} = 2 \log 2 - 1 - {2 (a + 1) \over p} > 0 \text{.} \label{eq:v3_4} 
\end{align}
By (\ref{eq:v3_0.5}), (\ref{eq:v3_0.75}), and (\ref{eq:v3_4}), for some ${\blat {}_{0}^{2}} = ( \lat _{0, k}^{2} )_{k = 1}^{p} \in (0, \infty )^p$, for some $0 < \ep _0 < 1 / 2$, for any $0 < \ep < \ep _0$, for any $0 < K < 1$, 
\begin{align}
1 &< \Big\{ \sum_{k = 1}^{p} {W_k ( {\blat {}_{0}^{2}} , \ep ) \over ( \lat _{0, k}^{2} )^{\ep }} \Big\} / \sum_{k = 1}^{p} {1 \over ( \lat _{0, k}^{2} )^{\ep }} \le {(P V_{2, \ep } ) (K {\blat {}_{0}^{2}} ) \over V_{2, \ep } (K {\blat {}_{0}^{2}})} \text{.} \non 
\end{align}
Suppose that the condition (\ref{eq:condition_bad}) holds for some $0 < \ep _1 < \ep _0$. 
Then there exist $0 < \de < 1$ and $\De > 0$ such that for any $0 < K < 1$, we have 
\begin{align}
&V_{2, \ep _1} (K {\blat {}_{0}^{2}} ) < (P V_{2, \ep _1} ) (K {\blat {}_{0}^{2}} ) \le \De + \de V_{2, \ep _1} (K {\blat {}_{0}^{2}} ) \text{.} \non 
\end{align}
Thus, 
\begin{align}
&V_{2, \ep _1} (K {\blat {}_{0}^{2}} ) < \De / (1 - \de ) \non 
\end{align}
for any $0 < K < 1$. 
This is a contradiction. 
\hfill$\Box$

\section{Proof of Theorem \ref{thm:general_ergodicity}}
\label{sec:proof_thm_second} 
Here, we use the drift and minorization technique (Jones and Hobert (2001)) to prove Theorem \ref{thm:lat}.

\subsection{Minorization conditions}
Fix $L > 0$. 
Suppose that $\tas , {\bet _{k}^{}}{}^2 / \sis , 1 / \sis , {\th _k}^2 / \sis \le L$ for all $k = 1, \dots , p$. 

Fix $k = 1, \dots , p$. 
Then we have 
\begin{align}
&p( \la _k | \bth , \bbet , \sis , \tas ) \non \\
&= \exp \Big[ - {\{ ( \tas )^{1 / 2} \la _k \bet _k - \th _k \} ^2 \over 2 d \sis } \Big] \pi _{\la } ( \la _k ) / \int_{- \infty }^{\infty } \exp \Big[ - {\{ ( \tas )^{1 / 2} \la \bet _k - \th _k \} ^2 \over 2 d \sis } \Big] \pi _{\la } ( \la ) d\la \non \\
&\ge \exp \Big[ - {\{ ( \tas )^{1 / 2} \la _k \bet _k - \th _k \} ^2 \over 2 d \sis } \Big] \pi _{\la } ( \la _k ) \ge \exp \Big[ - {\tas {\la _k}^2 {\bet _{k}^{}}{}^2 + {\th _k}^2 \over d \sis } \Big] \pi _{\la } ( \la _k ) \non \\
&\ge \exp \{ - (L / d) (1 + L {\la _k}^2 ) \} \pi _{\la } ( \la _k ) \text{.} \label{eq:general_minorization_1} 
\end{align}
Also, 
\begin{align}
&p( \tas | \bla , \bbet , \sis ) \non \\
&= {\rm{N}}_p (( \tas )^{1 / 2} ( \bla \circ \bbet ) | ( \X ^{\top } \X )^{- 1} \X ^{\top } \y , \sis ( \X ^{\top } \X )^{- 1} ) \pi _{\ta } ( \tas ) \non \\
&\quad / \int_{0}^{\infty } {\rm{N}}_p ( t^{1 / 2} ( \bla \circ \bbet ) | ( \X ^{\top } \X )^{- 1} \X ^{\top } \y , \sis ( \X ^{\top } \X )^{- 1} ) \pi _{\ta } (t) dt \non \\
&\ge \exp \Big[ - {1 \over 2 \sis } \{ ( \tas )^{1 / 2} ( \bla \circ \bbet ) - ( \X ^{\top } \X )^{- 1} \X ^{\top } \y \} ^{\top } ( \X ^{\top } \X ) \{ ( \tas )^{1 / 2} ( \bla \circ \bbet ) - ( \X ^{\top } \X )^{- 1} \X ^{\top } \y \} \Big] \pi _{\ta } ( \tas ) \non \\
&\ge \exp \Big[ - {1 \over \sis } \{ \tas ( \bla \circ \bbet )^{\top } ( \X ^{\top } \X ) ( \bla \circ \bbet ) + \y ^{\top } \P _{\X } \y \} \Big] \pi _{\ta } ( \tas ) \non 
\end{align}
and, since $\X ^{\top } \X \le \I ^{(p)} / d$, it follows that 
\begin{align}
p( \tas | \bla , \bbet , \sis ) &\ge \exp \Big\{ - {1 \over \sis } ( \tas \| \bbet \| ^2 \| \bla \| ^2 / d + \y ^{\top } \P _{\X } \y ) \Big\} \pi _{\ta } ( \tas ) \non \\
&\ge \exp (- L \| \y \| ^2 ) \exp (- \tas p L \| \bla \| ^2 / d) \pi _{\ta } ( \tas ) \non \\
&= \exp (- L \| \y \| ^2 ) \Big\{ \int_{0}^{\infty } \exp (- p L t \| \bla \| ^2 / d) \pi _{\ta } (t) dt \Big\} \non \\
&\quad \times \exp (- p L \tas \| \bla \| ^2 / d) \pi _{\ta } ( \tas ) / \int_{0}^{\infty } \exp (- p L t \| \bla \| ^2 / d) \pi _{\ta } (t) dt \text{.} \label{eq:general_minorization_2} 
\end{align}

By (\ref{eq:general_minorization_1}), a minorization condition holds for Algorithm \ref{algo:pc-two}. 
By (\ref{eq:general_minorization_1}) and (\ref{eq:general_minorization_2}), 
\begin{align}
&p( \tas _{\rm{new}} | \bla _{\rm{new}} , \bbet , \sis ) p( \bla _{\rm{new}} | \bth , \bbet , \sis , \tas ) \non \\
&\ge \Big\{ \exp (- p L \tas _{\rm{new}} \| \bla _{\rm{new}} \| ^2 / d) \pi _{\ta } ( \tas _{\rm{new}} ) / \int_{0}^{\infty } \exp (- p L t \| \bla _{\rm{new}} \| ^2 / d) \pi _{\ta } (t) dt \Big\} \non \\
&\quad \times \exp (- L \| \y \| ^2 ) \Big\{ \int_{0}^{\infty } \exp (- p L t \| \bla _{\rm{new}} \| ^2 / d) \pi _{\ta } (t) dt \Big\} \exp \{ - (L / d) (p + L \| \bla _{\rm{new}} \| ^2 ) \} p( \bla _{\rm{new}} ) \non 
\end{align}
for all $\tas _{\rm{new}} \in (0, \infty )$ and $\bla _{\rm{new}} \in \mathbb{R} ^p$. 
Thus, a minorization condition holds for Algorithm \ref{algo:pc-three}.

\subsection{Drift conditions}
For all $k = 1, \dots , p$, 
\begin{align}
E[ {\th _k}^2 | \bbet , \sis , \tas , \bla ] &= \sis \A ^{k, k} + [( \e _{k}^{(p)} )^{\top } \A ^{- 1} \{ ( \tas )^{1 / 2} ( \bla \circ \bbet ) / d + \b \} ]^2 \non \\
&\le \sis \A ^{k, k} + \A ^{k, k} \{ ( \tas )^{1 / 2} ( \bla \circ \bbet ) / d + \b \} ^{\top } \A ^{- 1} \{ ( \tas )^{1 / 2} ( \bla \circ \bbet ) / d + \b \} \non \\
&\le \sis \A ^{k, k} + \A ^{k, k} 2 \{ \b ^{\top } \A ^{- 1} \b + \tas ( \bla \circ \bbet )^{\top } \A ^{- 1} ( \bla \circ \bbet ) / d^2 \} \non \\
&\le d \sis + 2 d \{ d \| \b \| ^2 + \tas ( \bla \circ \bbet )^{\top } \I ^{(p)} ( \bla \circ \bbet ) / d \} \text{.} \non 
\end{align}
Therefore, 
\begin{align}
E[ {\th _k}^2 / \sis | \tas , \bla ] &\le d + 2 d^2 \| \b \| ^2 E[ 1 / \sis | \tas , \bla ] + 2 \tas E[ \| \bla \circ \bbet \| ^2 / \sis | \tas , \bla ] \non 
\end{align}
for all $k = 1, \dots , p$. 
We now show that the right-hand side is bounded. 
We have 
\begin{align}
E[ \| \bla \circ \bbet \| ^2 | \sis , \tas , \bla ] &= \sum_{k = 1}^{p} {\la _k}^2 E[ {\bet _{k}^{}}{}^2 | \sis , \tas , \bla ] \non \\
&= \sum_{k = 1}^{p} {\la _k}^2 \{ \Var ( \bet _k | \sis , \tas , \bla ) + (E[ \bet _k | \sis , \tas , \bla ])^2 \} \text{.} \non 
\end{align}
Let $D > 0$ be such that $D \X ^{\top } \X \ge \I ^{(p)}$. 
Then 
\begin{align}
\Var ( \bet _k | \sis , \tas , \bla ) &\le ( \sis / \tas ) ( \X ^{\top } \X )^{k, k} / {\la _k}^2 \non 
\end{align}
for all $k = 1, \dots , p$ and 
\begin{align}
\sum_{k = 1}^{p} {\la _k}^2 (E[ \bet _k | \sis , \tas , \bla ])^2 &= \| \bLa \{ \I ^{(p)} + \tas \bLa ( \X ^{\top } \X ) \bLa \} ^{- 1} \bLa \X ^{\top } \y ( \tas )^{1 / 2} \| ^2 \non \\
&= \tas \| \{ ( \bLa ^{- 1} )^2 + \tas ( \X ^{\top } \X ) \} ^{- 1} \X ^{\top } \y \| ^2 \non \\
&\le D \tas [ \{ ( \bLa ^{- 1} )^2 + \tas ( \X ^{\top } \X ) \} ^{- 1} \X ^{\top } \y ]^{\top } ( \X ^{\top } \X ) [ \{ ( \bLa ^{- 1} )^2 + \tas ( \X ^{\top } \X ) \} ^{- 1} \X ^{\top } \y ] \non \\
&= D \tas \y ^{\top } \X \{ ( \bLa ^{- 1} )^2 ( \X ^{\top } \X )^{- 1} ( \bLa ^{- 1} )^2 + 2 \tas ( \bLa ^{- 1} )^2 + ( \tas )^2 ( \X ^{\top } \X ) \} ^{- 1} \X ^{\top } \y \text{.} \non 
\end{align}
Therefore, 
\begin{align}
E[ \| \bla \circ \bbet \| ^2 | \sis , \tas , \bla ] &\le ( \sis / \tas ) \tr \{ ( \X ^{\top } \X )^{- 1} \} + D \tas \y ^{\top } \X \{ ( \tas )^2 ( \X ^{\top } \X ) \} ^{- 1} \X ^{\top } \y \non 
\end{align}
and 
\begin{align}
\tas E[ \| \bla \circ \bbet \| ^2 / \sis | \tas , \bla ] &\le \tr \{ ( \X ^{\top } \X )^{- 1} \} + D \y ^{\top } \X( \X ^{\top } \X )^{- 1} \X ^{\top } \y E[ 1 / \sis | \tas , \bla ] \text{.} \non 
\end{align}
Also, 
\begin{align}
E[ 1 / \sis | \tas , \bla ] &\le (n / 2 + a' ) / b' \text{.} \non 
\end{align}

Next, fix $k = 1, \dots , p$. 
Then 
\begin{align}
\Var ( \bet _k | \sis , \tas , \bla ) &\le \sis \non 
\end{align}
and 
\begin{align}
(E[ \bet _k | \sis , \tas , \bla ])^2 &= [( \e _{k}^{(p)} )^{\top } \{ \I ^{(p)} + \tas \bLa ( \X ^{\top } \X ) \bLa \} ^{- 1} \bLa \X ^{\top } \y ( \tas )^{1 / 2} ]^2 \non \\
&\le \tas ( \e _{k}^{(p)} )^{\top } \{ \I ^{(p)} + \tas \bLa ( \X ^{\top } \X ) \bLa \} ^{- 1} \e _{k}^{(p)} \y ^{\top } \X \bLa \{ \I ^{(p)} + \tas \bLa ( \X ^{\top } \X ) \bLa \} ^{- 1} \bLa \X ^{\top } \y \non \\
&\le \tas \y ^{\top } \X \bLa \{ \tas \bLa ( \X ^{\top } \X ) \bLa \} ^{- 1} \bLa \X ^{\top } \y \non \\
&= \y ^{\top } \P _{\X } \y \text{.} \non 
\end{align}
Therefore, 
\begin{align}
E[ {\bet _{k}^{}}{}^2 / \sis | \tas , \bla ] &\le 1 + \y ^{\top } \P _{\X } \y E[ 1 / \sis | \tas , \bla ] \text{.} \non 
\end{align}

Finally, for all $\ep > 0$, we have 
\begin{align}
E[ ( \tas )^{\ep } | \bla ] &\le \int_{0}^{\infty } t^{\ep } {\pi _{\ta } (t) / | \I ^{(p)} + t \bLa ( \X ^{\top } \X ) \bLa |^{1 / 2} \over ( b' )^{n / 2 + a'}} dt / \int_{0}^{\infty } {\pi _{\ta } (t) / | \I ^{(p)} + t \bLa ( \X ^{\top } \X ) \bLa |^{1 / 2} \over ( \| \y \| ^2 / 2 + b' )^{n / 2 + a'}} dt \non \\
&\le {( \| \y \| ^2 / 2 + b' )^{n / 2 + a'} \over ( b' )^{n / 2 + a'}} \int_{0}^{\infty } t^{\ep } \pi _{\ta } (t) dt \non 
\end{align}
by the covariance inequality. 
If $\tas \le \overline{\tas }$ for some known $\overline{\tas } > 0$, then $E[ ( \tas )^{\ep } | \bla , \bbet , \sis ] \le ( \overline{\tas } )^{\ep }$.

\end{document}